 \documentclass{amsart}
\usepackage{amsmath}
\usepackage{amsxtra}
\usepackage{amscd}
\usepackage{amsthm}
\usepackage{amsfonts}
\usepackage{amssymb}
\usepackage{eucal}
\usepackage{url}
\usepackage{xcolor}

\usepackage{enumerate}
 
\def\eqlas{\equiv_{Las}}

\def\0{\tilde{0}}
\def\cir{}

\def\Xx{\mathcal{X}}
\def\R{\mathcal{R}}
\def\Tt{\mathcal{T}}
\def\Mm{\mathcal{M}}
\def\ff{\mathfrak{f}}
\def\bQ{\bar{Q}}
\def\ha{\hat{a}}   \def\hb{\hat{b}}

\newtheorem{thm}[subsection]{Theorem}
\newtheorem{cor}[subsection]{Corollary}
\newtheorem{lem}[subsection]{Lemma}

\newtheorem{prop}[subsection]{Proposition}

\newtheorem{defn}[subsection]{Definition}

\newtheorem{rem}[subsection]{Remark}
\newtheorem{rems}[subsection]{Remarks}

\def\Om{\Omega}
\newtheorem{question}[subsection]{Question}
\newtheorem{example}[subsection]{Example}

\newtheorem{examples}[subsection]{Examples}
\newtheorem{problem}[subsection]{Problem}

\theoremstyle{definition}
\theoremstyle{remark}

\newcommand{\thmref}[1]{Theorem~\ref{#1}}
\newcommand{\secref}[1]{\S~\ref{#1}}
\newcommand{\lemref}[1]{Lemma~\ref{#1}}
\newcommand{\defref}[1]{Definition~\ref{#1}}
\newcommand{\propref}[1]{Proposition~\ref{#1}}
\newcommand{\corref}[1]{Corollary~\ref{#1}}

\newcommand{\remref}[1]{Remark~\ref{#1}}
\newcommand{\exref}[1]{Example~\ref{#1}}

\newcommand{\nc}{\newcommand}
\nc{\renc}{\renewcommand}
\nc{\ssec}{\subsection}
\nc{\sssec}{\subsubsection}
\nc{\on}{\operatorname}

\nc\ol{\overline}
\nc\wt{\widetilde}
\nc\wh{\widehat}
\nc\tboxtimes{\wt{\boxtimes}}

\renc{\d}{{\delta}}
\nc{\Aa}{{\mathbb{A}}}
 \nc{\Gg}{{\mathbb{G}}}  
\nc{\Hh}{{\mathbb{H}}}
 \nc{\Nn}{{\mathbb{N}}}
\nc{\Pp}{{\mathbb{P}}}
\nc{\Rr}{{\mathbb{R}}}
\nc{\BV}{{\mathbb{V}}}
\nc{\BW}{{\mathbb{W}}}
\nc{\Zz}{{\mathbb{Z}}}
\nc{\Qq}{{\mathbb{Q}}}
\nc{\St}{{\mathcal{S}}}
\nc{\Cc}{{\mathbb{C}}}
\nc{\Ff}{{\mathbb{F}}}
\def\inft{\bigwedge}
 
\def\eqlc{\equiv_{lc}} 
 
\nc{\bk}{{\mathbf{k}}}
\nc{\bb}{{\mathbf{b}}}
\nc{\bc}{{\mathbf{c}}}
\nc{\bd}{{\mathbf{d}}}
\nc{\be}{{\mathbf{e}}}
\nc{\bj}{{\mathbf{j}}}
\nc{\bn}{{\mathbf{n}}}
\nc{\bp}{{\mathbf{p}}}
\nc{\bq}{{\mathbf{q}}}
\nc{\bF}{{\mathbf{F}}}
\nc{\bu}{{\mathbf{u}}}
\nc{\bv}{{\mathbf{v}}}
\nc{\bx}{{\mathbf{x}}}
\nc{\bs}{{\mathbf{s}}}
\nc{\by}{{\mathbf{y}}}
\nc{\bw}{{\mathbf{w}}}
\nc{\bA}{{\mathbf{A}}}
\nc{\bK}{{\mathbf{K}}}
\nc{\bI}{{\mathbf{I}}}
\nc{\bB}{{\mathbf{B}}}
\nc{\bG}{{\mathbf{G}}}
\nc{\bC}{{\mathbf{C}}}
\nc{\bD}{{\mathbf{D}}}
\nc{\bP}{{\mathbf{P}}}
\nc{\bH}{{\mathbf{H}}}

\nc{\bM}{{\mathbf{M}}}
\nc{\bN}{{\mathbf{N}}}
\nc{\bV}{{\mathbf{V}}}
\nc{\bU}{{\mathbf{U}}}
\nc{\bL}{{\mathbf{L}}}

\nc{\bT}{{\mathbf{T}}}
\nc{\btau}{{\bar{\tau}}}
\nc{\bW}{{\mathbf{W}}}
\nc{\bX}{{\mathbf{X}}}
\nc{\bY}{{\mathbf{Y}}}
\nc{\bZ}{{\mathbf{Z}}}
\nc{\bS}{{\mathbf{S}}}

\nc{\sA}{{\mathsf{A}}}
\nc{\sX}{{\mathsf{X}}}
\nc{\sU}{{\mathsf{U}}}
\nc{\sB}{{\mathsf{B}}}
\nc{\sC}{{\mathsf{C}}}
\nc{\sD}{{\mathsf{D}}}
\nc{\sF}{{\mathsf{F}}}
\nc{\sG}{{\mathsf{G}}}
\nc{\sK}{{\mathsf{K}}}
\nc{\sM}{{\mathsf{M}}}
\nc{\sO}{{\mathsf{O}}}
\nc{\sQ}{{\mathsf{Q}}}
\nc{\sP}{{\mathsf{P}}}
\nc{\sZ}{{\mathsf{Z}}}
\nc{\sfp}{{\mathsf{p}}}
\nc{\sr}{{\mathsf{r}}}
\nc{\sg}{{\mathsf{g}}}
\nc{\sff}{{\mathsf{f}}}
\nc{\sfb}{{\mathsf{b}}}
\nc{\sfc}{{\mathsf{c}}}
\nc{\sd}{{\ltimes}}
\nc{\symd}{{\triangle}}
 \nc{\tU}{\widetilde{U}}
 \nc{\tH}{\widetilde{H}}
  \nc{\tG}{\widetilde{G}}
   \nc{\tE}{\widetilde{E}}
\nc{\tZ}{{\tilde{Z}}}
\nc{\tx}{{\tilde{x}}} \nc{\ty}{{\tilde{y}}}
\nc{\tq}{{\tilde{q}}}

\nc{\tfP}{{\widetilde{\mathfrak{P}}}{}}
\nc{\tz}{{\tilde{\zeta}}}
\nc{\tmu}{{\tilde{\mu}}}

 \def\e{\epsilon}

    \def\tL{{\widetilde{L}}}

\def\e{\epsilon}
 \def\Lam{\Lambda}

  \nc{\st}{{\mathop{\operatorname{\rm st}}}}     
    \nc{\cb}{{\mathop{\operatorname{\rm CB}}}}     
  \nc{\Ob}{{\mathop{\operatorname{\rm Ob}}}}
  \nc{\Sym}{{\mathop{\operatorname{\rm Sym}}}}
   \nc{\Aut}{{\mathop{\operatorname{\rm Aut}}}}
 \nc{\Spec}{{\mathop{\operatorname{\rm Spec}}}}
  \nc{\spec}{{\mathop{\operatorname{\rm Spec}}}}
\nc{\Ker}{{\mathop{\operatorname{\rm Ker}}}}
 \nc{\dom}{{\mathop{\operatorname{\rm dom}}}}
\nc{\End}{{\mathop{\operatorname{\rm End}}}}
 \nc{\Hom}{\on{\Hom}}
 \nc{\GL}{{\mathop{\operatorname{\rm GL}}}}
 \nc{\Id}{{\mathop{\operatorname{\rm Id}}}}
 \nc{\rk}{{\mathop{\operatorname{\rm rk}}}} 
 \nc{\length}{{\mathop{\operatorname{\rm length}}}}
\nc{\supp}{{\mathop{\operatorname{\rm supp}}}}
\nc{\val}{{\rm val}}

\nc{\res}{{\mathop{\operatorname{\rm res}}}}
\def\ind#1#2{ {#1} {\downarrow} {#2} }  
\def\Ind#1#2#3{{#1} {\downarrow}_{#3} {#2} }

\def\tensor{{\otimes}}
\def\meet{\cap}
\def\union{\cup}

\def\si{\sigma}
\def\g{\gamma}
\def\G{\Gamma}

\def\<{\begin}
 \def\>{\end}

\def\m{\smallsetminus}

\nc{\seq}[1]{\stackrel{#1}{\sim}}

\def\inv{^{-1}}
\def\claim#1{{\medskip \noindent \bf Claim #1.\ }}

\def\beq#1{\begin{equation} \label{#1}  }
\def\eeq{\end{equation}}

\def\Uu{\mathbb U}

\def\prf{\begin{proof}}
\def\pv{\end{proof} }
 \def\eprf{\end{proof} }

 \def\tp{\mathop{\rm tp}\nolimits}
\def\acl{\mathop{\rm acl}\nolimits}
 \def\dcl{\mathop{\rm dcl}\nolimits}

\def\lbl#1{   \label{#1}  }
\def\a{\alpha}

\def\ba{\bar{a}}

 \renc{\b}{{\beta}}

\def\eqtn{\underset{<n}{=}}
\def\eqthree{{\underset{<3}{=}}}
\def\simg{\sim_G}

\def\simv{\sim_V}

\def\std#1{{\widehat{#1}}}

 \def\2norm#1{\lVert{#1}\rVert_2}

 \title{Approximate equivalence relations}  
 \author{Ehud Hrushovski  }
 \thanks{Research   funded in part by the European Research Council under the European Unions Seventh Framework Programme (FP7/2007- 2013)/ERC Grant Agreement No. 291111.}

\begin{document}

 \begin{abstract}
 Generalizing   results for approximate subgroups, we study approximate equivalence relations up to commensurability,  in the presence of a definable measure.

      As a basic framework, we give a presentation of probability logic based on continuous logic.   
Hoover's normal form is valid here; if one begins with a discrete logic structure, it reduces arbitrary formulas of probability logic to correlations between quantifier-free formulas.  We completely classify  binary  correlations  in terms of the Kim-Pillay space, leading to strong results on the interpretative power of   pure probability logic over a binary language.  Assuming  higher amalgamation of independent types, we prove a higher stationarity statement for such correlations.  
 
We also give a short model-theoretic proof of a categoricity theorem for continuous logic structures with a
measure of full support, generalizing theorems of Gromov-Vershik and Keisler, and 
often providing a canonical model for a complete pure probability logic theory.  These results also apply to local probability logic, providing in particular a canonical model for a local pure probability logic theory with a unique 1-type and geodesic metric.

   For sequences of  approximate equivalence relations  with an 'approximately unique' probability logic 1-type,  we obtain a structure theorem    generalizing
  the `Lie model' theorem for approximate subgroups, \thmref{1}.  The models here are Riemannian homogeneous spaces, fibered over a locally finite graph.

 Specializing to definable graphs over finite fields, we show that  after a finite partition, a definable binary relation converges  in finitely many  self-compositions to an equivalence relation of geometric origin.  This generalizes the main  lemma for strong   approximation of groups.
 
   For NIP theories, pursuing a question of Pillay's, we prove an archimedean 
   finite-dimensionality statement for  the automorphism groups of definable measures, acting on a given type of definable sets.   
 This can be seen as an archimedean analogue of results of Macpherson and Tent on NIP profinite groups.

   \end{abstract}

   \maketitle
   
\begin{section}{Introduction}

Let $G$ be a group, with a translation-invariant, finitely additive measure on some Boolean algebra of subsets of $G$.   For instance $G$ may be an ultraproduct of finite groups with their counting measures.
A symmetric subset $X$ of a group $G$ is called
a {\em near-subgroup}\footnote{The definition of near-subgroup in  \cite{nqf} is a little more general.}
   
  if  both $X$ and the   triple product set $XXX$ 
are measurable and of finite, nonzero volume.      
These are very closely connected to  amenable  {\em approximate subgroups}, and arise in many branches of mathematics; 
see e.g. \cite{breuillard-lectures} for an introduction.

If $X$ is a Lie group, such as $GL_n(\Rr)$, a compact neighborhood of the identity is a near-subgroup;
all such neighborhoods are commensurable, i.e. each is covered by finitely many translates of the other.  
Conversely it 
  was shown in \cite{nqf} that a near-subgroup determines canonically a   connected Lie group $L$,  
 so that $X$  is commensurable with a pullback of a  compact neighborhood of the identity  in $G$.  This was used in particular to give a proof of (a strengthening of) Gromov's polynomial growth theorem,   based on measure-theoretic rather than metric properties of such a group.

Gromov was at that time writing \cite{gromov-ergo}.    He wrote:  
I  think there are many  ``almost structures"  which are far from ``actual structures" and these may play an essential role in how
 the brain generates math (...)  the ``dictionary structure"  contains algebraic
pattern e.g. of ``categories",  ``multi-categories" and also of ``2-categories" but these
patterns  are ``not perfect"  , e.g. some compositions may be  
``undefined/not implemented" and some may be non-associative.  (\cite{gromov-email}).
 
This comment precipitated the current work.   
The simplest kind of category is a groupoid, and the simplest groupoids, other than groups, 
are equivalence relations.    Transposing  the
  result on statistical recognition of approximate subgroups to approximate equivalence relations thus
 appears as a natural first step.    An equivalence relation is a symmetric relation $R$  
satisfying $R \circ R = R$; an approximate equivalence relation is one where $R \circ R$ is "commensurable" to $R$, i.e. all the $R \circ R$-neighbors of any point $a$ are   $R$-neighbors of some finite number of 
points; 
see  Definition  \ref{near-def}.    In case $X$ is an approximate subgroup
of a group $G$, one can define $R(x,y) := x \inv y \in X$; then it is easy to see that  $R$ is an approximate equivalence relation.
  It will turn out in fact that using this construction, the results we will obtain on approximate equivalence relations imply   
  the ones on approximate groups; the latter appears as the special case of approximate equivalence relation with a transitive
  automorphism group.      
  
   The language of the paper is that of probability logic; we will give details of that separately below.

{\bf {Stabilizer Theorem.}}
 
The dimension-theoretic stabilizer  was first introduced to model theory by Boris Zilber, in the setting of groups of finite Morley rank.  It was transformative to the subject,  enabling for instance the proof
of Zilber's indecomposability theorem, that seemed previously to belong to geometry rather than model theory.  The construction was  generalized to stable theories; the  uniqueness of independent pairs, given their individuual types,   was key.  Later (\cite{pac},\cite{CH},\cite{KP2}) it was realized that a good theory of the stabilizer exists without this uniqueness, 
under 
the weaker condition of the {\em independence theorem}.  
This statement asserts the  existence, under suitable conditions,  of a $3$-type with three prescribed restrictions to independent $2$-types.  
 
Still later, it was possible to transpose these ideas to groups defined in any theory, 
carrying a suitable measure.     This was the basis of a connection between
near subgroups and  the locally compact world:   for near subgroups $X$, the stabilizer is a $\bigwedge$-definable
group  $S$  contained in $X^4$, such that $X/S$ is compact.  Equivalently, one can find definable sets $Y$ commensurable to $X$ allowing a prescribed number of multiplications staying within $X^4$.

 In section \ref{4-stabilizer} we will prove a generalization of the stabilizer theorem to near equivalence relations.  \thmref{stabilizer} provides a canonical   $\bigwedge$-definable equivalence relation $S$ contained in $R^{\circ 4}$,  such that each neighbor set of $R$ is compact modulo $S$.

{\bf Riemannian homogeneous spaces.}  In the case of groups, one can go further and
describe approximate subgroups up to commensurability as approximations to finite dimensional Lie groups.  In section \ref{riemannian}
we obtain a similar theorem under an assumption of {\em approximate homogeneity}.      
  The   model spaces  are now 
  Riemannian homogeneous spaces,   with a `mesoscopic' graph relation connecting two points at distance at most $1$.
These are fibered over locally finite graphs; we obtain only a partial description of the total space, but a full description 
 of each connected component.  See \thmref{1}.

{\bf{Pseudo-finite fields.}}
Section \ref{grtogr}, ``From groups to graphs",   concerns approximate equivalence relations definable  in pseudo-finite fields.   
 
  In \cite{HP2}, a model-theoretic proof was given of
the `strong approximation' theorem of  Gabber, Matthews-Vaserstein-Weisfeiler, Gabber and Nori
(see \cite{mvw}, \cite{nori}) on   subgroups of   algebraic groups over the $p$-element field 
$\Ff_p$ for large $p$.    The main model-theoretic ingredient was the study of generation of groups by definable sets.  This
is   generalized in  section \ref{grtogr}  to the generation of equivalence relations by definable relations.    An arbitrary definable relation is decomposed, in each piece of a partition,   into relations that  generate an equivalence relation in finitely many steps, and relations of finite valency.  
For pseudo-finite fields this decomposition has an explicit algebraic form; but the general result is proved in the setting of simple theories with a well-behaved finite dimension theory.

{\bf{The measure stabilizer in NIP theories}}  In section \ref{NIPgr} we take an alternate route to finite dimensionality, under an assumption of NIP.
 The passage from locally compact spaces to finite-dimensional Riemannian manifolds in  section \ref{riemannian}  involved factoring out a large compact normal subgroup; commensurability is preserved, but little control over this compact kernel is available.   
In particular for a near-subgroup $X$ of a compact group $G$, i.e. a definable subset of positive measure, this procedure
loses all information.  In section \ref{NIPgr}
we factor out only the measure-theoretic stabilizer of $X$, and prove, assuming NIP, that up to possible profinite parts, the result is a finite-dimensional Lie group.   

More generally,  
let $\mu$ be a definable measure on $X$, let $q(u)$ be a type and let $\phi(x,u)$ be a NIP formula.  $\phi$ establishes a relation between
the space $\sX$ of weakly random global types on $X$, and the space of Kim-Pillay strong types extending $q$.  We obtain corresponding
quotients $\sX_{\phi,q}$ of $\sX$, and a canonical space $\sU_{\mu,\phi}$ of strong types compatible with $q$.   The automorphism 
group of any saturated model induces a compact group $G_{\mu,\phi,q}$ acting faithfully on both $\sX_{\phi,q}$  and $\sU_{\mu,\phi}$.
We prove that $G=G_{\mu,\phi,q}$ has finite archimedean rank.  This means that $G$ has a minimal normal subgroup $G^{00}$ with 
$G/G^{00}$ profinite, a maximal (up to finite index) profinite normal subgroup $G_{00}$, and $G^{00} / G_{00}$ is a finite dimensional Lie group.

This result can be transposed from automorphism groups to  definable groups. Let $G$ now be a definable group in a NIP theory, and $\mu$ a translation invariant definable measure on $G$.      $G$ has a minimal
$\bigwedge$-definable subgroup $G^{00}$, and $K=G/G^{00}$ is compact in the logic topology.  These compact groups 
were  the subject of Pillay's conjectures in the o-minimal case, see \cite{pillay-conj}, \cite{HPP}, showing that 
$K$ is a Lie group of the same dimension as $G$.
In the general NIP case, beyond compactness, the constraints on $K$ are unclear.  However one can define a canonical quotient $K_{P}$ of $K$ associated with a given   definable subset $P$ of $G$,
at least when $G$ carries a definable measure $\mu$;   namely identify two weakly random types of $G$ if they include the same set of 
translates of $P$.  (If $P=P_b$ is defined only with a parameter $b$, identify $p$ with $p'$ if for any $b' \models tp(b)$ and any $g,g' \in G$,
$g'P_{b'} \in p$ iff $gP_b \in G$.)   It was also Pillay who had the intuition that $K_{P}$ may be finite dimensional.  We show indeed that
$K_P$ has finite archimedean rank (\thmref{anand-nip-groups}).  In fact this is what motivated the more general \thmref{nipsymmetry}.  
 
 In view of results of Macpherson and Tent, it seems possible that a similar   finite-dimensionality phenomenon is valid for the totally disconnected part of $G$,
and in fact $G$ is  of adelic origin; in the most optimistic scenario, $G$ is interpretable in the model theoretic sum  finitely many $p$-adic and real fields.
 A generalization to the setting of approximate groups (viewed as
 piecewise-definable groups), using  Gleason-Yamabe theory in place of Peter-Weyl, would also be interesting.     
 See \ref{open7}.   
  
{\bf Probability logic}

In probability logic, existential     and universal quantifiers are augmented by {\em probability quantifiers}.   If as in \cite{keisler} they are entirely replaced by such quantifiers, we refer to {\em pure probability logic}.
Events  are taken to be definable subsets of a sort $X$, declared to be a {\em stochastic sort}.    $E_x \phi(x,y)$  is a formula with free variables $y$,  giving 
at $y=b$ the probability of the event $e_{\phi,b}=\{x: \phi(x,b) \}$.   

Since $E_x\phi(x,y)$ takes real values, it is natural to allow  all formulas to take bounded real values and employ continuous logic.   When $\phi$ takes values other than $0$ and $1$,  
$E_x\phi(b)$ is understood as the {\em expectation} of $\phi(x,b)$.
The axioms, going back to Kolmogorov, S. Bernstein, R. von Mises, Hilbert and Bohlmann, are just finite additivity and positivity.  Countable additivity is not assumed but is automatically obtained for the induced measure on the  type spaces over a model.     

By iterating the expectation quantifiers, we obtain measures on type spaces in several variables too.
The action of the symmetric group $Sym(n)$ on the space of $n$ types  is not assumed to preserve the measure; when it does, we say that Fubini holds.  Note that Fubini's {\em theorem} relates to the  algebra generated by rectangles; the Fubini {\em property} goes beyond this to arbitrary binary relations.  We do not assume Fubini at the level of the definition, but many results have stronger versions if Fubini is assumed.

We give a simple model-theoretic proof of a uniqueness theorem for models where every open set has positive measure.  This generalizes a theorem of Keisler's on uniqueness of hyperfinite models, and theorems of Gromov and Vershik on invariants for measured metric spaces.   When it exists the full-support model provides a canonical, homogeneous model for pure probability logic theories, replacing for some purposes the use of saturated models for first-order theories; this will be used in \secref{riemannian}.

 An elementary submodel $M_0$ is a kind of pool where everything not impossible has already happened.  Finite measure, like compactness, constrains the breadth of possible phenomena from above, and together they lead to a well-undestood theory over $M_0$(higher de Finetti theorems,  
  higher-dimensional Szemeredi lemma).  We present a model-theoretic version in Appendix B, either over an elementary submodel (following Towsner) or assuming qualitative higher amalgamation of types.   
  But we also pose the question 
  of finding the essential structures governing higher independence and hidden within $M_0$.    
   In this we try to emulate 
 Shelah's definability theorems for stable theories;  
 definability of types {\em over a model} is easier, but it was really the recognition of $acl^{eq}(0)$ and the proof of definability over that that enabled a useful theory of independence.   We obtain a satisfactory result for $n=3$, using auxiliary stable structures piecewise interpreted in the theory; so that an expectation statement can be referred to the stable structure, in this case Hilbert spaces.

These results 
will actually be required in the somewhat more general setting of local logic, where 
a large-scale metric is given and only balls of finite radius are assumed to have finite measure.   
 
  Appendix A develops basic stability theory for invariant relations, i.e.  disjunctions of $\bigwedge$-definable relations;
 the specialization to $\bigwedge$-definable relations is used in sections \ref{prelim} and \ref{problogic}.    
 Appendix C illustrates the use of probability logic in the setting of mixing results on  groups over pseudo-finite fields. 
 
 Many  open problems are  described throughout the text. 

 \smallskip
 
{\bf{Related work}}

While I thought at first that this was new territory, I soon learned that  approximate equivalence relations, by other names, are already very well studied.  I  talked about \thmref{1} in Aner Shalev's meeting on {\em Groups and Words}, in June 2012. 
 
Immediately afterwards, Nati Linial pointed out to me the relation of \thmref{stabilizer} with the work of Lovasz and Szegedy \cite{lovasz-szegedy}, \cite{lovasz}
on {\em graphons}.   Indeed while the language is different and the assumptions are slightly different, I believe
that basic methods of graphons yield an alternative proof of  \thmref{stabilizer}.

The `pure' probability logic we use, and the 
notion of ultraproduct that we obtain from it, are also closely related to Razborov's flag algebras \cite{razborov}.

Many variations on probability logic appear in the model-theoretic literature, implicitly and  explicitly.   The main results of \cite{keisler} are formulated within infinitary logic, $L_{\omega_1,\omega}$; whereas for us the use of compactness is essential.  Our semantics is in fact identical to that of definable Keisler measures \cite{NIP1} (named after a different work of Keisler's.)  One of the variations reported on in \cite{keisler}, due to Hoover, is finitary, as well as the treatment in \cite{gt};
they adjoin $\{0,1\}$-valued predicates $P_{x}^{>\alpha} \phi(x,y)$, intended to indicate  that the event defined by  $\phi(x,b)$ has
probability greater than or possibly equal to $\alpha$.  However in such a setting compactness
would dictate 
 
the existence of an event with probability $>0$ but $<1/n$ for each $n$; this is an intermediate state between emptyness and probability zero.
Such ghost predicates are difficult to control, and frequently lead to undecidability due to measure-zero phenomena that are not really intrinsic to the probabilistic viewpoint. The use of continuous logic is thus natural, being compatible with compactness and the standard interpretation of real-valued probability at the same time. 
 
   The papers \cite{JT}, \cite{marimon1}, \cite{marimon2}, \cite{ibarlucia}, \cite{BJM}
(in the setting of  $\aleph_0$-categoricity) are   closely related to our    independence theorem for probability logic 
\thmref{ind-prob}.  In particular \cite{JT}, Theorem  1.1 or 3.4 can be seen as special cases of \thmref{ind-prob} (1);  the measure-preserving action on $X$ assumed there can be viewed as data for definability of a measure on a new stochastic sort $X$.

In \cite{marimon}, \cite{jahel}, \cite{BJM}, the very interesting examples of two-graphs and kay-graphs are analyzed, showing in particular that a native
generalization of measure independence to higher amalgamation cannot hold, and on the other hand  (see \cite{marimon} 7.2.1, 7.2.3.)
that it does hold in certain circumstances, related to B.8.

 Ibarlucia in \cite{ibarlucia} employs a method of using auxiliary piecewise-definable stable structure, developed independently but very similar to ours; see also \cite{chevalier-H}.

     In the asymptotically finite setting, a statement equivalent to the stabilizer theorem for groups was
     independently proved in  \cite{sanders}, using a beautiful combinatorial argument. (See also \cite{bgt}.)   It is not clear if this method applies to finite approximate equivalence relations too.
 
As far as I know, the nearest result to the \thmref{1} on approximately homogeneous approximate equivalence relations is 
the paper \cite{bft} of Benjamini, Finucane and Tessera.  Their main  focus   is on finite approximate equivalence relations that are {\em exactly}
homogeneous for a group action; for these, they obtain results of the same strength as \cite{bgt}, in particular
showing that the phenomenon exists essentially for nilpotent groups and their homogeneous spaces.  But in the one-dimensional case, they also consider approximately homogenous  relations; again with technically somewhat different definitions than we use here.    

  A recent remarkable work of 
Gowers and Long \cite{gowers-long}  offers a wider interpretation of  `almost structure'.  

{\bf{Acknowledgements}}

This text is largely based on notes for a course given in Paris in Fall 2015, under the auspices of the 
Fondation Sciences Math\'ematique de Paris, and retains an expository aspect in the first sections.  I am very grateful to the participants and to the FSMP for their support.   
 Thanks also to Emmanuel Breuillard and Alex Gamburd  for many interesting conversations on the subjects of this paper, in particular at the time of the Groups and Words meeting when   Gamburd was working on \cite{BGS}.    
  \S 6 was motivated by these, as a kind of soft background result.  Thanks
to Alexis Chevalier and Mira Tartarotti for their  comments on section 3,     Pierre Simon on a previous version of Appendix A, and
to Shahar Mozes and Benjy Weiss for conversations around Keisler's relational law of large numbers.  
   And finally thanks to Dan Drai for many conversations that brought out, in particular, the relation of the stabilizer theorem or the graphon metric to measures used in practice in natural language processing.

The referees of this paper deserves exceptional credit.   In addition to many  local comments and references that  greatly contributed to the paper,  a serious issue was uncovered in  
 both \secref{riemannian} and \secref{grtogr} (and in a slightly different form in \S 7).  
 In both cases, the structure of  approximate equivalence relation is precisely described {\em piecewise}, but only a partial account is given of the
 interaction between the pieces.   See in particular Problem \ref{mixed} and Example \ref{unbounded-expansion}.

\end{section}

\begin{section}{Preliminaries }
\label{prelim}
\ssec{Real valued continuous logic.}  
Continuous model theory dates back to the book \cite{chang-keisler-c};  other roots lie in 
Robinson's nonstandard analysis of the same period, and their development (notably Henson's
study of nonstandard hulls of Banach spaces) in the  1970's.  Krivine studied a real-valued logic, and stability was introduced
to the area in \cite{krivinemaurey}.   A modern version, with a full-fledged stability theory as well as simplicity and NIP, was introduced by Ben Yaacov and coworkers in  \cite{byu}, \cite{benynip}, and other articles; see especially 
\cite{bybhu}.   

In continuous logic,  terms are defined in the same way as in first order logic.
But formulas $\phi$ are taken to take truth values in some compact   Hausdorff space $X_\phi$.  Any continuous
map on $c: X_\phi \times X_\psi$ to a compact space $ Y$ can be viewed as a connective, thus creating a
new formula $c(\phi,\psi)$ taking values in $Y$.  Any continuous map $q$  from the (Hausdorff) space of nonempty closed subsets of $X$
to a compact space $Y$ induces a quantifier, taking formulas $\phi(x,y)$  with   range $X$
to formulas $(qx)\phi(x,y)$ with free variables $y$ and range $Y$.  The interpretation of $(qx) \phi(x,b)$ in a structure
$A$ is $q(cl \{\phi(a,b): a \in A\})$.  (Chang and Keisler, 1966).

Specifically in continuous {\em  real-valued} logic,
the spaces $X_{\phi}$ will be closed intervals in $\Rr$, or occasionally in $\Rr \union \{\infty\}$.   The connectives
can be restricted to $+,\cdot,1$ (Stone-Weierstrass).     The quantifiers can be restricted to $\min$ and $\sup$,
though it is not always best to follow this religiously. 
We view two languages as having the same expressive power if a formula of one can be uniformly approximated
by a formula of the other, and vice versa.  (We do not seek formula-to-formula equality.)

A {\em complete theory} is a specified value for each sentence (formula with no free variables.)  
Similarly if $M$ is a structure, a {\em type $p(x)$ over $M$} is a specified value for each formula $\phi(x)$
with parameters from $M$.    

Let $A$ be a substructure of $M$.  
Formulas with parameters from $A$, and variable $x$, define functions on the set $S_x(A)$ of types $p(x)$ over
$A$; we topologize $S_x(A)$ minimally so that they are all continuous, in other words as a closed subset of the
  product topology.  Then $S_x(A)$ is compact, and any continuous function on $S_x(A)$  is uniformly approximated by formulas.     If $Z$ is a closed (respectively open) subset of   $S_x(A)$, we call $\{m \in M: tp(m/A) \in Z \}$ a 
   $\bigwedge$-definable (respectively $\bigvee$-definable) subset of $M$; these   
     notions should be used only in sufficiently saturated models, say ones where every type over $A$ is realized.

If $u$ is an ultrafilter on a set $I$ and $(a_i: i \in I)$ is an $I$-indexed family  of elements of a compact Hausdorff space $X$, there is always a unique $x \in X$ such that any neighborhood of $x$ contains almost all $a_i$ (according to $u$.)
This is denoted $\lim_{i \to u} a_i$.  

An {\em ultraproduct } along $u$ of structures $A_i$ for a real-valued language $L$ is defined in the usual way, except 
that the   value of a (basic) relation is the {\em limit along $u$} of the   values on the coordinates.  

\ssec{Metrics} \label{metrics}
There is often a distinguished binary formula $\rho$, whose interpretation is a metric $\rho: A^2 \to \Rr$,
and such that every term and every basic formula are uniformly continuous, by a prescribed modulus of continuity.
In this case, one modifies the definition of the  ultraproduct by identifying elements at distance zero.  This is  analogous to the situation with equality in 2-valued logic.  

\ssec{Localities}  \label{locality}
Note that the rules of continuous logic would force $\rho$ to have bounded image;  indeed for the   discussion 
so far, there is no harm
in replacing $\rho$ by $\min(\rho,1)$.     Ben Yaacov \cite{benyaacov-unbounded} defines an  unbounded version, where a fixed  unary function (the {\em gauge})   controls locality; it is similar to   many-sorted logic where quantifiers are restricted to finitely many sorts; but in place of a discrete set of sorts one has a continuous family.   We will however be interested in homogeneous structures, with a single $1$-type, and they are not compatible with   unary functions.  We will thus use a binary function $\rho^*$, satisfying the laws of a metric; it could be the same as the metric $\rho$, or distinct from it; in any case we are mostly interested in $\rho$ near $0$ (to determine a topology) and in $\rho^*$ near $\infty$
(to determine a coarse structure or, for us, the notion of locality, i.e. a family of sets where model-theoretic compactness will hold).

We allow relations $R$ to take unbounded values; but
 we assume that any basic relation comes not only
with a modulus of continuity (with respect to $\rho$) but also with a bound $b(R)$ on the support of $R$, so that 
\[ R(x_1,\ldots,x_n) \leq b( \max_{i,j} \rho^*(x_i,x_j)).\]   
Here $b$ is a continuous function $\Rr \to \Rr$ with compact support.
Similarly basic functions are assumed to take values at a bounded distance from 
their arguments.  
  We   redefine saturation by restricting to types at bounded distance.  
  
When we take an ultraproduct, 
we have to make an additional choice, beyond that of the ultrafilter.   Let $\Mm_0$ be the naive ultraproduct; then $\rho^*(x,y)$
can be infinite.  The relation:  $\rho^*(x,y)<\infty$ is an equivalence relation.   We make a choice of one class.    Thus
in an ultraproduct, $\rho^*(x,y)$ is finite {\em by definition}.    We refer to this as a local ultraproduct with locality relation $\rho^*$.    
 
We note that this logical structure depends on $\rho^*$ only up to coarse equivalence; replacing $\rho^*$ by $j \circ \rho^* $ where
$j: \Rr^+ \to \Rr^+$ is an order-preserving bijection will make no difference.

See \secref{local-structures} for more detail.

\begin{example} \rm \label{hilbert} The language of   Hilbert spaces is   taken to have sorts $S_r$ for any real $r \geq 0$, denoting the ball of radius $\leq r$.  
has terms $0,+,\cdot \alpha$ for any $\alpha \in \Cc$,  so $+: S_r \times S_{r'} \to S_{r+r'}$ and $\cdot \alpha: S_r \to S_{t}$ whenever $r|\alpha| \leq t$.    
There is one additional basic relation for the inner product, $(,): S_r \times S_{r'} \to [-rr',rr']$, the inner product; and the obvious axioms.  The metric is taken to be $|x-y|$, where $|x| = \sqrt{(x,x)}$.

The division into sorts adds somewhat artificial structure; a better approach is developed in \cite{benyaacov-unbounded}.
For the actual use of Hilbert spaces in this paper this will not be essential, we can take either one.

Suppose however we wish to consider $H$ as an affine space, without a distinguished $0$.  In this case it will not do to add sorts, whether discretely or continuously.  Instead we use the locality function $\rho^*(x,y)  = ||x-y||$; in this case it happens
to coincide with the metric.   The effect is again to limit quantifiers to bounded balls, but the balls can be anywhere on $H$.
\end{example}

\ssec{Cobounded equivalence relations and the logic topology} \label{logictop}  Let $X$ be a $\bigwedge$-definable set.  
A  $\bigwedge$-definable   relation $\Lam$  is called {\em a cobounded equivalence relation} if in any model  $M$, $\Lam$ 
defines an equivalence relation on $X(M)$, and $X(M)/\Lam$ has cardinality   bounded independently of $M$.    

We have $\Lam = \bigwedge_i \Lam_i$ with $\Lam_i$ definable, such that all antichains of $\Lam_i$ are finite, of size bounded by some $\beta_i \in \Nn$.   
If $M \models T$ and $N \succ M$, we can pick in $M$ a maximal antichain $c_1,\ldots,c_{\beta_i}$ for $\Lam_i$.   If
$a,b \in X(N) $ have the same type over $M$, or even just the same $\Lam_i$-type over $M$ for each $i$, then $(a,c_i) \in \Lam_i$ iff
$(b,c_i) \in \Lam_i$.  By maximality of the antichain, we do have $(a,c_i) \in \Lam_i$ for at least one $i \leq \beta_i$, and hence 
$(a,b) \in \Lam_i \circ \Lam_i$.  

 If $M \models T$ and $N \succ M$, and $a,b \in X(N) $ have the same type over $M$, then $a \Lam b$.
 Hence we have a natural map $S_X(M) \to X/\Lam$.  The image set is the same as $X(N) / \Lam$ for sufficiently saturated
 $N$, and we denote it simply by $X/\Lam$.   The surjective map $S_X(M) \to X/\Lam$  induces a topology on $X/\Lam$ - the {\em  logic topology} - which is compact and Hausdorff.   
 It does not depend on the choice of $M$.  Further, for any reduct $M'$ of $M$ such that $\Lam$ is still  $\bigwedge$-definable  in $M'$, since the compact $M$-induced topology refines the $M'$-induced Hausdorff topology, they coincide; expansions do
 not change the space $X/\Lam$.   
 
  For the same reason,  we have a well-defined map $S_X^{\Lam} (M) \to X/\Lam$ (where $S_X^{\Lam}$ is the space of $\Lambda$-types,
  meaning $\Lam_i$-types for each $i$);   and it induces the same topology on   $X/\Lam$.  Hence, if $Y$ is a closed subset of $X / \Lam$, then the pullback of $Y$ to $X$ is
  a $\bigwedge$-qf-definable set with parameters in $M$.

 Let $T$ be a complete theory, $X$ a sort.
There exists a unique {\em finest}  co-bounded  $\bigwedge$-definable equivalence relation $\Lam$ of $T$.  
  For this choice of $\Lam$, the space   $X/\Lam$   is called the space of compact Lascar types 
of $T$ in the sort $X$, or the  Kim-Pillay space $KP_T(X)$ (\cite{simple}).  Similarly for the qf KP-space.

The classes of $\Lam$  
 are called the compact Lascar types, or Kim-Pillay types.

\begin{rem}\rm \label{2top} Let $E$ be any co-bounded  $\bigwedge$-definable equivalence relation;
assume it is defined using formulas from a family $\Phi$, e.g. $xEy \iff \phi(x,y) = 0$ for each $\phi \in \Phi$.   Then
$X/E$ is a quotient of $KP_T(X)$, and we can compare the quotient topology to the topology defined as above using the space of $\Phi$-types alone.  
As the former is stronger and both are compact and Hausdorff, they must be equal.  
In other words, the image in $X/E$ of any $\bigwedge$-definable subset of $X$ is already the image of a  $\Phi-\bigwedge$-definable set.
\end{rem}

\ssec{Stability} 

\begin{defn}[\cite{krivinemaurey}]  A formula $\phi(x,y)$ is {\em stable} if for any model $M$ and any elements $a_i,b_i (i \in \Nn)$
of $M$, if $\lim_{i \to \infty} \lim_{j \to \infty}\phi (a_i,b_j)$ exists and equals $ \alpha $ and  $\lim_{j \to \infty} \lim_{i \to \infty} \phi(a_i,b_j)=\beta$,
then $\alpha=\beta$.  \end{defn}

The class of stable formulas $\phi(x,y)$ is easily seen to be  closed under continuous connectives.  

\begin{lem} \label{l2stable} Let $H$ be a Hilbert space,with elements $a_i,b_i (i \in \Nn)$ of the unit disk.  If $\lim_{i \to \infty} \lim_{j \to \infty} (a_i,b_j)$ exists then so does  $\lim_{j \to \infty} \lim_{i \to \infty} (a_i,b_j)$, and they are equal.  \end{lem}

This lemma means that $( , )$ is a  stable formula.    Using quantifier elimination for Hilbert spaces, this easily implies
that every formula is stable; but we will need this particular one.  The significance of this was realized in \cite{krivinemaurey} but also in  
\cite{grothendieck52}; see \cite{benyaacov-g}.   In the context of expectation quantifiers $E_t$ that will soon be introduced, it implies that
$E_t(f(x,t) g(x,t))$ is always a stable formula.
 
The following   is a continuous logic  version of Shelah's finite equivalence theorem (uniqueness of non-forking extensions over algebraically closed sets); see \cite{byu}.      The continuity in each variable is the open mapping theorem, 
or the definability of types.  Joint continuity does not hold in general.

(1) asserts that any value of $\phi(x,b)$ other than $\alpha$ causes forking, while (2) is a strong
converse asserting that the value $\alpha$ can be taken simultaneously for any family of $b$'s.

 As usual we write $p$ to denote the solution set of $p$; and $\alpha(a,b)$ for $\alpha(a/E, b/E')$.

\begin{thm} \label{uniqueness} \cite{byu} Let $\phi(x,y)$ be a stable formula on $P \times Q$.  
Then there exist  co-compact $\bigwedge$-definable equivalence relations $E$ on
 $P$
and $E'$  on $Q$,   
and a Borel  
 function $\alpha:  P/E \times Q/E' \to \Rr$, continuous in each variable and automorphism invariant, such that in any sufficiently saturated model $M$:
 \begin{enumerate}
\item In any prescribed $E'$-class  there exists a sequence $(b_j: j \in \Nn)$, 
such that for all $a \in P$, $\lim_{j \to \infty} \phi(a,b_j) - \alpha(a,b_j)=0$.  Equivalently, 
  for any $\e>0$, for some $k \in \Nn$, for any $J \subset \Nn$ with $|J| \geq k$,   and any $a \models p$, for some $j \in J$,
  $ | \phi(a,b_j) - \alpha(a,b_j)| < \e $
\item 
For any finite set $\{b_j: j \in J\} \subset Q$  and any $\e>0$,
there exists $a \in P$ in any prescribed $E$-class with $|\phi(a,b_j) - \a(tp(a),tp(b_j)| < \e$ for each $j$.
\item Let $M_0 \prec M$; assume $\phi$ is quantifier-free; then $\alpha(a,b)$ depends only on $qftp(a/M_0)$ and $qftp(b/M_0)$.\end{enumerate}
\end{thm}
 Thus $\alpha(p,q)$ gives the generic or expected value for $\phi(a,b)$ when $a \models p, b \models q$;
and any deviation from this value will cause dividing.  A more general statement will be proved in Appendix A.
(see \thmref{stable-i}.)

For each complete type $p \subset P$ 
 
  $\alpha$ is continuous as a function of two variables on $p \times Q$.   
   But in general it is not continuous on $P \times Q$,   even for the theory of pure equality augmented with infinitely many constants.  To see how this may arise in a probabilistic setting, consider the random graph, with infinitely many distinguished constants, and with a measure giving independent probability 1/2 to an edge.  Then for two types $p(x),q(y)$, and for $\phi(x,y)$ the probability that $z$ is a neighbor of both $x$ and $y$, we have $\alpha(p,q)=1/4$ {\em unless} $p  \vdash x=c_n$ and $q \vdash y=c_n$ for the same $n$; this is not bi-continuous.

     \ssec{Topologies}  \label{grouptop} \label{completion}

We discuss here some  elementary topology that will be needed later.

By a {\em pseudo-metric} on $X$ we mean a function $d: X^2 \to \Rr$ with $d(x,y)=d(y,x) \geq 0$,
$d(x,x)=0$, and $d(x,z) \leq d(x,y)+d(y,z)$.  There is a canonical map $j: X \to \bar{X}$ into a complete
metric space, preserving $d$, with dense image; we refer to $ \bar{X}$ as the 
 {\em completion} of $(X,d)$.     Typically $j$ is not injective.
 
Let  $f: X^n \to \Rr$ be a  function, uniformly continuous with respect to $d$; i.e. for all $\e>0$ there exists $\d >0$
such that  $|f(x)-f(y)| <\e$ whenever $x=(x_1,\ldots,x_n),y=(y_1,\ldots,y_n)  \in X^n$ and $d(x_i,y_i) < \d$.   Then $f$ induces
a function $\bar{f}: \bar{X}^n \to \Rr$.    Hence if $X$ carries an $L$-structure for some continuous logic language $L$,  an $L$-structure on  $\bar{X}$ is canonically induced.
All  continuous logic formulas are preserved, i.e $\phi(jx_1,\ldots,jx_n) ^{\bar{X}}= \phi(x_1,\ldots,x_n)^X$.   In particular,  the axioms for expectation quantifiers (\secref{problogic} (1-4))  are preserved. Hence if $X$ is a stochastic sort, then $\bar{X}$ becomes one too.  

Assume $X$ is a stochastic sort, with expectation quantifers $E$.  Write  $\Lam(x,y)$ if $j(x)=j(y)$; then $\Lam$ is a $\bigwedge$-definable.  Assume further as in \ref{logictop}
that $\Lam$ is co-bounded, so that $\bar{X}$ is compact.   Then the expectation quantifiers on $X$ induce a Borel measure on $\bar{X}$,
equivalently a positive linear functional $\int$ on the Banach space of real-valued continuous functions on $\bar{X}$.   Simply define 
$\int f = E (f \circ j)$.  We use here the fact that $f \circ j$ is a (uniform limit of) parameterically definable functions; this in turn can be seen using 
Stone-Weierstrass, as the definable functions into $\Rr$ form an algebra and separate points on $\bar{X}$.     
 
 \ssec{}  Let $d$ be a pseudo-metric  given by a formula of bounded real-valued continuous logic.   Assume 
  the equivalence relation $E$ defined by $d(x,y)=0$ is co-bounded.   Then $d$ induces
    a metric on $X/E$, interpreted in any sufficiently saturated model.  The metric is by definition continuous
    with respect to the logic topology on $X/E$.  The latter being compact, it follows that $d$ induces the logic topology;
    moreover,  $X/E$ is complete and hence coincides, as a metric space,  with the  completion $\bar{X}$.

This is valid locally in local continuous logic, with a metric $\rho$, when $d$ is definable  
and hence subordinate to $\rho$.   Namely, for any fixed $a \in X$, let $B$ be a ball of some finite radius $r$ around $X$.   Then as above, the logic topology on $B/E$ coincides with the 
topology induced by $d$.    It follows that globally,   the logic topology on $X/E$ is locally compact, and   induced by $d$.     Returning to $B$, by the same argument, the metric topology on $B/E$
induced by $d$ also coincides with the logic topology obtained by considering only subsets of $B$ defined with parameters in $B$.  We will use this remark (*) later on.  

\ssec{Isometry groups} \label{isometry_groups}  Let $Y$ be a locally compact metric space.   The   isometry group $G=Aut(Y,d)$ 
  is topologized by  the compact-open topology, or uniform convergence on compacts. On $Y$, the topology agrees with the topology of  pointwise convergence.   This is because
a compact $C$ admits a finite set $D_\e$ that is $\e$-dense in $C$; so if $f,g$ are isometries and $d(f(x),g(x)) < \e$ on for $x \in D_\e$,
then $d(f(x),g(x)) < 3 \e $ for all $x \in C$.  
  It is clear that left and right translations are continuous.  Thus to check continuity of inversion and multiplication, it suffices
  to verify it at the identity element.  
  If $g_i \to Id_Y$, then for any $a$ we have $g_i(a) \to a$, i.e. $d(g_i(a),a) \to 0$; since $g_i$ is an isometry, $d(a,g_i \inv(a)) \to 0$
  so inversion is continuous.   If also   $h_i \to 1$, then for any $a$ we have $h_i(a) \to a$; by local compactness we may take
  all $h_i(a)$ in some compact neighbourhood $C$ of $a$; since $g_i$ approaches $Id$ uniformly on $C$, for any $\e>0$,
  for large enough $i$ we have $d(g_i(y)), y) < \e$ for all $y \in C$ (for large enough $i$);  in particular  $d(g_ih_i(a), h_i(a)) < \e$;
  so $d(g_ih_i(a),a) < 2 \e$ for large $i$.  It follows that $G$ is a topological group.  
  
  The action $G \times Y \to Y$ is also easily seen to be continuous
    Moreover for $x_0 \in Y$,
  the map $G \to Y$, $g \mapsto gx_0$ is closed, since it suffices to check this after restricting
  to a closed bounded subset $Y'$ of $Y$, and the pre-image of $Y'$ is compact.  It follows that
  the stabilizer $G_{x_0}=\{g: gx_0=x_0\}$ is a closed subgroup, $Gx_0$ is closed in $Y$,
  and $G/G_{x_0}$ is homeomorphic to $Gx_0$

 \ssec{Graphs and metrics}
    A binary relation $R$ on $X$, viewed as a graph, is {\em connected} if for any $x ,y \in X$ there exist $n \geq 1$ and $x_1,\ldots,x_n \in X$ with
   $x=x_1$, $x_n=y$, and  such that    $R(x_i,x_{i+1})$  or $R(x_{i+1},x_i)$ hold for $ i<n$.   In this case, for the least such $n$, we define $d_R(x,y)=n-1$;
   $d_R$ will be referred to as the  metric associated to $R$.
    
 In the lemma below   two metrics appear, but all topological terms refer to $(X,d)$.
 
 \begin{lem}\label{grouptoplem}  Let $(X,d)$ be a locally compact metric space.   Let $\mathsf{R} \subset X^2$ be a closed binary relation, with $(X,\mathsf{R})$ connected; and let
 $d_\mathsf{R}$ be the associated   metric on $X$.   Assume some $d$-ball is contained in a $d_\mathsf{R}$-ball, and every $d_\mathsf{R}$-ball has compact closure. 
  
 Let $G=Aut(X,d,\mathsf{R})$ be the group of isometries of $(X,d)$ preserving $\mathsf{R}$.  
 Then $G$ is a locally compact topological group.  For any $x \in X$ and any compact $U \subset X$, $\{g \in G: g(x) \in U\}$ is compact.
 \end{lem}
  
\prf   
We saw above that    $G$ is a topological group, and that the topology given above agrees with  pointwise convergence; we will use the latter description in order to  reduce to compactness of product spaces.    
   
 To show that $G$ is locally compact, it suffices to show that  $1_G$ has a compact neighborhood.   Let $b_0$ be an open $d$-ball around $a_0$ contained  in a $d_\mathsf{R}$ ball.  Let $D_n$ be the $d_\mathsf{R}$-ball around $a_0$ of radius $n$; then $b_0 \subset D_{n_0}$ for some $n_0$.  
   Let $U_0 = \{g \in G:  g(a_0) \in b_0 \}$.  Then $U_0$ is an open neighborhood of $1$ in $G$.    If $g \in U_0$, then $g(D_n) \subset D_{n+n_0}$;
   indeed   $g(a_0) \in b_0$ so $d_\mathsf{R}(a_0,g(a_0)) \leq n_0$;  if $d_\mathsf{R}(x,a_0) \leq n$, then $d_\mathsf{R}(g(x),g(a_0)) \leq n$ (it is here that
   we use the assumption that $\mathsf{R}$, hence $d_\mathsf{R}$, are preserved by $G$.)  So $d_\mathsf{R}(a_0, g(x)) \leq n+n_0$.  
   Hence $U_0 \subset U_1$ where $U_1$ is the set of isometries of $X$  satisfying $g(D_n) \subset D_{n+n_0}$ for every $n$.
   Now $U_1$ is compact since it embeds homeomorphically into a closed suset of the product space $\Pi_n D_{n+n_0}^{D_n}$,
   mapping $g \to (g|D_n)_n$.    And $G$ is clearly a closed subgroup of the isometry group; so $G \meet U_1$ is compact too.       
      
\eprf 

\ssec{NIP}  Let $R(x,y)$ be a formula, and let $S_x(A)=S_x^R(A)$ 
be the space 
of quantifier-free $R$-types in the variable $x$ over a set $A$.   For  $R(x,y)$ taking values in $\{0,1\}$,
$R$  has NIP (=does not have the independence property)
if $|S_x(A)|$ grows at most polynomially  in $|A|$, i.e. $|S_x(A)| \leq C \cdot |A|^k$ for some $C,k$ and for all $A$.
By a theorem of Sauer, Shelah, and Vapnik-Cervonenkis,   this is equivalent to:  $|S_x(A)| < 2^{|A|}$ for $|A| > m$;
the minimal such $m$ is the Vapnik-Cervonenkis dimension, and we have $|S_x^R(A)| \leq  2 |A|^m$. 

   For an $\Rr$- valued formula $R(x,y)$, one says that $R$ has NIP if for any fixed $\e>0$, 
     $S_x^R(A)$ grows at most polynomially  in $|A|$ up to $\e$-resolution; in other words $S_x^R(A)$ can be covered
     by polynomially many $\e$-balls, or equivalently admits at most polynomially many (in $|A|$) pairwise disjoint $\e$-balls.
  See \cite{benyaacov-c-nip}.  This notion generalizes to general continuous logic (with values in compact
  spaces.) 
  
   On the other hand, we will say that $R$ has pNIP (of degree $d$)  if for $|A|,|n| \geq d$, 
    $S_x^R(A)$ can be covered
     by at most  $ (|A| n)^d$  $1/n$-balls.   Thus the growth is polynomial not only in the base size but also in the resolution.
     pNIP relations are closed under connectives corresponding to Lipschitz functions $\Rr^k \to \Rr$.     
 
  \begin{rem}  If the relation $\phi(x; yz) = R(x,y)\leq z$ has NIP, then $R$ has pNIP.
  Indeed the $R(x;y)$ types, up to $1/m$-resolution, over a set $b_1,\ldots,b_n$ can be viewed 
  as $\phi$-types over $b_1,\ldots,b_n, 0,1/m,\ldots,1$ so their number is polynomial in $mn$.     
  \end{rem}
       
 We will later need an  effective version of the {\em uniform law of large numbers} of Vapnik-Cervonenkis.  What is essential for us, to obtain finite packing dimension,  is that
 $N$ in \propref{evc} should grow at most polynomially with $n$.  This already follows from  \cite{vc}, Theorem 2.   (In the notation there, set $\e_0 = 1/(2n)$; one looks for a lower bound on $l$ that ensures that the right hand side is $<1$; this ensures not only existence but a nonzero percentage of    $N$-tuples $c_1,\ldots,c_N$, such that for each $b$ there are at least $1/(2n)$ values of $i$ with $R(c_i,b)$.)   However we quote more precise bounds.

 \begin{prop} [\cite{haussler-welzl}] \label{evc}    Let $(U,\mu)$ be a probability space, and  $R(u,b)$ be a relation  with Vapnik-Cervonenkis dimension $d$
 (i.e. the family of events $\{u: R(u,b)\}$ has  Vapnik-Cervonenkis dimension $d$.)  
  Then for any $n$ there exist $c_1,\ldots,c_N$, $N= 8dn \log(8dn) \leq 24(dn)^2$, 
 such that for any $b$ with $\mu R(u,b) \geq 1/n$, for some $i$ we have $R(c_i,b)$.  \end{prop}
 \prf  This follows from   Cor. 3.8 of \cite{haussler-welzl}, with $\delta=1/2$.    In \cite{haussler-welzl} the result is stated for $\mu$
a normalized counting measure on a large finite set, but as the bound does not mention the size of this set, the result immediately extends to all  probability spaces by a standard approximation argument, see \cite{hps} around 2.7.  
\eprf
             
  \begin{rem}\label{evc2}    
  Let $\mu$ be a   measure on $\phi$-types, not necessarily generically stable.  
  It is shown in \cite{NIP1} (Lemma 4.8(i)) 
  that for any $n$ and any model $M$ there exist $c_1,\ldots,c_N$ in an elementary extension $M^*$ of $M$, such that 
  for any $b$ with $\mu R(u,b) \geq 1/n$, for some $i$ we have $R(c_i,b)$.   The proof yields    a polynomial bound $N \leq O(n^\delta)$ for some $\delta$.
 
  \end{rem}

\ssec{Weakly random types}   \label{weaklyrandom} Let $M$ be an $\aleph_1$-saturated model, and let $\mu$ be a finitely additive measure on formulas
over $M$, or just on Boolean combinations of formulas $\phi(x,b)$.  
A $\phi$-type $p$ over $M$ is called {\em weakly random} if any formula $\psi$ in $p$ has $\mu(\psi)>0$.  
 Let $\sX$ be the compact Hausdorff space of weakly random global $\phi$-types; $\mu$ induces a Borel  probability measure on $\sX$.  
In any theory, a formula dividing (or forking) over $\emptyset$ has measure zero for any $0$-definable measure.  Thus a 
weakly random $\phi$-type cannot fork over $\emptyset$.   In particular if we fix a model $M_0$, it cannot fork over $M_0$. 
In a NIP theory, if a type $p$ over  a saturated model $M \succ M_0$ does not fork over $M_0$, then for any $\theta$ there exists a set of types over $M_0$ 
$I(\theta)$ such that $\theta(x,c) \in p$ iff $tp(c/M_0) \in I(\theta)$.   (See e.g. \cite{NIP1}, 2.11.)      
Hence the set of weakly random types has cardinality at most $\beth_2(|M_0|+|L|)$.   In fact  $\sX$ is  separable (see \cite{chernikov-simon} 2.9 and 2.10)  but I am not sure if it is in general metrizable.   It is so in the case of a smooth measure, or a generically stable type.

\end{section}
\begin{section}{Probability logic} \label{problogic}

Many versions of probability logic were investigated by Keisler and Hoover, see \cite{keisler}, following work of Carnap, Gaifman,   Scott and Krauss; most of these were based on $L_{\omega_1,\omega}$.  We will use here a first-order, real-valued  version based on continuous logic. This enables the use of compactness, and in particular studying families of finite structures via their probability logic limits.

   While it is possible to mix $\inf,\sup$ quantifiers with probability quantifiers, we will be mostly interested  in
{\em pure probability logic}, where only  probability quantifiers are used.

The flavor of this logic is determined by  Hoover's quantifier-elimination theorem \ref{hoover} and, over a model, the independence
theorem   \thmref{overamodel}.   Roughly speaking the first result gives a quantifier-elimination to one block of quantifiers; the latter  says that unexpected interaction among events can occur finitely many times, but is no longer unexpected once seen often enough.  
In the case of binary interactions,  a much more precise description  is available,   \thmref{ind-prob}.   
 What  {\em can} be   interpreted is a compact structure, the 'core', with an action of a compact group on it.
It can be viewed as   the space of Lascar types of singletons.   Each binary relation gives rise to a binary function on this  core, making it into a compact structure;  and the probability quantifiers  induce a measure on it.   Given this core, along with the natural map of the universe into it,  
  all values of all  formulas obtained using probability quantifiers  are completely determined.    See \corref{binary}.
  
  For  binary languages, these results  indicate that probability logic (at least over binary languages) has substantial descriptive value but limited 
  interpretative strength.    Interesting unary and `almost unary' relations can be interpreted using
  probability quantifiers (the almost unary ones are just unary if the Galois group of the theory is trivial.)
  However no new binary or higher relations can be defined, beyond combinations of almost unary ones and the
  originally given quantifier-free formulas.   
  This is a severe
  restriction on the interpretative power of  pure probability logic.    It stands in   contrast to the   limitless interpretative abilities of first-order logic.

 We will extend this to local probability logic, involving a locally compact core and a locally compact group acting on it.      

By a {\em probability quantifier in $x$} we mean a syntactical operation from formulas $\phi(x,y)$ (with $y$ a sequence of   variables distinct from $x$) to formulas $E_x \phi$ in the variables $y$ \footnote{we assume at the syntactical level, that if $x'$ is another variable of the same sort as $x$, and $\phi'$ is obtained from $\phi$
by using $x'$ in place of $x$, then $E_x \phi = E_{x'} \phi'$.}
satisfying: 
\begin{enumerate}
\item  $E_x(1) =1$ for the constant function   $1$ (viewed as a function of any set of variables);
\item   $E_x(\phi+\phi')= E_x(\phi)+E_x(\phi')$,  and \\
 
  $E_x(\psi \cdot \phi) = \psi\cdot E_x(\phi)$ when $x$ is not free in $\psi$;
\item  $E_x(|\phi|) \geq 0$. 
 
\end{enumerate}
 
Note that (1-3) are universal, first order axioms.  

If  quantifiers are used,
(3) (applied to $(\sup_x \phi - \phi)$ becomes equivalent to 

(4) \ \  $E_x(\phi) \leq \sup_x \phi$.  

However we will be interested especially in formulas of pure probability logic $\phi$, that do not involve quantifiers, so that it is preferable to have (4)
explicitly.

\begin{rem} \rm If axioms (1-3) hold for pure probability logic formulas,   (4) will hold in   existentially closed models for this sublanguage.   To see this
let $L_{pr}$ denote the formulas obtained from basic ones using expectation quantifiers alone; view all formulas of $L_{pr}$ as basic.  Let $M$ be an existentially closed model for an $L_{pr}$-universal theory $T_{pr}$ including 
 axioms (1-3).    If  (4) fails,  then  for some $\e>0$,
with $\rho(y)=  E_x(\phi(x,y))-\e$, $\phi(x,y) \leq \rho(y) $ must follow from some $\theta (y,y')$
where $\theta(b,b')$ holds for some $b,b'$ from $M$.  But then $E_x \phi(x,y) \leq E_x \rho(y) = \rho(y) = E_x(\phi(x,y)) - \e$, a contradiction.  
Thus (1-3) suffice to axiomatize the pure probability logic validities, i.e. the 
universal sentences applied to formulas using $E_x$ but  no other $x$-quantifiers.
\end{rem}
 
 The name {\em probability quantifier} arises from the situation where $\phi$ takes values in $\{0,1\}$; in general one might also call it an
 {\em expectation quantifier}.   

\ssec{Semantics}

  If $M$ is any   model 
  of axioms (1-4), and $S_x(M)$ is the type space over $M$ in variable $x$, we obtain a positive linear functional 
  on  a dense subset of  $C(S_x(M))$, namely the interpretations of formulas $\phi(x)$ with free variable $x$ and parameters in $M$; it follows from the last axiom that if $\phi(x,b) $ defines the same function as $ \phi'(x,b')$,  then $E_x\phi(x,b)=E_x \phi'(x,b')$.  
  
  By the Riesz representation theorem,
  there exists a unique regular Borel measure $\mu_x$ on $S_x(M)$ with $\int \phi(x) d\mu(x) = (E_x \phi)$.
  {\em This is the intended semantics, and  constitutes the completeness theorem for expectation logic.}
  
  For pure probability logic, we take $S_x(M)$ to be the quantifier-free type space.
  
For example, the  volume of an $r$-ball $B$ around $b$ is determined by the pure probability logic type of $b$.   This is evident from the semantics
as the Borel measure on the type space over $M$ includes this information; $B$ is a $\bigwedge$-definable set, and corresponds to a closed
subset of $S_x(M)$.   One can also see this directly, but more computationally, by expressing ${\rm vol}(B)$ as the   limit
of $Ex \theta(d(x,b))$,   where $\theta$ is a continuous function into $[0,1]$
supported on $[0,r]$,  approximating the characteristic function of $[0,r]$ in the uniform norm.  
   
  By a {\em stochastic sort}, we mean a sort endowed with such an operation $\phi \mapsto E_x \phi$.  
 We will not necessarily assume that every sort is stochastic.

  \ssec{Fubini} If $X$ and $Y$ are stochastic sorts, with corresponding   probability quantifiers 
 $E_x$ and $E_y$,  we obtain two measures on the variables $(x,y)$, arising from 
$E_xE_y$ and $E_yE_x$.   They agree on  formulas obtained by connectives from formulas in $x$ and formulas in $y$.  
On {\em compact} sorts, this suffices to force the two measures to commute.  In general 
they may not, even if $X=Y$, since a function defined by $\phi$ on $S_{x,y}$ may not be
measurable for the product measure.   We will say that $E_x,E_y$ commute
if   
$E_xE_y \phi = E_y E_x \phi$ for all $\phi$; see \cite{hps}.   We say that Fubini holds if any two stochastic sorts commute.

On the other hand, in the foundations of NIP theories notably, one encounters Keisler measures that do not commute.   Thus we do not include Fubini in the list of axioms, but invoke the assumption  when needed.   
  
 \ssec{Pseudo-finite semantics}  
 The above treatment of probability logic takes as a starting point a family of formulas, closed under an expectation operator, as well as continuous connectives.   This is analogous to a view of logic as a family of formulas, closed under quantifiers and connectives.    In another approach,
 one forms the family of formulas formally by closing the basic relations under continuous connectives, the $\inf_x$ operator and the $(Ez)$ operator
 for stochastic sorts.    Each basic formula $R$ comes with a real interval $I_R$ (so that $R$ takes values in $I_R$) and a 
  uniform continuity modulus $\mu_R$, so that $|R(x_1,\ldots,x_n)- R(y_1,\dots,y_n)| \leq \mu_R (\max d(x_i,y_i))$.   These are propagated 
  to general formulas $\phi$ in the natural way; in particular the  interval and uniform continuity modulus of $(Ez)\phi$ and of $\sup_z \phi$ are defined to be those of $\phi$.  
  Let us explain, given a finite structure or an ultraproduct $M$ of finite structures, how to evaluate each formula $\phi$.
  This is done by induction on the complexity of formulas.   The value $(Ez)\phi(z,m)$ is defined to be be the mean value of $\phi(z,m)$
  (with respect to the counting measure on the relevant sort of $M$.)  The values of $\inf_x \phi$ and $C(\phi_1,\ldots,\phi_k)$ (where $C:\Rr^k \to \Rr$ is continuous) are also defined in the obvious way.

 Note that in a saturated model, where every type avoiding measure-zero formulas is realized, it may
be impossible to avoid nonempty parametrically definable measure zero sets.  In this case, 
the pure probability logic type of an element in this theory need not determine the isomorphism type.

\ssec{Hoover's normal form}  We return to the general setting of probability logic.  
Here is Hoover's theorem on  reducing expectation quantifiers to a single block; it is valid in general in our setting with several stochastic sorts.  
See  \cite{keisler} (in a slightly different setting.)

\begin{thm}[Hoover] \label{hoover} Any formula $\psi(y)$ built using connectives and expectation quantifiers can be approximated by ones of the form
$E_x \phi$, where $\phi(x,y)$ has no (probability) quantifiers, and $x$ is a sequence of variables.   \end{thm}

\prf  Let $\Psi$ be the class of formulas that can be so 
 approximated.  Clearly $\Psi$ contains the quantifier-free formulas, and is closed under probability quantifiers; we have to show
 in addition that $\Psi$ is   closed under 
connectives corresponding to continuous functions $c$.  We give two proofs of this.    The first works directly for
any $c$.  Let  $\bar y$ be a sequence of a large number $N$ of copies of $y$.  By the law of large numbers, 
$E_y \phi$ is approximated by $\frac{\sum_j \phi(y_j)}{N}$, uniformly in  the remaining free variables of $\phi$.  ($\phi$ takes values in a bounded interval, say $[0,1]$; so    $|\phi - E \phi| \leq 1$ and thus $(\phi - E\phi)^2$ has expectation at most $1$.  A weak version of the  law of large numbers now states that $|\frac{\sum_j \phi(y_j)}{N} - E(\phi) | \leq \lambda$ with probability at least $1-1/(N\lambda^2)$.  Taking $\lambda = N^{-1/4}$ will do.)
 So $c (E_y \phi)$ is approximated by 
$E_{  \bar y}   c( \frac{\sum_j \phi(y_j)}{N}) $, uniformly in the other variables.

The second proof was explained to me by Ita\"i Ben Yaacov. It does not require Fubini.  
Using  the Stone-Weierstrass theorem  we may take $c$ to be a polynomial.  This decreed,
 no further approximations are needed; 
the normal form becomes valid purely algebraically.  We have to consider the sum or product of two
expressions $E_x \phi$, $E_y \psi$ where we may assume the quantified variables $x,y$ are disjoint from each other and
the free variables.  In this case the sum is $E_xE_y (\phi+\psi)$ and the product is $E_xE_y \phi \cdot \psi$. 
  \eprf

\ssec{Stability of binary correlations.}

Ben Yaacov proved the stability of the theory of measure algebras in  \cite{cat}. 
Taking the viewpoint of piecewise interpretable structures - in this case measure algebras-
this immediately implies  
stability of $E_x \phi(x,y) \& \psi(x,z)$  in any theory with a real-valued expectation operator. 
The implication was not immediately noticed however, and stability of  $E_x \phi(x,y) \& \psi(x,z)=0$
was  reproved directly in \cite{CH} in a restricted environment, in order to prove the  independence theorem there. This was then transposed to the forking ideal in place of the measure $0$ ideal in \cite{pac} (in finite S1-rank)  and   \cite{kim-pillay} for general simple theories,
yielding the independence theorem for  simple theories.   Here we return to measure correlation and
give the simple proof from \cite{cat}.

\begin{prop}\label{stablei} For any $\phi(x,y)$ and $\psi(x,z)$ valued in $\{0,1\}$,  the formula $\theta(y,z)$ defined by 
$E_x \phi(x,y) \& \psi(x,z)$ is stable.

 More generally,  any $\phi(x,y)$ and $\psi(x,z)$,  the formula  $E_x (\phi(x,y) \cdot \psi(x,z))$ is stable.

In fact,  for any formulas $\phi(x,y)$ and $\psi(x,z)$ valued in a compact subset $C$ of $\Rr$, and any 
 continuous
function $c:C^2 \to \Rr$, the formula $E_x(c(\phi(x,y),\psi(x,z)))$ is stable.  

 \end{prop}

\prf  We prove the second statement first.  
Let $M$ be a model, and let $b_i,c_j \in M$.  Let $S$ be the type space over $M$ in the variable $x$.
The expectation quantifiers induce a measure $\mu$ on $S$, such that $\int \phi(x) d\mu(x) = (E_x \phi)$
for any formula $\phi(x)$ over $M$.  
Now $\phi(x,b_i)$ defines a continuous, bounded real-valued function $f_i$ on $S$; while
$\psi(x,c_j)$ defines $g_j$.  So $f_i,g_j \in L^2(X,\mu)$, and $E_x (\phi(x,b_i) \cdot \psi(x,c_j)) = \int f_ig_j = (f_i,g_j)$.
Thus stability follows from \lemref{l2stable}.  

The first statement is a special case, since $\& = \cdot$
on $\{0,1\}$. 

As for the third statement, we can approximate $c$ 
uniformly by a polynomial; so we may take $c$ to be a polynomial.  Since $E_x$ is additive, we may take $c$
to be a monomial $p_m(u)p_n(v)$, where $p_n$ denotes
the $n$'th power map.    Upon replacing 
$\phi$ by $p_m \circ \phi$ and $\psi$ by $p_n \circ \psi$,
the statement now follows again from the 1st paragraph. 
\eprf

Let $Y$ be the sort of the variables $y$, and let 
  $\std{Y}$  denote the  associated strong type spaces, i.e. $\std{Y}=Y/E$ with $E$ the smallest $\bigwedge$-definable
  co-compact equivalence relation.  Let $\hat{b}$ denote the image of $v$ in $\std{Y}$; similar notation for $z,Z$.  
By \lemref{uniqueness}, it follows that there exists a  function $\alpha: \std{Y} \times \std{Z} \to \Rr$, continuous in each variable,  such that for any $b \in Y, c \in Z$
with $\ind{c}{z}$ we have $E_x (\phi(x,b) \& \psi(x,c)) = \alpha (\hat{b},\hat{c})$.

\ssec{The independence theorem and statistical independence}

 \propref{stablei} and \thmref{uniqueness} combine to yield a basic principle
of probability logic, {\em the independence theorem}.   A qualitative version is true in greater generality for ideals with a certain saturation 
property, (S1),  enjoyed by the measure-zero ideals of measures.    As noted in \cite{nqf}, when we actually have a measure 
it is possible not only to assert that the   value of $\mu(R(a,z) \meet R'(b,z))$ is uniquely determined,  
but   to give an explicit formula for it.    We include a proof here, though the independence theorem will only be used in qualitative form later on.  

 Towsner \cite{towsner} noted the relation to combinatorial results and gave a proof of a related statement 
 of `triangle-removal' type for $n$- amalgamation over a model,  by $L^2$-methods.

The   proof  given below reconciles these two approaches; the structure interprets (piecewise) a Hilbert space, where stability reigns in the qf part; this can be viewed as the true source of the stability of the formula in question.
The parameterizing sorts are not required to carry a measure, and the exceptional set is recognized explicity.

    We assume $Z$ is a sort with expectation operators,
$X,X'$ are two other sorts, $R \subset X \times Z$ and $R' \subset X' \times Z$ are two relations.
For some purposes we will assume that  $X,X'$ also carry expectation operators, and that Fubini holds; this will
be stated explicitly.

Let us say that two functions $f,f'$ on a measure space $(X,\mu)$ are  independent if 
for any Borel  $B,B' \subset \Rr$, $f \inv(B)$ and $(f') \inv(B')$ are statistically independent events.
Equivalently, for any two bounded Borel functions $e,e'$ on $\Rr$, $E((e \circ f) \cdot (e' \circ f') )= 
 E(e \circ f) E(e' \circ f')$.  If $f$ and $f'$ are characteristic functions of two events,    this is the usual notion of 
 statistical independence.   
 
 We say that $f,f'$ are  independent {\em over} a $\si$-subalgebra $\mathcal{B}$ of the measure algebra if for each such $e,e'$,
 the conditional probabilities relative to $\mathcal{B}$ satisfy:  
 $E((e \circ f) \cdot (e' \circ f') : \mathcal{B} )=  E(e \circ f:  \mathcal{B}) E(e' \circ f' : \mathcal{B})$.  
 
A suggestive case occurs topologically when $\pi: X \to Y$ is a continuous map of Polish spaces,  $\mu_y$ is a Borel family of measures on the fibers, $\nu = \pi_* \mu$,    $\mathcal{B} $ is the measure algebra of $\nu$, and 
 $\mu= \int_{y} \mu_y$; this means that for any continuous function $\phi$ on $X$ we have $\int \phi d\mu (x)  = \int (\int \phi(x) d\nu_y(x))  d\nu (y)$.
 In this case $f,f'$ are independent iff  
  for almost all $y \in Y$, $f,f'$ are independent with respect to $\mu_y$.    We will also say in this case that $f,f'$ are
  statistically independent over $Y$.

\begin{defn}  \label{ind-def}Write $\Ind{A}{B}{C}$ if for any stable continuous logic formula $\phi(x,y)$ over $C$, and tuples $a$ from $A \union C$, $b$ from $B \union C$, $tp(a/b)$ does not $\phi$-divide over $C$.  In other words, if $\phi(a,b)=\alpha$, then for any   indiscernible sequence $(b,b_1,b_2,\cdots) $ over $C$ and any $\e>0$ and $n$, there exists $a'$ with $| \phi(a',b_i) - \alpha | < \e$ for $i \leq n$.   If $C = \emptyset$, we write $\ind{A}{B}$.

If we restrict to stable formulas $\phi(x,y)$ with $\phi$ defined over $\emptyset$, write $\Ind{A}{B}{0;C}$.
\end{defn}

We use continuous logic formulas in this definition even if $T$ is a first-order theory. 

Let $R(x,y) = \bigwedge_{i} R_i(x,y)$, where the  family $\{R_i: i \in I\}$ can be taken to be closed under finite conjunctions.   
Assume $R(x,b)$ divides, i.e. there exists an   indiscernible sequence $(b,b_1,b_2,\cdots) $ such that $\meet_j R(x,b_j)$ is inconsistent.
Let $\mu$ be a definable  measure, or more generally an invariant measure (i.e. $\mu( \phi(x,c))$ depends only on $\phi$ and on $tp(c)$.)  
Then $\mu(R(x,b))=0$, in the strong sense that $\mu(R_i(x,b))=0$ for some $i \in I$.  See \cite{nqf} 2.9.    Thus if $\psi(x)$ is a definable set of positive measure, then for any $B$, in some elementary extension, there exists $a$ with $\psi(a)$ and such that no $R(a,b)$ holds, with $b \in B$, if $R$ is a $\bigwedge$-definable relation such that $R(x,b)$ divides.   This can be applied to 
any specific value  $\alpha$ (or any closed range of values) of a continuous logic formula $\phi(x,y)$, letting $R(a,b)$ hold if $\phi(a,b)=\alpha$.

Let $cl(A)$ denote the bounded closure, or continuous-logic algebraic closure in the sense of \cite{byu}.
Thus a type over $cl(A)$ is the same as a Kim-Pillay strong type over $A$.

 \begin{lem}  \label{l2}  Stable independence has the properties of: \begin{itemize}
 \item Symmetry:  $\Ind{A}{B}{C}$ implies $\Ind{B}{A}{C}$.
 \item   trivial monotonicity:   $\Ind{A}{B}{C}$ implies $\Ind{A}{B'}{C}$ if $C \subseteq B' \subseteq B$. \\
 Also, $\Ind{A}{B}{0;C}$ implies $\Ind{A}{B}{0;C'}$ if $C \subseteq C' \subseteq B$.
 \item   finite   character, 
 \item small bases:  for any $A,B$ there exists $C \subseteq B$, $|C| \leq |A|+|L|$ with $\Ind{A}{B}{C}$.
 \item   Existence and stable stationarity:   For any $b$ and 
any Kim-Pillay type $Q$ (over $C$), there exists $a \in Q$ with $a,b$ stably independent over $C$.    Moreover
for any stable formula $\phi(x,y)$ over $C$, the truth value of $\phi(a,b)$ is the same for all such $a$.  
 \item Transitivity: if $a,b$ are stably independent over $A$ and $a,c$ are stably independent over $cl(A \union \{b\})$, then $a$ is stably independent over $A$ from $(b,c)$.  
 
 \end{itemize}
\end{lem}
\prf   All but transitivity follows directly  from \cite{byu}.   (Transitivity is proved there under the assumption of global stability; we check it here under our more local assumptions.)     
 
Transitivity: Work over $A$. Let
 $\phi(x;y,z)$ be stable.  Note that any instance $\phi(x;b,z)$ is stable
(there are no sequences $(a_i,b_i,c_i)$ with $\phi(a_i;b_j,c_j)$ iff $i<j$; in particular no such sequences with all $b_j=b$.)   Let 
 $p=tp_{KP}(a)$, and let $\psi(y,z)$ be
the $p$-definition of $\phi$.   Then  
  $\psi_b = \psi(b,z)$ is the $tp(a/cl(b))$-definition of $\phi_b(x,z) = \phi(x;b,z)$.  This is because there exists a ('Morley') sequence $a_i$
  such that for any $(b',c')$, $\lim_i \phi(a_i,b',c') = \psi(b',c')$; in particular this holds for $b'=b$.  
  On the other hand the existence of such an indiscernible sequence implies that for any $\beta \neq \psi(b,c')$,
  $  \phi(a,b,y)=\beta $ divides for $y \models tp(c'/b)$; so the  $tp(a/cl(b))$-definition of $\phi_b(x,z)$
  must be $\psi(z)$.
   
 Now assume $a,b$ are stably independent over $A$ and $a,c$ are independent over $cl(Ab)$.   Then 
 $\phi(a;b,c) $ holds iff $  \psi_b(c)$ iff $ \psi(b,c)$.   This shows the stable independence of $a$ and $(b,c)$.
 \eprf

\ssec{Analytic structures viewed as interpretable}

 Let $\mu$ be a definable measure (in variable $x$) for a theory $T$, for instance obtained using expectation quantifiers.  
We will view the Hilbert space $L^2(\mu)$ as  piecewise-interpretable in   $T$.  For any model $M \models T$,
we have $L^2(\mu)(M) = L^2(S_x(M),\mu)$.

Then the Hilbert space formulas provide us with stable formulas of $T$, in the sense of continuous logic; and the results
on stable independence apply.   

For our purposes, we could use the theory of  probability algebras in place of the theory of Hilbert spaces;
the probability algebra $B(\mu)$
 can be identified with the   elements of $L^2(\mu)$ represented by $\{0,1\}$-valued Borel functions, but we have not only the induced norm from $L^2(\mu)$ but also {\em multiplication} as part of the structure.  
 Stability of $B(\mu)$ can in any case be deduced from that of $L^2(\mu)$.   
 For other applications, the Banach lattice $L^1(\mu)$ will be needed.     
 The Hilbert space picture is   appealing 
 in particular in connection with the Peter-Weyl theorem, and the representation of the automorphism group of the algebraic (bounded) hyperimaginares.  
 For definiteness  we will talk
about Hilbert spaces below, but the discussion would be the same for the others.

Let $\Uu$ be a large saturated model of $T$.  Let $H$ be the Hilbert space $H=L^2(S_X(\Uu))$.

For any small substructure $A$, possibly including (hyper)imaginaries, we define $H_A$ to be the subspace of $H$
fixed by $Aut(\Uu/A)$.  

\begin{rem} \rm  There is a canonical embedding of $L^2(S_XA)$ into $H$, falling into $H_A$
(namely, $f \mapsto f \circ r$, where $r$ is the restriction $S_X(\Uu) \to S_X(A)$.)      This will be surjective   assuming a certain 'strong germ' property;
without such an assumption, there may exist a family $(D_c: c \in Q)$ of definable sets such that  $\mu(D_c \triangle D_{c'} ) = 0$ for all $c,c' \in Q$, but no $A$-definable set is equivalent to any $D_c$.  
 We will not make this assumption,  so  the image 
of $L^2(S_XA)$ in $H$ may be smaller than $H_A$.
    \end{rem}

The elements of $H$ can be viewed as hyperimaginaries of $\Uu$; in fact $H$ is piecewise interpretable
in $\Uu$, in a sense that will now explain.      An element $\xi$ of $L^2(S_X(\Uu))$ can be   approximated by continuous functions, given as the value of a formula $\phi(x,a)$ of real-valued continuous logic.  Thus 
 
\[\xi = \lim_n \phi_n(x,a_n) \]
with the limit taken in the $L^2$-norm. Let $\kappa= \2norm{\xi}$
and $\phi_n(x,a_n)$ of $L^2$-norm $\leq k+1$. 
 Moreover we can choose the sequence with  $\2norm{ \phi_n(x,a_n) - \xi } \leq 2^{-(n+1)}$, so that 
\[  \2norm{ \phi_n(x,a_n) - \phi_n(x,a_{n+1})}  \leq 2^{-n}, \ \ \  \2norm{ \phi_n(x,a_n) } \leq \kappa +1 \]
Let $\bar{a}$ be the sequence $(a_n)$, and $\bar{\phi}$ the sequence $\phi_n$, and write 
  $L_{\bar{\phi}}$ for the set of sequences $\bar{a}$ satisfying the displayed formula.
The number of possibilities for $\bar{\phi}$ is bounded (by $|L|^{\aleph_0}$); the  sets 
$L_{\bar{\phi}}$ are easily seen to be directed, under inclusion (by using disjunctions, and using parameters to
choose the appropriate disjunct.)  
For each such $\bar{\phi}$,  $L_{\bar{\phi}}$ is a  $\bigwedge$-definable set.  In this sense, $\union_{\bar{\phi}}L_{\bar{\phi}}$ is piecewise $\bigwedge$-definable.   The equivalence relation 
$\lim \phi_n(x,a_n) = \lim \phi_n(x,b_n)$ is also  $\bigwedge$-definable; it is equivalent to 
\[ \bigwedge_n E_x ((\phi_n(x,a_n) - \phi_n(x,b_n))^2)   \leq 2^{-(n-2)} \]
where $E_x$ is the expectation.  Moreover, within $L_{\bar{\phi}}$, the relations $a+b=c$, $\alpha a = a'$ (for $\alpha \in \Rr$) are $\bigwedge$-definable, and so is the formula giving the inner product; it is   approximated
by $E_x \phi_n(x,a_n) \phi_n(x,b_n)$, uniformly in $(a,b)$ for $a,b$ such that $\2norm{a}, \2norm{b}$ is bounded.

 Let $L_H$ be the language of $H$, i.e. in our case the Hilbert space language. 
 
   \begin{lem} \label{7.0}Assume $A \leq B \meet C$,    $A=cl(A)$ and $B,C$ are stably independent over $A$.   
Let $H$ be a stable structure, piecewise interpretable   in the sense considered above.  
   Assume $H_A$ is closed in $H$, i.e. every $H$-hyperimaginary element that is bounded over $H_A$ lies in $H_A$.  
Then $H_B,H_C$ are $L_H$- independent  over $H'_A:= H^{eq} \meet A$.
\end{lem}
 
\prf    Let  $b \in H_B$, $c \in H_C$  and let $\phi(x,y)$ be a (stable) $L_H$-formula.  We view $b,c$ as hyperimaginary elements of $M$.  Then
$tp(b/A)$ is consistent with   an $A$- definable $\phi$-type $p(x)$.  Note that the canonical base of $p$,
in the sense of $L_H$, is then defined over $A$ and lies in $H^{eq}$, so $p$ is defined over
$H'_A$.  By the  stable independence of $B,C$ over $A$, $\phi(b,c)$ holds iff  $\phi(x,c) \in   p |C$.
As this is true for every $L_H$-formula $\phi$, $b,c$ are $L_H$-independent over $H'_A$.

\eprf

 In the case of Hilbert spaces, it is known \cite{byb} that for any 
imaginary $h$, we have $cl(h) = cl(e)$ with $e$ coding a finite dimensional subspace $E$ of $H$.   
(Note that the lemma refers to hyperimaginaries of $H$ itself, not to the induced structure from the piecewise-interpretatoin.)  
   Moreover,  $E \subseteq cl(h)$.  (It suffices to show
 that $E_1 :=\{a \in E:  \2norm{a} = 1\}  \subset cl(e)$.      But $E_1$ is compact, and $e$-definable, so by definition of `closed' for continuous logic, 
 $E_1 \subseteq  cl(e)$. )  Thus we obtain a simpler form, \lemref{7}.     
   
\begin{lem} \label{7}Assume $A \leq B \meet C$,    $A=cl(A)$ and $B,C$ are stably independent over $A$.
Then $H_B,H_C$ are independent as Hilbert spaces over $H_A$.
\end{lem}

Let $Z$ be a stochastic sort.

We have a regular Borel measure $\mu$ on $S_Z(\Uu)$ induced by the $Z$-expectation operators.
Let $H$ be the associated piecewise-hyperdefinable Hilbert space, and $B$ the associated  probability algebra.   
Let $H_{\0},B_{\0}$ be the subspace (respectively subalgebra) of $H,B$ consisting of elements with   bounded orbit
over $\emptyset$ (almost invariant).   We will also write simply $\0$ for $B_{\0}$; so $H_{\0}$ can be identified with $L^2(B_{\0})$.  
  $B_{\0}$ consists of a bounded number of elements, each of which lives in some hyperdefinable piece of $B$ and hence can be identified
with a hyperdefinable element.

Hence two elements with the same Kim-Pillay type have the same type over $B_{\0}$; the  equivalence relation $tp(x/B_{\0})=tp(y/B_{\0})$ is $\bigwedge$-definable, and with boundedly many classes.   

Let $\hat{X}$ denote the set of Kim-Pillay strong types of elements of $X$; so $\hat{X}=X/E$ for a certain co-bounded $\bigwedge$-definable equivalence 
relation on $X$.  
For $a \in X(\Uu)$, let $\hat{a} $ denote the image of $a$ in $\hat{X}$.  

Since $B_{\0}$ consists of boundedly many hyperimaginaries,  $\hat{a}$ determines $tp(a/B_{\0})$.

Let $\Xi=\hat{Z}$ be the maximal bounded quotient of $Z$ by a $\bigwedge$-definable equivalence relation.  $\Xi$   carries a Borel probability measure, the pushforward $\mu_0$ of $\mu$ on  $S_Z(\Uu)$ (or any $S_Z(M)$.)

      Recall that we write $\ind{a}{b}$ if  $tp(a/b)$ does not divide over $\emptyset$ with respect to any stable formula.

\begin{thm}[Independence theorem for probability logic] \label{ind-prob} \, \\  Let $(Z,\mu)$ be a stochastic sort.  
\begin{enumerate}
\item    Let $a \in X(M)$, $b \in Y(M)$, with 
  $\ind{a}{b}$.  Then   $R(a,z), R'(b,z)$ are statistically  independent over $\0$.

 \item    Assume   $X$ carries a definable measure  $\nu$,  and $\nu, \mu$ commute (i.e. Fubini holds).   
  Then there exist a $\bigwedge$-definable set  $X^* \subset X$ of full measure,  such that if $a \in X^*$, $b \in Y(M)$, and 
  $\ind{a}{b}$, then   $R(a,z), R'(b,z)$ are statistically  independent  over $\Xi$.
 
\end{enumerate}
   \end{thm} 
   
   In the proof below, we will use the measure space  $(X,\0,\mu|\0)$;  we write an integral of a measurable function $f$ with respect to this space
  simply as $\int_{\0} f$.  We will also use the compact space $\hat{X}$;  as a measure space,  it is always understood to carry  the   measure derived by pushforward from $\nu$ (and also denoted $\nu$).  $\Xi$ always carries the pushforward measure from $S_Z(\Uu)$.  
  We denote conditional expectation with respect to a subalgebra $A$ by $E(f|A)$; when $A$ is the measure algebra of a space $S$,
  e.g. of $S=\hat{X} \times \Xi$,  we will also  write $E(f | S)$ for $E(f|A)$.

 \begin{proof}[Proof of \thmref{ind-prob} (1)] 
  
 For $a \in X(M)$, define a function $\ha: S \to \Rr$,  $q \mapsto R(a,c)$ where $c \models q$; similarly
$\hb$.  
Let $e,e'$ be   bounded Borel functions on $\Rr$, and write ${a^e}$ for $e \circ \ha$.  
let $E_r(\phi) =E(\phi  | \0)$ denote the conditional expectation to $\0$.   
 We have to show that as $\0$-measurable functions, we have  
$ E_r({a^e} \cdot {b^{e'}} )= E_r({a^e}) E_r(b^{e'}) $.
Equivalently, for any  bounded   $\phi \in B_{\0}$, integrating with respect to $(\0,\mu|\0)$ we have:
\[  \int \phi(r) E_r({a^e} \cdot {b^{e'}} )  = \int \phi(r) E_r{a^e} E_r{b^{e'}}   \]
We can restrict to the characteristic functions of clopen sets, and view $\phi$ as a $\{0,1\}$-valued,
almost invariant Borel function $S_z(M)$.

  Let $a_0$ be the orthogonal
projection to ${H_{\0}}$ of ${\phi{a^e}}$, $b_{\0}$ the orthogonal
projection to ${H_{\0}}$ of ${{b^{e'}}}$. 
For any element $g \in {H_{\0}}$,   
we have 
$   \int g(r) \phi(r) E_r {a^e}  = \int {g \circ \pi}(r)  (\phi \circ \pi){a^e}  = (g,\phi {a^e})_H = (g,a_0)_{{H_{\0}}}$.  
so $  a_0(r) =  E_r (\phi(r) {a^e}) $.   Similarly with the natural notation, $  b_0(r) =  E_r({b^{e'}}) $.   Thus what we have
to show is:
\[  \int \phi(r)  {a^e} \cdot {b^{e'}}    = \int   a_0(r) b_0(r)   \]
or in Hilbert space notation
\[ (\phi   a^e, b^{e'}) _E = (  a_0, b_0)_{{H_{\0}}} \]

Now ${\phi{a^e}} - a_0 \perp {H_{\0}}$, and   ${{b^{e'}}} - b_0 \perp {H_{\0}}$.  In particular
$(\phi {{a^e}} - a_0,b_0) = (a_0,{{b^{e'}}} - b_0) = 0$.  Thus what we need is
that $\phi{{a^e}} - a_0$ and ${{b^{e'}}} - b_0$ are orthogonal vectors.    This follows from Hilbert space independence
of $a^{e},b^{e'}$, that we have by \lemref{7}.   
\end{proof}

The proof of \thmref{ind-prob}(2) involves three steps.    (a)   Define a  subset $X^*$ of $X$, consisting of the elements
whose $R$-interaction with $\0$ factors through $\Xi$, and show it is $\bigwedge$-definable.

(b) Show $\nu(X^*)=1$.   It will be useful to prove this under a weaker assumption than stated, namely   that $\nu$ is
$\hat{X}$-definable rather than $0$-definable.

(c)  For $a \in X^*$, and any $b \in Y$ with $\ind{a}{b}$, $R(a,z)$, $R(b,z)$ are statistically independent over $\Xi$.

Towards (a), let $ \beta_a $ be the orthogonal projection of $R(a,z)$ to $H_{\0}$.   

 We noted earlier, using the fact that the elements of $B_{\0}$ are hyperdefinable, that $B_{\0} \subset \dcl(\Xi)$.  Hence $\beta_a$ depends only on $q=tp(a/\Xi)$, and thus only on $\hat{a}$.  We may also write $\beta_{\hat{q}}$.
 
 Let $\alpha_{\hat{a}}$ be the  orthogonal projection of $ R(a,z)$ (or equivalently of $\beta_a$)  to $L^2(\Xi)$.

Define  

\[ X^* := \{a \in X:  \beta_a \in L^2(\Xi) \}  \]
Equivalently, $a \in X^*$ if $\alpha_{\hat{a}} =\beta_{\hat{a}}$.

\begin{lem} \label{ind-prob-a}   $X^*$ is $\bigwedge$-definable.   In fact it given by a conjunction of  pure probability logic formulas with parameters.
 
  \end{lem}
\prf   Let $(c_i)$ be a basis for the orthogonal complement of $L^2(\Xi)$ in $H_{\0}$.   
We have $a \in X^*$ iff $(R(x,a),c_i)=0$ for each $i$.   By piecewise $\bigwedge$ interepretability of the Hilbert space structure,
each map $a \mapsto (R(x,a),c_i)$ is (a uniform limit of) definable relations on $X$.  So $X^*$ is $\bigwedge$-definable. 

To see that the relevant formulas use only expectation quantifiers, we may assume $||R(a,z)||_2 \leq 1$.  
Approximate $c_i$ by continuous functions $c_{i,j} \in C(S_Z(\Uu))$,  so that $||c_i - c_{i,j}||_{L^2(\Xi)} < 2^{-j}$.   Then the condition  $(R(x,a),c_i)=0$ can be written as the conjunction of $(R(x,a),c_{i,j}) \leq 2^{-j}$.   Now  $(R(x,a),c_{i,j}) = E_z R(x,z) c_{i,j}(z)$
is visibly obtained by an expectation quantifier.  
\eprf 

(Is $X^*$ $\bigwedge$-definable by pure probability logic formulas without parameters?)

Next we prove (b).  This concerns only $(X,Z,R)$, so $Y,R'$ are not involved, and in the proof of \lemref{ind-prob-lem2} we will use the letter $y$ 
for a second variable  ranging over $X$. 

Using the two projections from $S_{xz}(\Uu)$, we can view both   measure algebra of $\hat{X}$   and $\0$
as subalgebras of the measure algebra $B_{xz}(\Uu)$.
  Let $\0[\hat{X}]$ denote the $\si$-subalgebra
of $B_{xz}(\Uu)$ they generate; this is just the measure algebra of the product measure space  of $(Z,\0) \times (\hat{X},\nu)$. 

\begin{lem}\label{ind-prob-lem2}  $\nu(X^*)=1$.
\end{lem}

\prf  We let $x,y$ range over $(X,\nu)$, while $z$ ranges over $(Z,\mu)$.  Since $\mu,\nu$ commute, we have, as our principal use of Fubini,
\[E_xE_yE_z (R(x,z) \wedge R(y,z)) = E_z E_x E_y  (R(x,z) \wedge R(y,z)).\]
Now $E_x E_y  (R(x,z) \wedge R(y,z)) = E_x R(x,z) E_y R(y,z) = ( E_x R(x,z))^2$ by property (2) of probability quantifiers; so 

\begin{equation}   \label{b1}  E_xE_yE_z (R(x,z) \wedge R(y,z)) = E_z ((E_x R(x,z))^2 )  \end{equation}

 Let $\beta= E( R(x,z) |  \0[\hat{X}] )$ be the conditional expectation of $R(x,z)$ (as an element of $L^2(S_{x,z}(\Uu)$) to 
 an element of  $L^2(  \0[\hat{X}] )$.  

We will express  each side of \eqref{b1}  in terms of   $\beta$.   First, using \thmref{ind-prob}(1), for any $a,b \in X$,
\[E(R(a,z) \wedge R(b,z) | \0) = E(R(a,z):\0 ) E(R(b,z) | \0) =    \beta_a \beta_b\]
so 
\[E_z(R(a,z) \wedge R(b,z) ) = \int_{\0} (E(R(a,z):\0 ) E(R(b,z) | \0) =  \int_{\0}   \beta_a \beta_b   =  \int_{\0}   \beta_{\hat{a}} \beta_{\hat{b}} \]
where $\int_{\0}$ denotes the integral of a $\0$-measurable function (note all our functions are bounded and so integrable.).   
The last step expresses the fact that   $\hat{a}$ determines $tp(a/\0)$, and in particular $\beta_{a} = \beta_{\hat{a}}$ depends only  on $\hat{a}$.

  Note that  $\int_{\hat{a} \in \hat{X}} \int _{\hat{b} \in \hat{X}}  \beta_{\hat{a}} \beta_{\hat{b}} =
 ( \int_{\hat{a} \in \hat{X}}   \beta_{\hat{a}})^2$.  But $\int_{\hat{a} \in \hat{X}}  \beta_{\hat{a}} = E(\beta|\0)$.  
 Thus   $E_xE_yE_z (R(x,z) \wedge R(y,z)) =$
\beq{b2}    \int_{\hat{a} \in \hat{X}} \int _{\hat{b} \in \hat{X}} \int_{\0} \beta_{\hat{a}} \beta_{\hat{b}} = 
  \int_{\0}  
  E(\beta | \0)^2  = || E(\beta|\0)||_2^2   \eeq
where the norm is taken in $H_{\0}$.   Here we used ordinary Fubini, allowed since the integrand is measurable for the product measure.

On the other hand, for $c \in Z$, $E_x R(x,c) $, as the value at $c$ of a pure probability formula with respect to $\nu$,
  depends only on the pure probability type of $c$ over $\nu$;  and hence only on $\hat{c}$, the image of $c$ in $\Xi$, which determines
  $tp(c/\hat{X})$.  (Here we use definability of $\nu$ over $\hat{X}$-parameters.)    By factoring $E(\beta | \Xi)$ through $E(\beta |  \hat{X} \times \Xi])$, we see that for almost all $\hat{c}$,    $E_x R(x,\hat{c}) =E(\beta| \Xi)(\hat{c})$.   
   
  Thus, squaring this equality and integrating now over $\hat{c} \in \Xi$, 

\beq{b3} E_z ((E_x R(x,z))^2 ) =  ||   E(\beta| \Xi)    ||_2^2            \eeq

where now the norm is taken in $L^2(\Xi)$.  By (1,2,3) we have:
\[ || E(\beta |\0) ||_{L^2(\0)} =    || E(\beta| \Xi) ||  _{L^2(\Xi)}      \]

Now $  E(\beta| \Xi)   $ is the orthogonal projection to $L^2(\Xi)$ of $(E(\beta | \0) \in H_{\0}$.
Since they have the same norm,  we must have $E(\beta | \0) \in L^2( \Xi)$.

Let $\psi$ be any continuous  `test function' on $\hat{X}$.   Then $\psi \nu$ is another   measure on $X$; it may not be $0$-definable 
but it is $\hat{X}$-definable, so that the above result applies.  We obtain:
$E( \psi  \beta  | \0 ) \in L^2(\Xi)$.   By factoring first through the product algebra $B_{\hat{X}} \times \0$, 
it follows that for almost all $q \in \hat{X}$, $E(\beta_q) \in L^2(\Xi)$, concluding the proof.

\eprf

    \begin{proof}[Proof of \thmref{ind-prob}(2)]
   It remains to prove  (c).  Let $a \in X^*$.   Recall $\beta_a \in L^2(\0)$ is the  conditional expectation of $R(a,z)$
  relative to $\0$,  equivalently the orthogonal projection of $R(a,z) \in L^2(S_M(Z))$ to $L^2({\0})$;  while $\alpha_a$ is 
   the orthogonal projection  of $R(a,z)$ to $L^2(\Xi)$;  $\alpha_a$ depends only on $tp(a/\Xi)$ and hence on $q= \hat{a}$.  
    Similarly define $\beta'_b,\alpha'_b$ for $R'$.      Then by \thmref{ind-prob}(1), 
  \[ E_z(R(a,z) \wedge R'(b,z)) =      \int_{r \in S(\0)} \beta_a \beta'_b = (\beta_a,\beta'_b) . \]
  The last term is the inner product in $L^2(\0)$.  
      By definition of $X^*$, we have   $\beta_a = \alpha_a$,  so   
   \[ (\beta_a ,\beta'_b)  =   (\alpha_a, \beta'_b) =  (\alpha_a,\alpha'_b)    \]
  ( The last equality is by the characteristic property of   orthogonal projections to closed subspaces of Hilbert space, $(P(u),v) = (P(u),P(v))$.)
Thus
    \[ E_z(R(a,z) \wedge R'(b,z))      =  (\alpha_a,\alpha'_b)      \]
 where now the inner product is computed in $L^2(\Xi)$, and hence proves independence over $\Xi$.

\end{proof}

 \begin{rem} \label{ind-prob-rem}  \rm   Assume $L$ is countable.
 \begin{enumerate}  \item   Asides from the sharper conclusion, 
\thmref{ind-prob} has a considerably wider domain of applicability than a purely $L^2$-based statement such as \thmref{overamodel}, which applies only to a random $2$-type.
 For example \thmref{ind-prob}  applies when $tp(a/\Xi)=tp(b/\Xi)$, frequently an important situation, though $tp(a,b )$ is certainly not random in this case.  An example of this is given in Appendix C.
 
 \item Note $X^*$ is defined in terms of $(X,Z,\mu,R)$ alone.  It is shown to have full $\nu$-measure for   any commuting $\nu$.   And if $a \in X^*$,
 statistical independence over $\Xi$ is proved for {\em{any }}$(Y,R' \subset Y \times Z)$.
 
 \item A variation: Let $(Z,\mu)$ be a stochastic sort.   Let $a \in X(M)$, $b \in Y(M)$, with 
  $\ind{a}{b}$.  Assume $X$ has a $\Xi$-definable measure $\nu$ commuting with $\mu$,  and concentrating on $tp(a)$.  
   Then  
  $R(a,z), R'(b,z)$ are statistically  independent over $\Xi$.      It suffices for $\nu$ to be Borel-definable, in the sense of \cite{NIP1}.   
  We   do not use self-commutation of $\nu$!   
  This is proved in the same way as \thmref{ind-prob} (2), but more easily; in  \lemref{ind-prob-lem2} integration over $\hat{X}$  becomes unnecessary, 
  since only one strong type is involved.

\item \label{forkingrare} Let $(Z_1,\mu_1)$ and $(Z_2,\mu_2)$ be   stochastic sorts.
For a measure one set of types $q_2$
 on  $Z^2$, if $(a,b) \models q_2$ then $\ind{a}{b}$.   Here `types' can be taken to be $\Phi$-types
 where $\Phi$ is the family of all 
 stable probability logic formulas.   
  
 \item \label{ind-prob-rem-n}  Let $(Z,\mu)$ be a self-commuting stochastic sort, and $R_i \subset Z^2$ a definable binary relation.    Then for almost all  types $q$ on $Z^n$, if 
 $(a_1,\ldots,a_n) \models q$ then  
  
 the events $R_i(a_i,z)$ ($i=1,\ldots,n$) are independent over $\0$, and over $\Xi$ in case $\mu$ is self-commuting.    This follows inductively from \ref{forkingrare} and \thmref{ind-prob},
 taking at first $X=X^{n-1}, Y=Z=X$ to obtain that $\bigwedge_{i \leq n-1} R(a_i,z)$ is statistically independent from $R(a_n,z)$ over $\Xi$.   
  
 \item  \label{fullspace}We used here the full Kim-Pillay space, without restricting the level of definability of the implied  $\bigwedge$- definable equivalence relations.  
 This is inevitable due to the starting data;
our notion of independence uses the complete type of $a$ and of $b$; in particular if $a$ is $0$-definable  or lies in the bounded closure of $0$, via a  
formula involving quantifiers, 
then $\ind{a}{a}$ holds.   
On the other hand the deduction of (2) from (1) uses probability quantifiers only.   Since the Hilbert space is PPL interpretable, it should be possible to formulate a version of (1) and hence of the full theorem with definability in terms of probability quantifiers, given a stronger assumption of independence at the quantifier-free level.
 
  \end{enumerate}
 \end{rem}

\ssec{Interpretative power of probability logic in a binary relational language} 
 
Let $M$ be an $L$-structure with all sorts stochastic with commuting expectation quantifiers;  for simplicity  take a single sort $X$, and assume the measure on $X$ is self-commuting.   

Recall $\hat{X}$ is the biggest bounded quotient of $X$.  If $f: \hat{X} \to \Rr$ is a continuous function, then $\alpha(x)=f(\hat{x})$ is an $M$-definable function.
Let  $M_{\hat{X}}$ be the result of adding to the language all such functions $\alpha$, or equivalently, a countable subset that separates points.
 $M_{\hat{X}}$ is the expansion by all parameterically definable relations that are definable in $M_A$ for {\em any} elementary submodel $A$ of $M$;
 $M_A$ being the expansion of $M$ by constants for $A$.   Moving to $M_{\hat{X}}$ is known as {\em working over the algebraic closure of the empty set.}
 (Here in the sense of continuous logic.)

Assume that $L$ consists of unary and binary relations, and possibly unary function symbols.   Let $L_{qf}$ denote the set of quantifier-free formulas, and $L_{prob} $ the
result of closing $L_{qf}$ under continuous connectives and expectation quantifiers.

\begin{cor} \label{binary} Assume $L$ is binary, and countable.   Away from a measure zero set of $n$-types, the $L_{prob}$ type of a tuple $(a_1,\ldots,a_n)$  of $M$ is determined   by the  quantifier-free type of $(a_1,\ldots,a_n)$ along with
the values of {\em unary} $\hat{X}$-definable formulas $\alpha(a_i)$. 
\end{cor}

Equivalently, for each $\phi \in L_{prob}$    
 there exist  formulas $\alpha_1,\alpha_2,\cdots$ of $L_{\hat{X}}$, each taking a single variable
from among $x_1,\ldots,x_n$,  $L_{qf}$-formulas
$\beta_k(x_i,x_j)$ - each taking two variables from among the $x_i$, and a Borel function 
$\Psi$,  such that $E(|\phi - \Psi(\alpha_1,\alpha_2,\ldots,\beta_1,\beta_2,\ldots)|)=0$;  in other words for almost all $x_1,\ldots,x_n$, 
$\phi = \Psi(\alpha_1,\alpha_2,\ldots,\beta_1,\beta_2,\ldots)$.

\prf   Let us first see that the second statement follows from the first.  Let $Z_n$  be the space of $L_{prob}$-types on $X^n$, $Z=Z_n \times_{(Z_1)^n} \hat{X}^n$,  $W$ the space of $\Delta$-types where $\Delta$ consists of all 
qf formulas along with $\hat{X}-$ definable unary formulas.  We have a natural restriction map $r: Z \to W$.  By the first statement, there exists a measure-one set $Z' \subseteq Z$ such
that $r$ is injective on $Z'$.  We may take $Z'$ to be an $F_\si$ set, i.e. $Z' = \union_n Z_n$ is a countable union of compacts (seeing that $Z$ is compact.)   Let $W_n = r(Z_n)$.
Then $W_n$ is a closed subset of $W$, and $r \inv$ is continuous on $W_n$.   It follows
that $r \inv$ is Borel on $\union_n W_n$; and any continuous function $\phi$ on $Z$
can be expressed as $\Psi(r(z))$ where $\Psi = \union_n (r \inv | W_n)$.

Next let us prove that the $L_{prob}$ type is indeed determined by the given data.  
The    unary relations $\alpha$ arising from continuous functions on $\hat{X}$ can be recombined to give  the map $x \mapsto \hat{X}$.
So it suffices to show that for   $\phi \in L_{prob}$,   the value $\phi(a_1,\ldots,a_n) $ is determined a.e.
by the quantifier-free type of $(a_1,\ldots,a_n)$ along with the elements $\hat{a_i} \in \hat{X}$.
Using Hoover's theorem \thmref{hoover}, we may take $\phi$ to have the form $E_w \psi(w,a_1,\ldots,a_n)$
 where $w$ may be a tuple , and $\psi$ is quantifier-free. As usual in quantifier-elimination, 
 working inductively, we may assume $w$ is a single variable.  By Stone-Weierstrass 
 we can take $\psi$ to be a polynomial
 in basic formulas $R(w,x_j)$.   Since $E_w$ is additive, it suffices to determine the value of each monomial,
 i.e. of products of such basic relations.   In the presence of function symbols, we view a relation
 $R'(fw,gx_j$) simply as another relation $R''(w,x_j)$.
 We may collect together all relations belonging to a given variable $x_j$ to obtain a single relation $R_j(w,x_j)$.   The value $E_w R_j(w,a_j)$ is determined by $\hat{a_j}$.    Finally the value of $E_w \Pi_j R_j(w,x_j)$ is just the product of these last,  by
    \remref{ind-prob-rem} (\ref{ind-prob-rem-n}).

\eprf

\begin{rem} \begin{enumerate}
\item  Here $\widehat{X}$ should be viewed as a topological {\em structure}.   The relations are these:    for each stable q.f. formula $\phi(x,y)$,
we have a map on $\widehat{X}^2$ giving the generic value of $\phi$ at $(p,q) \in \widehat{X}^2$.
\item It would be interesting to determine when $\Psi$ can be taken to be continuous and not just Borel.   If one is content with quantifier-elimination up to 99\%, rather than almost everywhere, $\Psi$ can be taken to be a continuous function of finitely many variables:    For each $\phi=\phi(x_1,\ldots,x_n) \in L_{prob}$,  there exists a quantifier-free $L_{\hat{X}}$ formula $\phi'$ such that 
$E(|\phi-\phi'|) \leq 0.001$.  
\item  A (real-valued) $\Uu$-definable formula $\psi$ is a  {\em matrix coefficient} if the set of $Aut(\Uu)$-conjuates of $\psi$ spans a finite-dimensional space; equivalently, $\psi$ factors through a  definable map from $X$ to a  $\bigwedge$- interpretable finite-dimensional Hilbert space.  
In place of working over the algebraic closure, one can make 
similar statements in terms of  matrix coefficients maps or in terms of definable maps into $\bigwedge$- interpretable finite-dimensional Hilbert spaces.
\item  We did not restrict the definability level of the  unary maps $\alpha_i$   in \corref{binary}.  
In case we are working over an elementary submodel, it suffices to take qf-definable ones.   If the Galois group
of $\hat{X}$ is trivial, one can take qf-definable maps over a saturated model $M$, with the property that they are invariant
under $Aut(M/\hat{X})$.  In general it should be possible to describe the quotient of $\hat{X}$ we require using probability logic definable functions;  we do not take it up here, but 
see \remref{ind-prob-rem} (\ref{fullspace}).  
\end{enumerate} \end{rem}

This can be read as saying that with pure probability logic, over a binary language\footnote{I.e. the signature has only binary relation and unary function sybmols}, interesting finite or finite-dimensional structures are interpretable
along with a map from $M$ into them; and given these, nothing else can be interpreted that is not visible at the level of basic relation symbols.

The fundamental problem here is to extend the theory  \thmref{ind-prob} to $4$-  amalgmation and higher.   The following weak version would already be useful.    Recall that in the presence of a notion of independence of two substructures over a third, an {\em independent system} of substructures is a family $(A_u: u \in S)$ where $S$ is a simplicial complex, such that $A_u$ is independent from
$\union \{A_v: \neg u \leq v \}$ over $\union \{A_w: w < u \}$.
 
\begin{problem} \label{higher}    
Assume $\mu$ is a strictly definable measure on $X$.    
Does  there exist  a canonical piecewise-interpretable independent system of measure algebras 
 ($S_{\widetilde{u}}:  u \subset [n])$  containing the measure algebras $F(u)$  of formulas in variables from $u$? \end{problem}
  
 Part of the above statement is {\em existence} of such an independent system.  This should be possible essentially using the result over a model $M$
 (\thmref{overamodel}), but replacing $M$ by a probability space of possible interactions with the variables; this only provides a highly 'almost everywhere' result.      
 
 By \thmref{stationarity}, at least a  measure stationarity is obtained assuming higher amalgamation.  For stable theories, the 
 expansion required to obtain higher amalgamation is understood, see  \cite{cigha}.   Could this be combined with stability of the measure algebras  so as to give a more precise construction bringing out the geometry?

 \ssec{Stability and NIP}

The following proposition - for stability and NIP - is a very special case of  a powerful general theory 
of {\em randomization}, due to  Ben Yaacov
and Keisler.      
 The proof we give for all three is a simple application of the Vapnik-Cervonenkis  uniform law of large numbers.
 
 \begin{prop}  \label{pl} Stability, NIP and pNIP are preserved by probability quantifiers: assume  $\psi(u,x;y)$ is stable (respectively NIP,pNIP).  Then
   $(E_u)\psi(u,x,y)$ is stable (NIP,pNIP). 
  \end{prop} 

 \prf  Suppose   $(E_u)\psi(u,x,y)$ is unstable.  Then there exist $\alpha < \beta \in \Rr$ and
  $(a_i,b_i)$ ($i \in \Nn$)  such that $(E_u)\psi(u,a_i,b_j) <\alpha$ when $i<j$ while $(E_u)\psi(u,a_i,b_j) > \beta$
  when $i>j$.  By  
  \cite{vc}, for some $N$ there exist $c_1,\ldots,c_N$ such that for any $a,b$, 
  \[ | (E_u)\psi(u,a,b) -  \frac{1}{N}\sum_{k=1}^N \psi(c_k,a,b)  | < \frac{\beta-\alpha}{3} \]
  Let $\alpha' = \alpha + \frac{\beta-\alpha}{3}$ and $\beta' = \beta - \frac{\beta-\alpha}{3}$.
  By refining the sequence we may assume $\lim_{j \to \infty} \lim_{i \to \infty} \psi(c_k,a_i,b_j) = \gamma_k$
  and $\lim_{i \to \infty} \lim_{j \to \infty} \psi(c_k,a_i,b_j) = \gamma'_k$ both exist.  
Now for $i<j$ we have $\frac{1}{N}\sum_{k=1}^N \psi(c_k,a_i,b_j) < \alpha'   $ while for $i>j$ rather
$\frac{1}{N}\sum_{k=1}^N \psi(c_k,a_i,b_j) > \beta' $.  Thus $\frac{1}{N}\sum_{k=1}^N \gamma'_k < \frac{1}{N}\sum_{k=1}^N \gamma_k $.  But by   stability of $\psi(c_k,x,y)$ we have $\gamma_k=\gamma'_k$; a contradiction.

A similar proof works for NIP, once we know that the value of the sample size  $N$ in the Vapnik-Cervonenkis 
theorem can be bounded polynomially;  in the case of pNIP, we need the bound to depend polynomially on both
the desired approximation and on the pNIP degree (equivalently, on the Vapnik-Cervonenkis dimension.)

To simplify notation take the special case of a $\{0,1\}$-valued relation $\psi(u,x,y)$.  
Let $d'$ be the Vapnik-Cervonenkis dimension of $\psi(u;x,y) \& \psi(u;x',y')$ viewed as a relation between
$u$ and $x,y,x',y'$.   Let $d$ be the Vapnik-Cervonenkis dimension of $\psi(u,x,y)$ 
viewed as a relation between
$x$ and $u,y$.

  Let $n,m \in \Nn$; let $B$ be a set of size $m$ (in the $y$-sort), and let $a_1,\ldots,a_r$ 
have distinct $(Eu)\psi(u;x,y)$ types over $A$, at resolution
$1/n$; in other words if $i \neq j$ then for some $b \in B$, $|(Eu)\psi(u;a_i,b) - (Eu) \psi(u;a_j,b) | > 1/n$.
In particular (possibly after interchanging $i,j$) we have $\mu (\{u: \psi(u;a_i,b) \& \neg \psi(u; a_j,b) \}) > 1/2n$.  
We have to bound $r$ polynomially in $m,n$.   
 
By \lemref{evc}, there exist a set $C$ of size $ \leq (16d'n)^2$ such that for any $i \neq j$, for some $b \in B$ and some $c \in C$
we have $\psi(c;a_i,b) \& \neg \psi(c; a_j,b)$ (or vice versa).  Thus the elements $a_i$
have distinct $\psi$-types over $B \union C$.  By assumption, if $|B| \geq d$ we have 
$r \leq (|B|+|C|)^d $; this gives the required polynomial bound.  
\eprf

 From the above (either using Hoover's normal form, or induction on complexity of the formula) we obtain:

\begin{cor}  Let $L$ be a   language, possibly of continuous logic.  Let $M$ be an $L$-structure.   Assume each basic formula   is stable under $Th(M)$, with respect to any partition of the variables into two nonempty sets.   Then every 
pure probability logic formula is stable.
\end{cor}

A metric space is said to have {\em finite packing dimension} if for some $C,\alpha>0$, for all sufficiently large $n$,
any set of disjoint balls of radius $1/n$ has size at most $C n^\alpha$.    The following is Theorem 4.1 (c) of
\cite{lovasz-szegedy}.

\begin{prop}[Lovasz-Szegedy] \label{ls4.6} Let $\phi(x,y)$ be a $\{0,1\}$-valued NIP formula on $X \times Y$.
 
  Assume given an invariant, generically stable measure $\mu(y)$, with associated 
   expectation operator
$E_y$.     Define a pre-metric $d$ 
on $X$ by 
\[ d(a,b) = E_y (|\phi(a,y)- \phi(b,y)|)  = \mu( \phi(a,y) \triangle \phi(b,y))\]
Let  $M$ be a model, and $\bar{M}$ the completion of $X(M)$.  Then $\bar{M}$ has finite packing dimension,
depending only on the Vapnik-Cervonenkis dimension of $\phi$. 

The same is true for the $L^2$-distance  $d_2(a,b)=\mu( \phi(a,y) \triangle \phi(b,y))^{1/2}$.
\end{prop}

\prf

Let $\delta$ be the Vapnik-Cervonenkis dimension of $\phi$.  By the Sauer-Shelah lemma, 
the number of $\phi$-types over an $N$-element set is bounded by $O(N^{\delta})$.  
Assume the   $1/n$-balls around $a_1,\ldots,a_k$ are disjoint.
We have to bound $k$ polynomially in $n$.  
For $i \neq j$ we have $d(a_i,a_j) \geq 1/n$, so the measure of either $\phi(x,a_i) \& \neg \phi(x,a_j)$
or the dual set, is $\geq 1/2n$.  Let $N= 16 \d n \log(16 \d n)$ and let $ c_1,\ldots,c_N$ be as in 
\lemref{evc}.  Then for each $i \neq j$ for some $\nu \leq N$ we have $\phi(a_i,c_\nu) \& \neg \phi(a_j,c_\nu)$
or vice versa.  Thus the $a_i$ have distinct $\phi$-types over $c_1,\ldots,c_N$.  The number of such types
is at most $O( N^\d)$.  So $k \leq O(N^\d)  \leq O((n \log(n))  ^{ \d})$.

If we use $d_2$ then  $d_2(a_i,a_j) \geq 1/n$ implies $d(a_i,a_j) \geq 1/n^2$, so the same argument gives
$k \leq O((n^2)^{2\d}) = O(n^{4 \d})$.   \eprf

\begin{rem} \label{ls4.6b}  \rm The proof of \propref{ls4.6} is valid for any definable measure $\mu$ on $\phi$-types, using
 \remref{evc2}.   Since the measure is definable, it suffices to consider $a_1,\ldots,a_k$ in a model $M$; while $c_1,\ldots,c_N$
 may be taken in an elementary extension $M^*$.    
 
\end{rem}

\ssec{A categoricity theorem, follwing  Gromov, Vershik, Keisler} 
  
We now formulate a uniqueness theorem for probability logic structures carrying a metric and a definable measure of full support.  
For compact measure spaces, this is a theorem of Gromov's;    Vershik gave a simpler  proof, \cite{vershik}.  All compact structures have  categorical 
 continuous logic theories; the point here is that only   expectation quantifiers are used.

The result also bears a close relation with the 
 uniqueness  theorems for pseudo-finite structures
of Keisler (\cite{keisler77} p. 34, \cite{keisler} 3.2.9; but note the strong property B4 assumed there, and not necessarily valid in our setting, e.g. for the random graph.  (It is valid however when the measure is on the model itself, as is the case in Gromov's theorem.)

In our application to approximately homogeneous approximate equivalence relations,  
the theory itself ensures full support, i.e. when the volume of a ball of a given radius $r>0$ is bounded above $0$.

We   prove the theorem without   uniform full support, compactness or $\si$-additivity  assumptions.    In this case the result may be thought of as  a probability logic analogue to uniqueness theorems for prime models, rather than a categoricity theorem.   Note that it gives in particular a 'soft'   proof of the Gromov-Vershik theorem, different from Vershik's, using a basic model theoretic `preservation theorem':   if the universal theory of $M$ contains that of $N$, then $M$ embeds into an elementary extension of $N$.
 
 Let $L$ be a continuous logic language; it has in particular a formula $d(x,y)$ for a metric, and various additional real-valued relations, uniformly continuous with respect to the metric.      Adjoin expectation operators, and 
let $T$ be a  pure probability   logic theory of $L$;  thus we have a class $C$ of formulas $\phi$ including all quantifier-free formulas, 
and closed under expectation operators.   We say $T$ is ppl-complete if for every $\phi(x) \in C$, $T$ determines $(Ex)\phi$.  

Let $M \models T$.  Recall that the expectation quantifiers induce a measure
on the type space $S_x(M)$, so that any $\bigwedge$-definable set over $M$  is assigned a measure.  
$M$ is said to have {\em full support} if the measure of any ball is positive.     $M$ is {\em complete} if it is complete as a metric space.
 
\begin{thm} 
 \label{gromov} 
Let $T$ be a complete theory of pure probability logic.    If $M,N$ are two complete models of $T$ with full support,
then $M \cong N$. 

  Moreover, any two  tuples in $M$ with the same pure probability
  logic type are conjugate by an automorphism of $M$.
  \end{thm}
  
 \prf   Let us view every formula using connectives and  expectation quantifiers (the class $C$ above) as basic. 
   Write  $M \leq N$ to mean  that any basic formula $\theta$
  satisfies $\theta^{N} (a) =\theta^{M}(a)$ whenever $a \in M^n$.    
   
\claim{1}  Let  ${M_1},{M_2}$ be two complete models of $T$ with full support.  Then the universal theories  of ${M_1},{M_2}$ are equal.  
\prf  Let $\phi(x)$ be a basic formula of $L$,    where $x=x_1,\ldots,x_n$.   It suffices to show that
  if $\phi(a)^{M_1} \geq 0$ for all $a \in {M_1}^n$, then  $\phi(b)^{M_2} \geq 0$ for all $b \in {M_2}^n$.  Suppose for contradiction that $\phi(b)^{M_2} < 0$.  By (uniform) continuity, for some $\e>0$,
  for any $b' \in \Pi_{i=1}^n  B_\e(b_i)$, we have $\phi(b') <0$.  By the full support assumption the measure of each of the balls $B_\e(b_i)$ is nonzero; thus the same is true of their product.  Let $\psi = \min(0,\phi)$.  Then  $E_x \psi^{M_2} < 0$.  But clearly $E_x \psi^{M_1} \geq 0$, contradicting the assumption that the pure probability theories are the same. \eprf
  
  \claim{2}  Let $M \leq N$ with $M$ complete, and $c \in N \m M$.  Then for some $\e>0$,  the ball $B_{\e}(c)$ is disjoint from $M$.  
   \prf   If there were no such $\e$,
  we could find a sequence of elements of $M$ approaching $c$; but $M$ is complete, so $c \in M$ would follow.    \eprf
   
     \claim{3}  Let $M \leq N$.   Let $B=B_{\e}(c)$ be a ball in $N$ with no points in $M$.   Then $\mu(B)=0$.
  \prf      Let $\e' = \mu(B)$.   
    Find in $M$ elements $a_1,\ldots,a_k$ such that 
   $\beta:=\mu(\union_{i=1}^{k} B_\e(a_i))$ is as large as possible, to within $\e'$; so  that  the union of $k+1$ $\e$-balls of $M$ has volume
   $< \e' + \beta$.  By Claim 1, the same is true in $N$.  However  $\mu( \union_{i=1}^{k} B_\e(a_i) \union B_\e(c)) = \beta+\e'$,
   a contradiction.  \eprf

   Let us now prove the theorem.   By Claim 1, $M,N$ have the same universal theory; so $N$ embeds into an elementary extension $M^*$ of $M$;
   we view it as so embedded.   Let $c \in N$.   
   By Claim 2, if $c \notin M$ then some ball  $B_{\e}(c)$ is disjoint from $M$; and by Claim 3, $\mu(B_{\e}(c)) = 0$.  But this contradicts
   the full support assumption on $N$.  Thus $N \subseteq M$.  Similarly, $M \subseteq N$, so $M=N$ and in particular $M \cong N$. 
   
   For the `moreover',   if $a',a''$ have the same type,    enrich $M$ by additional real-valued 
  relations $\phi(x,a')$ (respectively $\phi(x,a'')$), for  $\phi$ a  probability  logic formula,  to obtain structures $M',M''$
  with the same pure probability logic theory, and with full support.  By the main part of the theorem, there exists
  an isomorphism $M' \to M''$, hence an automorphism of $M$ with $a' \mapsto a''$.
 
  \eprf
  
 \begin{rem}  The statement and proof of \thmref{gromov} remain valid for many-sorted theories.  Each sort is assumed to be endowed with a metric,
 and with expectation quantifiers; $M$ and $N$ are assumed to be complete and of full support in each sort separately.
 \end{rem}
  
    \ssec{Local probability logic}
We will require a slight variant, {\em local probability logic}.

We work with local continuous real-valued logic as in \secref{locality}.   Recall that a  local relation
$\phi(x_1,\ldots,x_n)$  has bounded support, determined by  $\rho^*$ and some compactly supported continuous function $b=b_\phi: \Rr \to \Rr$; we guarantee that 
\[ |\phi( x_1,\ldots,x_n)| \leq b_{\phi}(\max_{i,j} \rho^*(x_i,x_j)).\] 

Our description will depend in addition on a choice of positive reals $C_1 \leq C_2 \leq \ldots$; $C_k$ should be thought of as a bound for the measure of a
$\rho^*$-ball of radius $2k$.  

Given a formula $\phi(x,y_1,\ldots,y_n)$ with $n+1$ variables,  $n \geq 1$, we allow an expectation quantifier, so that we can form $(E x) \phi(x,y_1,\ldots,y_n)$.    

By a {\em probability} or {\em expectation quantifier in $x$} we mean a syntactical operation from formulas $\phi(x,y)$ (with $y$ a nonempty sequence of   variables distinct from $x$) to formulas $E_x \phi$ in the variables $y$, satisfying (2-4) of \secref{problogic}, and this generalization of (1):

$(1_{loc,k})$  For any continuous $\beta: \Rr \to [0,1]$ supported on $[-k,k]$, $E_x \beta(\rho^*(x,y)) \leq C_k$.

The numbers $C_k$ are also used in the inductive definition of the syntactic bound for the modulus of continuity of a formula; namely the modulus of $E_x \phi(x,y)$ is $C_{b_{\phi}+1} $ times the modulus of continuity of $\phi(x,y)$.

In practice, we 
concentrate on the case where  the locality relation  is induced by a 
 two-valued relation $R$;  namely $\rho^*=d_R$.      
 \thmref{stabilizer} will  be formulated in this setting (though it could be generalized).      The idea is that 
quantification and expectation can only be taken within $d_R$- balls of some bounded radius.

\ssec{Semantics} 

A model $M$ for local probability logic is a   model  $M$ for the underlying local continuous logic theory, along with a definable measure
on the type space $S_x(M)$, such that for any local formula $\phi(x,y)$, we have $\int \phi(x,b) = (E_x \phi)(b)$.    In particular,  the measure of a $d_R$-ball of radius $k$ is at most $C_k$.  

In locally pseudo-finite semantics, we begin with a family of locally finite graphs $G_i$, letting $\rho^*$ be the graph distance,
and using a multiple $\mu = c_i \mu_{count}$ of the counting measure; such that the volume of a ball of radius $1$ is at most $C_1$.

By {\em pure} (local) probability logic we mean the local formulas obtained from the basic ones using local connectives
(\secref{local-structures}) and expectation quantifiers $E_z$ alone.   
Due to the locality  stipulation in the formation of formulas, there may be no pure probability logic formulas without free variables.  
  We   define the pure probability logic {\em theory} of a structure $M$ by allowing universal quantifiers on the left.   Thus to give this theory is equivalent to 
 determining the closure of $\phi(M)$ for any tuple $\phi$ of formulas.  
 Of course once a constant is added sentences do appear; in the Proposition below, where a constant is assumed, the theory can be taken to be the set of values of {\em sentences}.
  By Claim 1 of \thmref{gromov}, the additional universal quantifiers do not add information here, when $M$ has  full support.

\begin{prop}  \label{gromov2} Let $(X,a)$ be a complete pointed model of a local probability logic theory, such that 
any ball has finite, nonzero measure.      Then the isomorphism type of $(X,a) $  is uniquely determined by the probability
  logic theory of $(X,a)$.  In other words if $(Y,b)$ is another structure with the same properties, and
  the type of $a$ in $X$ equals the type of $b$ in $Y$, then $(X,a) \cong (Y,b)$.   Similarly for $k$-pointed structures.

    \end{prop}

\prf   
Note that  \thmref{gromov} applies to each closed neighborhood of the distinguished point.
 The proof is the same as of \thmref{gromov}, but to begin with choose a   type $q$
  in variables $x_{i,j},  i,j \in \Nn$, such that $x_{k,j} $ lies at distance $\leq k$ from $a$;
with $q$ random in the product space of the balls of radii $1,2,3,\ldots$ around $a$.  
Use also  the additional relations allowed there of the form $\phi(a,x)$, for $\phi$ a  probability  logic formula.

The $k$-pointed case follows from the $1$-pointed case, as the language may include constants.
 Compactness of a metric space implies separability and completeness and so only strengthens the hypothesis.
\eprf

\begin{rem}\rm \label{gromov2r} \rm We will obtain $X$ as the completion of a (locally) saturated probability logic structure $M$, with 
respect to a definable pseudo-metric $d$.   This includes a quotient with respect to the equivalence relation $d(x,y)=0$, which is assumed (locally) co-bounded (this is equivalent to the  (local) compactness assumption on $X$.)   Let $P$ be a $1$-type
of $M$ with respect to pure probability logic, and let $\bar{P}$ be the image of $P$ in $X$.
Then  \propref{gromov2} assures us that the  (isometric) isomorphism group $G$ of $X$ is transitive on $\bar{P}$.
{\em{This does not, in itself, mean  that $Aut(M)$ is transitive on $P$}}, since 
the induced map $Aut(M) \to Aut(X)$ may not be surjective; the pure probability logic type may not generate a complete type.

If we use full continuous logic, including the expectation quantifiers, we can  
enrich $X$ by predicates for all the images on $X$ of $0$-definable relations on $M$.   
They are all closed in the logic topology, and hence in the metric topology.  
  Also  take  $M$ is $|L|^+$-saturated and homogeneous.  In this case, 
 the natural map $G \to Aut(X)$ is surjective.  To see this,  let $\ba$ be a random sequence from $X$ as in the proof of \thmref{gromov}, 
and $\bar{b} =g(a)$.   Lift $\ba$ to $a \in M$.  Then $\bar{b}$ lifts to $b \in M$ satisfying the same type.  By saturation there exists an automorphism of $M$ taking $a$ to $b$.    \propref{gromov2}
is similar. \end{rem}

    \begin{defn} \label{inch}
  A sequence of finite graphs $(\Omega_n,R_n)$  is {\em approximately homogeneous} if   for any pure probability formula in 
  one variable $\phi(x)$, the value of $\phi$ becomes constant as $n \to \infty$:  
  \[ \lim_{n,n' \to \infty}  \sup_{x \in \Omega_n, x' \in \Omega_{n'}} |\phi(x) - \phi(y)|  = 0 \]
  
    A similar definition applies in local probability logic.
  
  The sequence is  {\em approximately homogeneous a.e.} if any pure probability formula in 
  one variable $\phi(x)$, for some $v=v(\phi)$, for any $\e>0$, for all sufficiently large $n$ we have 
  \[ \mu(  \{x \in \Omega_n:   |v - \phi(x) | > \e ) < \e \]
  
In local probability logic, the sequence $\Omega_n$ is approximately homogeneous a.e.  if for  any local pure probability formula in 
  one variable $\phi(x)$, for some $v=v(\phi)$, for any $\e>0$ and $m$,   for    for all sufficiently large $n$ and
  any ball $B$ of $\Omega_n$ of radius $m$, 
    \[ \mu(  \{x \in B   |v - \phi(x) | > \e  \}) < \e \]
   Equivalently, for  any continuous $\beta: \Rr \to [0,1]$ with compact support, 
   for all sufficiently large $n$ we have 
   \[ (E_x \beta(\rho^*(x,t)) |v - \phi(x) |)^{\Omega_n}  <  \e \]
 
       \end{defn}

\begin{rem}\rm   
  Let us  formulate the notion of a sequence
of graphs approaching a $1$-homogeneous graph in probability logic, in terms used in combinatorics (\cite{razborov}, \cite{lovasz-szegedy};
compare also \cite{benjamini-schramm}.).   

    In particular the measure of the set of neighbours $R(a)=\{b: (a,b) \in R\}$ approaches some real number $\varpi$.  
Let $N$ be the set of connected graphs on  $m+1$ vertices.    
 Given $a \in \Omega$,  and $\gamma \in N$, let $C(\gamma,a)$ be the set of graph embeddings $\gamma \to \Omega$ with $0 \mapsto a$.   Define  the {\em  local statistics function} $LS_m : \Omega \to [0,1]^N $
\[LS_m(a) (\gamma) =    \mu_m(C(\gamma,a)) = |C(\gamma,a)|/\varpi^m \]

 Say $(\Omega,R)$ is {\em $(m,\e)$-homogeneous} if the range of $LS_m$ is concentrated in an $\e$-ball (for sup metric on $\Rr^N$.)   If $(\Omega,R)$ and $(\Omega',R')$ are both $(m,\e)$-homogeneous, we say that they
  are $(m,\e)$-close if the respective ranges intersect.  
  \end{rem}

\end{section}

\begin{section}{Stabilizer theorem for approximate equivalence relations}
\label{4-stabilizer}

 Two metrics $d,d'$ are {\em commensurable at scale $\a$} if an $\a$-ball of $d'$ is contained in finitely many $\a$- balls
 of $d$, and vice versa; $k$-commensurable at scale $\alpha$ if the number of balls needed is $\leq k$.
 
 A metric space is {\em $k$-doubling at scale $\alpha$} if $d, (1/2) d$ are $k$-commensurable at scale $\a$.

         \begin{defn}\label{near-def}  Let $\G=(\Omega,R)$, where $R$ is a symmetric, reflexive binary relation.  
 $R$ is a {\em $k$-approximate equivalence relation}  if condition (1) holds.
    It is a {\em near equivalence relation} if for some finitely additive measure $\mu$ on $\Omega$, (2,3) hold.   $R$ is {\em an amenable approximate equivalence relation} if  (1-3) hold.            
    
\begin{enumerate}
 
\item (Main axiom; `doubling').   For all $a$, a 2-ball $R^2(a)$   is a union of  at most $k$ $1$-balls $R(b)$.  
\item   For some $\varpi >0$ and $\kappa>0$, for all $a  \in \Omega$, 
  $(1/\varpi)  \leq \mu(R(a)) \leq   \varpi$; and $\mu(R^3(a)) \leq \kappa$.

\item  (Weak Fubini).  For some $\vartheta>0$, for all $a$,  $\mu(\{b:  \mu ( R( a) \meet R(b) ) \geq \vartheta  \}) > \vartheta$.
\end{enumerate}
  \end{defn}
  
 Given (2), weak Fubini  follows from Fubini, applied to $\{(x,y) :   R(a,x) \wedge R(x,y) \}$, a subset of $R(a) \times R^2(a)$ of   measure at least equal
 to  $1/\varpi^2$.   We note in the lemma below that it automatically  holds (assuming 2) if $R$ is replaced by the distance-two relation $R^2$; 
 or given (1,2),  for   $(\G,R)$ with transitive automorphism group.

\begin{lem} \label{near-lem} \begin{enumerate}
 \item If  $(\Om,R)$ satisfies (2)  then $(\Om,R^2)$ is an amenable approximate equivalence relation. 
 \item If $(\Om,R)$ satisfies (1,2) and has a transitive automorphism group, then  $(\Om,R)$  is an amenable approximate equivalence relation.  
 \end{enumerate}
 \end{lem}
 
 \prf  (1) We use the graph analogue of ``Rusza's trick".   
      Namely,

    a maximal
   disjoint set of balls $R(a_i)$ contained in $R^3(a)$ must have size at most $\kappa \varpi < \infty$.
   By maximality, for any $b \in R^3(a)$ we have $R(b) \meet R(a_i) \neq \emptyset$ for some $i$, so
   $b \in R^2(a_i)$; thus  $R^3(a) \subset \union_i R^2(a_i)$, i.e.  $R^3(a)$
   is a union of at most $\kappa \varpi$ two-balls.   From this, inductively, $R^{2+m}(a)$ is
   the union of $(\kappa \varpi)^m$ two-balls.   In particular taking $m=2$, 
  we obtain \defref{near-def}  (1) for $R^2$.   For (2) we use the same measure; since every two-ball contains a one-balls 
  we have the lower bound on $R^2(a)$; and since every three-ball of $R^2$ is contained in at most $(\kappa \varpi)^4$ two-balls, we have the upper bound on $3$-balls of $R^2$.   Finally to check (3), if $b \in R(a)$ then $R(a) \subset R^2(b)$ so
    $\mu(R^2(a) \meet R^2(b)) \geq \mu(R(a)) \geq (1/\varpi) $.   This shows that $(\Om,R^2)$ is an amenable approximate equivalence relation.

    (2)  Note   that (1,2) imply that (3) holds for {\em some} $a$:     if $R^2(a) = \union_{i=1}^k R(a_i)$, then (3) must hold for some $a_i$, since for any $b \in R(a)$ we
have $\mu(R(b) \meet R^2(a))>0$.   Hence (3) holds for all $a$ if we have homogeneity.

  \eprf

 Note   that subspaces of Euclidean space are doubling at every scale.  
 Let $\Omega$ be an $\e$-sphere packing of  $\Rr^n$, i.e. a maximal set (of 'centers')
 such that any two are at distance at least $\e$.  So any $2\e$-ball contains at least one point of $\Omega$.
 It follows that $(\Omega,R)$ is a $k$-approximate equivalence relation, for appropriate $k$ on the order of $2^n$;
 where $R$ is the 'distance at most  $1$' relation.
 
 If $(\Omega,R)$ is given to us but not $\Rr^n$, can we recover the relations  corresponding to radius
 one half balls, or smaller balls?   \thmref{stabilizer} gives an affirmative result in this direction.  We begin
 with $d_R$ which makes sense for distances $\geq 1$, deduce $\rho$ which is meaningful in distances between
 $0$ and $1$, and show that we still have some doubling, and $\rho, d_R$ more or less fit together at the scale
 $1$.

By definition, a set $Z$ has measure zero iff $B \meet Z$ has measure 0, for all (finite measure) balls $B \subset \Omega$.
 
  Let $R$ be an amenable $k$-approximate equivalence relation on $\Omega$.  We obtain a countably additive measure on
  the type space, or on a sufficiently saturated elementary extension of $(\Om,R)$.  Let $\Phi_0$ be the set of
  pairs $\phi(x),(\alpha,\beta)$ where $\phi$ is a formula of pure probability logic in one variable, and $(\alpha,\beta )$ is a rational interval in $\Rr$, i.e. $\alpha<\beta \in \Qq$; {\em such that}  $S(\phi,\alpha,\beta):=\{x: \alpha < \phi < \beta\}$ has measure $0$.;
  equivalently in terms of expectations, $E_x C(\phi(x)) =0$ for any continuous function $C$ supported on $(\alpha,\beta)$.  Let $\Om_0$
  be the union of all $S(\phi,\alpha,\beta)$, with $(\phi,(\alpha,\beta)) \in \Phi_0$.  This is a countable union, so 
  $\Om_0$ has measure $0$.  Note that allowing $  \phi'(x),(\alpha,\beta)$ where $\phi'$ is a uniform limit of 
  formulas would not change this union.  Let   $\Om ^*$ be the complement of $\Om_0$.  Then
   $\Om ^*$ is a partial type of pure probability logic, has full measure, and is the smallest such.

 \<{thm} \label{stabilizer}
 Let $R$ be an amenable $k$-approximate equivalence relation on $\Omega$.   
 Assume the graph $(\Omega,R)$ is connected.  Then
 there exists a formula $d(x,y)$ of local probability logic, without parameters, such that:
  \begin{enumerate}

\item \label{4.2.2} $d$ defines a pseudo-metric.  
 $d$ and $d_R $ are $k'$ - commensurable at scale $1$, where $k'$ depends only on $k$.  
 
   $d$ is $k''$-doubling at any scale $s \leq 1/2$; where $k''$ depends only on $k$ and $s$.
 
  \item     In the completion of $(\Omega,d)$  (modulo $d(x,y)=0$), all closed balls of radius $\leq 1$ are compact.
The images in the completion of all  $d_R$-balls  any radius are thus  compact.

 \item for any $m \in \Nn$, let   $S_m$   be the distance $\leq 1/m$ - graph of $d$.  
 Then  $S_m^{\cir m} \subset    R^{\cir 4}$, and

 \item   \label{4.2.4}
  for some $C=C_m>0$, for all $a \in \Om ^*$ we have:  $\mu S_m(a) \geq C \mu R(a).$

\end{enumerate}
 
   \>{thm}     \footnote{Without assuming weak Fubini, one obtains the same theorem but with $S_m^{\cir m} \subset R^{\cir 8}$ at worse; it suffices to replace $R$ by $R^2$,  in Claim 1, using \lemref{near-lem}.} 
\prf  
 
(1) 
 
Define $d_0(x,y)$  by the expression:
 \[ E_z  (  | E_t( {R}(t,x) \& {R}(t,z))  -    E_t( {R}(t,y) \& {R}(t,z)) | )   \]
 
 The expectation quantifiers $E_t$ are clearly local.  The quantities  $  E_t( {R}(t,x) \& {R}(t,z)) $ and $E_t( {R}(t,y) \& {R}(t,z)) | $ both vanish
 unless $d_R(x,z) \leq 2$ or $d_R(y,z)\leq 2$; thus the $E_z$ is also a legitimate local probability logic quantifier.
  
 It is clear that $d_0$ defines a pseudometric. 
 
 By \propref{stablei},  $\psi(x,z)= E_t( {R}(t,y) \& {R}(t,z))$ is.a stable formula; we will use this below.   We remark that by the same argument, $d(x,y)$ is stable.

 Let $\vartheta'$ be as in the weak Fubini axiom \defref{near-def}(3), and choose $0<\vartheta<\vartheta'$.

 Define $d =  d_0/\vartheta$. 

\claim{1}  $d(x,y) \leq 1$ implies $d_R(x,y) \leq 4$.    

\prf Pick $z$ with   $E_t( {R}(t,x) \& {R}(t,z)) \geq \vartheta'$.  
If $d_R(x,y) >4$ 
  we have $E_t( {R}(t,y) \& {R}(t,z))=0$  so 
 $E_z  (  | E_t( {R}(t,x) \& {R}(t,z))  -    E_t( {R}(t,y) \& {R}(t,z)) | ) \geq \vartheta'$.
 This shows that $d_0(x,y) < \vartheta'$ implies $d_R(x,y) \leq 4$.   Since $\vartheta$ was slighlty decreased from $\vartheta'$,
  we obtain the stated version with the weak inequality.

\eprf

  So a $1$-ball of $d$ is contained in a $4$-ball of $d_R$, and hence also in finitely many $1$-balls
of $d_R$.

 \claim{2}   For any $s>0$, a $1$-ball of $R$ is covered by a bounded number   of 
 $s$-balls of $d$.  Moreover the bound depends only on $s$ and $k$.
 
 \prf Otherwise, in an ultraproduct, we have a $1$-ball of $R$ containing an infinite $s/2$-discrete
 subset for $d$; in particular there exist   $b,c$ with $d(b,c) \geq s/2$ and
 with the same Kim-Pillay type.

   In the same ultraproduct, consider an element $a$ 
   avoiding
 any given countable set of $b,c$-definable measure zero formulas; in particular it does not divide over $b$ or $c$.  
 By \thmref{uniqueness}   for such $b,c$, for   any such   random $a$ we have 
 \[  E_t( {R}(t,a) \& {R}(t,b)) =    E_t( {R}(t,a) \& {R}(t,c))  \]
   So $d(b,c)=0$, a contradiction.  \eprf

Putting together Claims 1 and 2, we see that $d,d_R$ are commensurable at scale $1$.     Moreover,
 a $1$-ball of $d$ is contained in a $4$-ball of $d_R$ and hence in $k^{3}$ 1-balls of $d_R$; each of these is 
   contained in a bounded number of $s$-balls of $d$, by Claim 2; so certainly a $2s$-ball of $d$ is contained in a bounded number
   of $s$-balls of $d$, if $s \leq 1/2$.      This proves (1).     
  
(2)\,     follows from the total boundedness of the balls of $d$ (Claim 2), and completeness.   
    
(3)\,   Clear

(4)\,        Recall $\Om^*$ from just above the Theorem.

\claim{3}  For any $m >1 $ and any $c \in {\Om^*}$, 
  $\mu S_m(c) > 0$.   
  
  \prf  We may replace $(\Omega,R)$ by an extension saturated for local formulas.  
   Let $c \in \Om^*$, and suppose for contradiction that   $\mu S_m(c) =0$.    
 Let $P=\{c': \mu S_m(c') =0\}$.  Since $c \in \Om^*$, 
  
  for some $R$-ball $B=R^l(a)$, 
$R^l(a) \meet P$ is not contained in a $\bigvee$-definable set of measure zero.  
Let $A=\{a_i: i \in I\}$ be a maximal subset of $P \meet  R^l(a)$ such that the $S_{2m+2}(a_i)$ are disjoint, $a_i \in A$.  Then the
 $d$-balls $S_{m+1}(a_i)$ cover $P \meet R^l(a)$.  
$A$ cannot be infinite; otherwise, by Claim 2, infinitely many elements $a_i$ of $A$ are contained in a single $1/(2m+2)$-ball of $d$,
say around $x_0$; but then $x_0$ lies in each of the $S_{2m+2}(a_i)$, contradicting their disjointness.    Thus
  $A=\{a_1,\ldots,a_\nu\}$ is a finite set, and the $d$-balls $S_{m+1}(a_i)$ cover $P \meet R^l(a)$.  
 
   By saturation,  for some $\e>0$, for all $c' \in R^l(a)$, $\mu S_m(c') < \e$ implies that $c \in S_m(a_i)$ for some $i \leq \nu$.
Since $\mu S_m(c')< \e$ is a  $\bigvee$-definable set, it cannot have measure $0$ in $R^l(a)$; so
some $S_m(a_i)$ has measure $>0$.  This contradicts $a_i \in P$.      
  \eprf

We have shown that   $\mu S_m(c)>0$ for $c \in \Om^*$.  By compactness, 
   $\mu S_m(c) $  must be bounded strictly above $0$, uniformly for all   $c \in {\Om^*}$.

\eprf

\begin{rem}\rm   Fix $\e>0$, and define an $\e$-slice to be a $\bigvee$-definable set $Z$ such that for any ball $R(a)$
we have $\mu(Z \meet R(a)) < \e \mu(R(a))$.  
 
Then for some $C=C_{m,\e}$,  for all $a$ away from an $\e$-slice,  $\mu S_m(a) \geq C \mu(R(a))$.  
 Moreover,  $C_{m,\e}$ depends on $m,\e,k$ alone.  
This follows from \thmref{stabilizer} (4)
by   a compactness argument.  \end{rem}

\begin{question} \rm While $d$ is definable from $R$, almost contained in 
  the distance-two relation  
$R^2$ and  intuitively much finer, it is not clear under what conditions
the image of $R^2$ in the completion of $d$ has the same  probability logic theory  as $(M,R^2)$.   When does the
 image of $R^2$, and similar definable sets, have boundary of measure zero?  
 \end{question}

\begin{rem}\rm    Initially, the  stability theoretic proof produced an $\bigwedge$-definable equivalence relation $S= \meet_m S_m$ corresponding to $d(x,y)=0$, 
  implying the existence of the   relations $d(x,y) \leq 2^{-m}$ with explicitly exhibiting them.  However, an analysis of the proof
   (in the case of ideals arising from a measure)
  shows  that $(x,y) \in S$ iff for almost every  $z$, $\mu(R(x) {\meet} R(z)) = \mu(R(y) {\meet} R(z)) $.  
Thus we can take $S=\meet_n S_n$ with 
  \[    x S_n y  \ \iff  \mu \{z:  | \mu(R(x) {\meet} R(z)) - \mu(R(y) {\meet} R(z)) | \geq 2^{-n} \} \leq 2^{-n} \]
 Putting this into real-valued based probability logic  
 naturally leads to the smoother form used above.  Note that a bounded  $L^1$ function $f$ on a finite probability space has small $L^1$-norm iff it is small in the sense of convergence in measure, i.e. for a small $\e$
 it takes value $>\e$ only on  a set of measure $\e$.  
\end{rem}

 \begin{example}  \rm  Let $G$ be a group, and $X$ an approximate subgroup.   
 Define $R(x,y)$ iff $x \inv y  \in X$.   Then $R$ is an  approximate equivalence relation.
 $G$ acts on $(G,R)$ on the left, by automorphisms.  Since $d(x,y)$  of \thmref{stabilizer}  is definable
 without parameters, $G$ also preserves $d(x,y)$; it follows that $d(x,y)=0 \iff x \inv y \in S$ for some $S$.  
 This recovers the  stabilizer theorem of \cite{nqf}. 
 
  \end{example}

\ssec{Comparison to  Lovasz-Szegedy}   

  After talking about this material in the Groups and Words meeting in Jersualem in 2012, Nati Linial pointed out the relation to \cite{lovasz-szegedy}.
   indeed, their definition
of a {\em graphon} uses precisely the definition of $d$ in \thmref{stabilizer}.  (Strictly speaking, they work in the case of a finite measure, or doubling constant one, whereas ours is only locally finite.)  
 Moreover a more recent paper \cite{lovasz-szegedy2015}  concerns automorphisms of graphons so the overlap must be very considerable.  
 From a model-theoretic viewpoint, the {\em graphon} 
 is isomorphic to a quotient
of the space   of compact Lascar types of a theory of graphs; namely the canonical quotient topologized
by   formulas $\phi(x,b)$ where $\phi$ is the graph edge relation.  In the presence of probability quantifiers,
the KP space, and hence this quotient, carry a canonical measure.  
 
Model-theoretically,  one constructs first a saturated model, then a type space, 
the Lascar compact quotient, and (in the presence of probability quantifiers) measure on it.    The graphon approach, by contrast, begins with the   measure theory; but it is still able to construct parallels  of the above  objects.
The proof that these objects  
are in fact the same,  insofar as pure probability logic goes,  requires  the independence theorem \thmref{ind-prob}.

 \end{section}

\begin{section}{Approximately symmetric approximate equivalence relations, and Riemannian models} \label{riemannian}

\begin{defn} \label{rhs}     
\begin{enumerate}
\item 
A {\em Riemannian homogeneous space} is a connected Riemannian manifold with transitive isometry group. 
\item   
 A metric space  $X$ is {\em locally finite} if    for each point $x \in X$, and any $r>0$, the $r$-ball around $x$ contains only
finitely many points.   A graph $(X,R)$  is locally finite if the associated metric is; equivalently $R(a)$ is finite for all $a \in X$.   It is {\em homogeneous} if the isometry group is transitive.
\ 
 
\end{enumerate}
\end{defn}
 
 We say that  a metric space   is {\em (1-)proper} if each ball (of radius $1$) is compact.    
  A 1-proper metric space is automatically complete.   A connected Riemannian homogeneous space, or a connected (in the graph sense)
  locally finite metric space, are proper and hence complete and separable.  

Riemannian homogeneous  spaces carry a cannonical measure (where balls have finite measure).  Namely
given   $f: U \to \Omega$ where $U$ is an open ball in $\Rr^n$, $\mu(f(U))=\int_U |J(f)| dx$, where $\int_U dx$ is the usual integral, and $J(f)$ is the Jacobian of $f$, i.e. $\det(\partial_i f / \partial_j x_i)$, computed in any orthonormal basis.
Likewise, homogeneous   locally finite metric spaces carry a cannonical measure:   the counting measure, normalized
so that a unit ball has measure $1$.

 Let $X$ be a Riemannian homogeneous space.   The isometry group of $X$ is then a Lie group $L$; and the stabilizer of a point is a compact subgroup of $L$.   See \cite{kobayashi}, Theorem 1.2.

Conversely, let $L$ be a  Lie group  and assume given a transitive action of $L$ on a manifold $X$, with compact point stabilizer $K$.  
Then $X$ can be given a Riemannian structure, such that $L$ acts by isometries.   (Weyl trick:  pick a point $p$; the 
the stabilizer $K$ of  $p$ is compact.  Pick any inner product on the tangent space $T_pX$, and average over $K$ so that we obtain a $K$-invariant inner product.   For any other point $q$, there exists a unique inner product structure on $q$
such that for any $g \in L$ with $g(p)=q$, $g$ induces   an isometry of tangent spaces $T_pX \to T_qX$.)   The invariant  Riemannian structure
is not unique but there is a finite-dimensional space of choices, namely a choice of a $K$-invariant inner product on the tangent space $T_p$.  The tangent space splits as a direct sum of  finitely many invariant subspaces, and
the invariant Riemannian structure is determined up to a scalar renormalization on each of them.    In any case all
these metrics are commensurable, and all commensurable to  any $L$-invariant metric on $X$.

Like Lie groups, Riemannian homogeneous spaces are rather special creatures; they tend to have no deformations, except 
for some freedom in constructing nilpotent ones.   Two canonical Riemannian homogeneous spaces are associated with each Lie group $L$:  the Lie group itself, and the quotient $L/K$ where $K$ is a maximal compact subgroup (which is unique
up to conjugacy.)    In case the stabilizer of a point acts  irreducibly on the tangent space,  they were fully classified by Cartan and Wolf.

This leads to  examples of  homogeneous approximate equivalence relations.

\begin{examples} \label{aer-examples} \rm
 \begin{enumerate}
  \item   Assume $(\Om,R)$ has bounded valency, i.e. $1 \leq |R(a)| \leq k$.
  Then $|R^2(a)| \leq k^2$ so  $(\Om,R)$ is a $k'$-approximate equivalence relation, $k'\leq k^2$.  
   
  \item Assume $|\Om| \leq k^* |R(a)|$ for any $a$.       Then $(\Om,R)$ is a near
 equivalence relation, and $(\Om,R^2)$ is a $k'$-approximate equivalence relation.     Conversely, a $k$-approximate equivalence relation
  {\em of bounded radius $r$} is of this type, $|\Om| \leq  k^* (|R(a)|)$ with $k^*=k^r$.  This is because, 
  inductively,  
  the $l$-ball $R^l(a) = \union_{c \in R^{l-1}(a)} R(c)$ is a union of $\leq k^{l-1}$ $1$ balls.   When $k^*$ is large compared to $k$
  (say $k^*>10k$), this is still an interesting source of examples for us; for instance if $R(a) = [a-1,a+1]$, it is hard to see locally whether
  $\Omega = \Rr$ or $\Omega = \Rr / k^*$ for some very large $k^*$.  
  
  \item Let $(\Omega,d)$ be a Riemannian homogeneous space.   
  Then balls are not finite but they have finite measure with respect to the canonical measure 
  induced by the Riemannian structure.    The metric is $k$-doubling at scale $1$ (and all other scales) for an appropriate $k$, and all $1$-balls have the same measure.    We view this as a {\em model} for $k$-approximate equivalence relations.

  \item  (See \ref{inch}).   \label{aer-examples-4} $(\Omega,d)$ be a Riemannian homogeneous space, and $R=\{(x,y): d(x,y) \leq 1\}$.       Let $(X_i,R_i)$  be a sequence of symmetric 
  binary relations whose pure probability theory approaches that of $(\Omega,d)$.  Then using \lemref{near-lem}, noting that condition (2) of
  \defref{near-def} holds by approximate homogeneity, we see that 
   $(X_i,R_i^2)$ is   an  a.e. approximately homogeneous sequence of $k'$-approximate
  equivalence relations.  We explain below how to obtain such approximations  that are {\em locally finite} (we present this in case $\Omega=L$
  has the form $G(\Rr)$, with $G$ a semi-simple algebraic group over $\Rr$.)   
  \end{enumerate}
  
\ssec{Sprinkling}   
  Let us see how to obtain  sequences of locally finite approximations in Example \ref{aer-examples} (\ref{aer-examples-4}), by choosing them at random (similar procedures are known as  `sprinkling' in   some physics literature.)
   
   Let $L$ be a Lie group.   A {\em lattice} is a discrete subgroup $\Lambda$ of $L$ of finite covolume, i.e. $L/\Lam$ admits a translation invariant finite measure.   Assume   that 
  for any compact $K \subset L$ (with $1 \in K$) there exists a   lattice $\Lambda_K$ with $\Lambda_K \meet K = (1)$.  This is the case for simple Lie groups $L$; for $L=G(\Rr)$ an algebraic group over $\Rr$, where $G$ has no $G_m$-quotients, one can take the arithmetic lattice $G(\Zz)$ or a congruence lattice therein (Borel-Harish Chandra).
  
 Any formula $\phi$ of local  probability logic can quantify only to a bounded  distance, i.e. at $x_1,\ldots,x_k$ expectation quantifiers are applied only
 on $\union_i Kx_i$.  In this case,  if $\Lam$ is a sufficiently small lattice (i.e. $KK \meet \Lam = (1))$, the value of $\phi$ in $L$ and in $L/\Lam$ is the same.  Thus it suffices to show how to find finite approximations $L/\Lam$;   the pullback to $L$ will be an equally good locally finite approximation to $L$.

  This reduces the problem to the finite volume case.   Here  a random choice of $n$ points, for large $n$, provides a good approximation.  
 This is Keisler's relational law of large numbers:  in a probability model $M$ (with the measure on $M$ itself, as is the case for $L/\Lambda$), the initial sections $M_n$ of a random sequence with the normalized counting measure approach $M$ in the sense that for any sentence $\phi$
 of probability logic, $ \phi^M - \phi^{M_n}$ approaches $0$.  See  \cite{keisler77} 6.13, \cite{keisler} 3.1.3.   In case $M$ has a unique probability logic $1$-type, it follows that the sequence $M_n$ is a.e. approximately homogeneous.
 
 Another approach gives a stronger result, at least in the case of an arithmetic lattice, or a co-compact lattice.  In the co-compact case,
$L/\Lam$ along with the metric is easily seen to be interpretable in $\Rr_{an}$.  In the  arithmetic case, \cite{bkt} prove the existence of 
an open, semi-algebraic {\em Siegel set} $A$ for $L$.  Within $A$ one can find a fundamental set $F$ for $L/\Lam$, and any element
of $A$ can be translated to an element of $F$ via a finite subset of $\Lam$.     If $C$ is a compact semi-algebraic neighborhood of the identity in $L$, then $CA$ is contained in finitely many translates of $A$; this implies that the action of $C$ on $L/\Lam$ is also definable semi-algebraically.
The image in $L/\Lam$ of the closed unit ball for the metric may not be semi-algebraic, but  it is definable in $\Rr_{an}$.  Thus our structure
is interpretable in an o-minimal one, and in particular it is NIP.   

Now for a NIP structure $M$, with a measure $\mu$, a fundamental theorem of Vapnik-Cervonenkis (see \cite{vc}, \cite{benyaacov-c-nip}, \propref{evc}) 
 shows that  if $n$ points of are chosen at random to give a subset $M_n$, the law of large numbers holds not just for a given event but uniformly for all definable events.
 Thus up to $\e$-resolution,   $M_n$ look like $M$ not only to an observer outside both, but also to an internal one at point $p$  who takes into account relative positions to $p$, or several such points $p$.  The relational version of this follows inductively (much more readily than in \cite{keisler} 3.1.3, where more careful estimates are needed.)    It is also easy to see here that when $M$ is homogeneous, the $M_n$ will be approximately homogeneous (not only a.e.).

  \end{examples}

 \begin{rem} \rm It is natural to ask for a purely   probabilistic construction, replacing the use of   lattices above.  With our present definition of an 
 amenable approximate equivalence relation, sprinkling  points on $\Om$ using a Poisson process will not work;  while in some sense rare, 
 there would be infinitely many.   
 \end{rem}

\begin{thm} \label{1} Let $(\Omega,R)$ be a local ultraproduct of a sequence of approximately homogeneous, amenable, $k$- approximate  equivalence relations.  Then there exists a  metric  space $(X,d_X)$, and a surjective $\a: \Omega \to X$, both canonically defined,  such that, letting $R_X=\{(x,y) \in X^2:  d_X(x,y) \leq 1 \}$, 
we have:
\begin{enumerate}

\item   The distance between connected components is bounded strictly above $0$; each connected component is clopen, and the space $\Xi$ of connected components is discrete.     The graph induced by $R_X$ on $\Xi$ is   locally finite.

\item   Each connected  component $C$ of $X$ is a Riemannian homogeneous space ; $d_X$ is an invariant metric on $C$.

\item  $R$ is  commensurable with  $\alpha^* R_X = \{(x,x'): R_X( \alpha(x),\alpha(x'))\}$.

\item  For $0 \leq r \leq 1$, the   relation $d_X(\alpha(x),\alpha(y)) \leq r$ is $\bigwedge$-definable on $\Om$.   So is the relation
asserting: \\
  $d_X(\alpha(x),\alpha(y)) \leq r \& \  \a(x),\a(y)$ lie in the same connected component.

\end{enumerate}
 If we assume only a.e. approximate homogeneity, the same result holds but the domain of $\alpha$ is a full measure $\bigwedge$-definable set $\Omega^*$.
\end{thm}

 A few comments, preliminary to the proof.  
 
 1)   The local ultraproduct is taken with respect to the locality relation $d_R$ (\secref{locality}).   Thus $(\Omega,R)$ is a connected
 graph by construction.  In case only a.e. homogeneity is assumed, we choose a component from $\Om^*$.   Thus
 the pure probability logic theory does not depend on this choice.
 
 2)  The theorem describes the situation in scales close to $1$, and refers to $d_X$ only at such scales.  
 To further study the large-scale structure, one needs to combine the geodesic metric on the Riemannian
 manifold with the graph metric of $R_X$.

3)  The proof will   begin with the metric $\rho$ of \thmref{stabilizer}, and the completion $\bar{Y}$ with respect to this metric.
 We would like to find a locally compact group acting on $\bar{Y}$, and connect to the theory of locally compact groups.
The proof would be  simpler if we assumed full first order homogeneity, i.e. a unique $1$-type, or even a unique $1$-type of nonzero measure.   Then we could make use of the automorphism group of the saturated model $\Omega$   and the 
induced action  on $ \bar{Y}$.
  But this would give a result of quite a different nature, applicable only to finite approximations whose full first order theory approaches a given limit.
Full first-order approximate homogeneity is not really a graph-theoretic condition; it concerns not so much the given graphs, but all graphs interpretable within them.    We prefer therefore to assume  only convergence in  the   sense of probability logic.
We must then  accept that even though all elements have the same probability logic type, their full types may differ, and the automorphism group of the ultraproduct $\Omega^*$ may act trivially on the completion.  This obligates us to   work with automorphisms of the completion that are not necessarily induced by automorphisms of $\Omega^*$.  We will use      \thmref{gromov} to obtain  automorphisms of a full measure subset $Y$ of $\bar{Y}$.

\bigskip

   In the proof of \thmref{1}, we will pass from the given approximate equivalence  relation $(\Om,R)$ to the completion with respect to a metric $\rho$, and obtain an induced measure.    For uniformly continuous functions on the completion, expectation quantifiers can be computed 
  either on the completion or on the original structure, giving the same result. This need no longer be the case for discontinuous functions, such as the characteristic function of $R$-balls.    We thus prepare a smoother version of this function.

\begin{lem}   \label{smoothingR}  There exists $R^*$ definable from $R$ in pure probability logic, such that:

 \begin{enumerate}
 \item    $R^*$ is  uniformly continuous    with respect to the metric $\rho$.
  \item For some $\b^2>0$, if  $R(a,b)$ holds, then $R^*(a,b) \geq \beta^2$.
  \item If $R^*(a,b) >0$, then  $R^{17}(a,b)$.
 
  \end{enumerate}
\end{lem}

\prf     Let ${\alpha_0}: \Rr \to \Rr^{\geq 0}$ be the continuous function 
  with value $1$ on $(-\infty,\frac{1}{4}]$,  value $0$ on $[\frac{1}{2}, \infty)$, and linear on $[\frac{1}{4}, \frac{1}{2}] $.  
   Let $\beta>0$ be a lower bound on the
   volume of a $\rho$-ball of radius $\frac{1}{4}$.  

Define 
  \[ R^*(x,y) := E_{u,v} ({\alpha_0}(\rho(x,u)) \cdot R^9 (u,v) \cdot {\alpha_0}(\rho(y,u))) \]

So $R^*$ is definable from $R$ in pure probability logic.

(1)   Uniform continuity is clear since ${\alpha_0}$ is uniformly continuous, and $\rho$ is a metric; so if  $\rho(x,x') \leq \e$ then $| \rho(x,u)-\rho(x',u) |\leq \e$,
 
and similarly for $\rho(y,v)$.

(2) Assume $R(a,b)$ holds.    Whenever $\rho(a,u) \leq 1$ and $\rho(b,v) \leq 1$, we have
$R^4(a,u)$ and $R^4(b,v)$  by \thmref{stabilizer} (3),   so $R^9(y,v)$.  When $\rho(a,u),\rho(b,v) \leq 1/4$ we have ${\alpha_0}(\rho(x,u))={\alpha_0}(\rho(y,v))=1$
so the expectation in the definition of $R^*$ is at least the volume of the product of the balls $\rho(a,x) \leq 1/4, \rho(b,y) \leq 1/4$.

(3)  
If $R^*(a,b) >0$, then for some $u,v$ we have $R^9(u,v)$ and ${\alpha_0}(\rho(a,u)), {\alpha_0}(\rho(b,v))>0$, so $\rho(a,u) < 1/2$, $\rho(b,v) < 1/2$; thus again  $R^4(a,u)$ and $R^4(b,v)$ so $R^{17}(a,b)$.  
\eprf

\begin{proof}[Proof of \thmref{1}]

By definition of the local ultraproduct, $(\Omega,R)$ is connected as a graph.  
Let $\rho$ be the metric of \thmref{stabilizer} and let ${(\bar{Y},\bar{\rho})}$  be the completion.

 By \thmref{stabilizer} (2),  $\bar{Y}$ is $1$-proper.  
We have
a surjective map $h: \Omega \to {\bar{Y}}$, such that the pullback of a closed bounded subset of ${\bar{Y}}^n$ is a
bounded $\bigwedge$-definable subset of $\Omega^n$, in pure probability logic, with parameters.  Indeed $\rho$ is definable in pure probability logic;
if we view $\rho$ as quantifier-free, then 
 for closed $Z \subset \bar{Y}^n$, the pulback $h \inv(Z)$ is quantifier-free definable
with parameters    in $(\Omega, \rho)$ alone (see \secref{logictop}.)   

Let $R^*$ be as in \lemref{smoothingR}.   We consider the structure $(\Omega, \rho,R^*)$; it is a reduct (generally a proper reduct)
of $(\Omega, R)$.   For this reduct, $\rho$ can serve as a metric, since both $\rho$ and $R^*$ are uniformly continuous with respect to $\rho$.
By   the discussion in  \ref{completion}, a structure  $(\bar{Y}, \bar{\rho},\bar{R^*})$ is induced; and further we have local expectation quantifiers
on this structure, and can speak of the expectation of a pure probability logic definable set  (so that a measure on the local type spaces is induced.)  

\claim{0} There exists a smallest  pure probability logic $\bigwedge$-definable  subset $\Omega^*$ of $\Omega$ of full measure;
it determines a unique $1$-type of pure probability logic.   

\prf   The  language is countable,  and has countably many formulas $\phi$ in one variable.
(We can take $\phi$ to be in Hoover normal form.)  Let $v=v(\phi)$ be the generic value in the sense of  \defref{inch}. 
Then $\phi(x)=v(\phi)$ has full measure by the a.e. almost homogeneity assumption.   
Let $\Omega^*$ be the intersection of $\phi(x) = v(\phi)$ over all $\phi$.   Then $\Omega^*$  has full measure;
and the value of any pure probability logic formula is determined on $\Omega^*$.

\eprf

  Let $Y$ be the image of $\Omega^*$ in $\bar{Y}$.   Then $Y$ is closed (this can be checked   within a given small closed ball; and the topology there is the logic topology); and $\bar{Y} \setminus Y$ has measure $0$.
 
Now $\rho$
is definable from $R$ in local probability logic, so any two elements of $\Omega^*$ have the same  pure local probability logic type in the language including  $\rho$.  It follows that the same is true for their images in $Y$, in a language including $\bar{\rho}$ and $R^*$,
using the fact that the pullback of $\bar{\rho}$ is $\rho$, that expectations computed  in $Y$ and in $\bar{Y}$ are the same (as $\bar{Y} \setminus Y$ has  measure $0$), and
that the expectation of a function on $\bar{Y}$ equals the expectation of the pullback  on $\Omega$, by definition of the measure on $\Omega$
(\ref{completion}).  

  In case the sequence is approximately homogeneous, we have $\Omega=\Omega^*$ and so $Y = \bar{Y}$,
since there is only one pure local probability logic type in $\Omega$.  

Let $\bar{R}$ be the image of $R$ in $\bar{Y}$, and let $R_1$ be the image of $R \meet (\Omega^*)^2$ in $Y$.  Form    
  the metrics $d_{R_1}$ on $Y$ and $d_{\bar{R}}$ on $\bar{Y}$.   
The pullback of a finite $d_{\bar{R}}$-ball $B_r(\bar{a})$ to $\Omega^*$ is contained in $B_{5r}(a)$ using \thmref{stabilizer} (3), which controls the 
thickening caused by taking the $\rho$-completion and pulling back.   Thus all $d_{\bar{R}}$-balls in $\bar{Y}$ have finite measure.

\claim{1}$Y$ is locally compact; all $\bar{R}$-balls have compact closure.

\prf   Local compactness of $\bar{Y}$ follows from  \thmref{stabilizer} (\ref{4.2.2}).  
Since $Y$ is a closed subset of $\bar{Y}$,  it is locally compact too.   Now it suffices to show that a $\bar{R}$-ball $b$ is totally bounded
with respect to $\bar{\rho}$; i.e. that for any $\e>0$, $b$ may be covered by finitely many $\e$-balls for $\bar{\rho}$.  For this it suffices to see that a maximal
set of $\e/2$- $\bar{\rho}$- balls is finite.  This in turn follows from the finiteness of $\mu(b)$, and \thmref{stabilizer} (\ref{4.2.4})(along with the fact that $\mu(\Om^* \setminus \Om)=0.)$ 
\eprf

By local compactness the local $n$-type spaces
can be identified with $\bar{Y}^n$ itself.  Hence we  obtain an induced measure on $\bar{Y}$.

Let $G = Aut(Y,\bar{\rho},R^*)$ be  the group of isometries of $(Y,\bar{\rho})$ that preserve $R^*$.  
     $G$ is transitive on ${Y}$, by \propref{gromov2}.

Define a topology on $G$ as in \secref{grouptop}.

\claim{2}  $G$  is locally compact.     for $a \in Y$, the point stabilizer  $G_a = \{g: g a =a \}$ is compact;   for $a \in Y$, the natural map $G/G_a \to Y$ is a homeomorphism.  
  
\prf By \lemref{grouptoplem}, taking the relation $\mathsf{R}$ there to be $\{(x,y):  R^*(a,b) \geq \beta^2\}$, and $d=\rho$. 
 (The assumptions of \lemref{grouptoplem} hold by Claim 1 and  \lemref{smoothingR} (2,3).)    \eprf

  According to Gleason-Yamabe
 \cite{yamabe}, $G$ has an open subgroup $H$, and a compact normal subgroup $N$ of $H$ (contained in any
  desired neighborhood of $1$), such that $H/N$ is a finite dimensional Lie group.

 Dually to \cite{nqf},   let  ${\b_0}$ be the set of pairs $(H,N)$ with $H$ an open subgroup of $G$,
 and $N$ a compact normal subgroup of $H$.   
Let $\b$ be the set of pairs $(H, N) \in {\b_0}$ such that if $(H',N') \in {\b_0}$ and $N \leq N' \leq H' \leq H$ then 
$H=H'$ and ${N}={N}'$.  Equivalently, $(H, N) \in {\b}$ iff the locally compact group $H/N$ is connected, with no nontrivial compact normal subgroups.   (Hence by Yamabe, is a Lie group.)   Non-emptyness of $\beta$ follows
from a chain condition on closed subgroups shown in the second paragraph of \cite{nqf} 4.1.    
It is also shown there that if $(H,{N})$ and $(H',{N}') \in \beta$ then $H \meet H'$ has finite index in $H$ and in $H'$,
and $(H \meet H', H \meet {N}') \in \beta$.   And $H$ determines $N$ uniquely,
i.e.    for any $H$, there is at most one $N$ 
with $(H,N) \in \beta$.

Fix $a \in Y$.  Since $G_a $  is compact while $H$ is open, $G_a \meet H$ has finite index in $G_a$;
in particular there are only finitely many $G_a$-conjugates of $H$.  Taking their intersection,
we see that there exists $(H_a,{N}_a)  \in \beta$ normalized by $G_a$.

Let $(H_b, {N}_b) = (g  H_a g \inv, g  {N}_a g \inv) $ for any $g$ with $g(a)=b$.  Define an equivalence relation $E$ on ${Y}$ 
by $bEb'$ iff $H_b=H_{b'}$ and $H_b b = H_{b'} b'$.  Then $E$ is $G$-invariant.  If $a E a'$, then $a' = h a$ for some $h \in H_a$.
Conversely if $a' = h a$ with $h \in H_a$, then  $H_{a'} = h H_a h \inv = H_a$; and $H_a a = H_a a' = H_{a'} a'$.  Thus the $E$-class of $a$ is just $H_a a$; it is the image of $H_a$ under the natural map $G \to Y$, $g \mapsto ga$.  We saw in Claim 2 that this map induces a homeomorphism $G/G_a \to Y$.   Thus the $E$-class of $a$ is   open;
by transitivity each $E$-class is open, and hence each $E$ -class is clopen.

Give  $Y/E$ a graph structure via ${R_1}$.  

\claim{3}    $Y/E$ has finite valency.  The  topology on $Y/E$ induced from $Y$ 
is   discrete.    
\prf  Since a $d_{\bar{R}}$-ball $b$ is compact while each $E$-class is open, $b$ meets only finitely many $E$-classes.
\eprf

 Since the $E$-classes are open, there exists $r_0>0$ such that the $r_0$-ball around $a$  
 is contained in the $E$-class of $a$; by transitivity of $G$, the same holds for all points.    
 Thus two points in distinct $E$-classes are at $\rho$- distance $\geq r_0$.    We may choose $r_0 \leq 1/2$, so that
 the pullback to $\Om$ of a $\bar{\rho}$ of radius $r_0$ is contained in finitely many $d_R$-balls of radius $1$ (\thmref{stabilizer} (1)).

 Define a finer equivalence relation $\be$ on ${Y}$ by $b \be b'$ iff $b Eb'$ and ${N}_b b = {N}_{b'} b'$.  
 (Note that $b E b'$ implies $H_b=H_{b'}$ and hence $N_b=N_{b'}$.)
 Then for each $E$-class $Y' \subset Y$,
 $Y'/\be $ is connected, since a conjugate of the connected group $H_a/{N}_a$ acts transitively.  
 
 Let $X = Y /\be$, and $\alpha: \Omega^* \to X$
 the composition $\Om^* \to Y \to X$.  
 
Define \[d_X(x,y) := \inf \{ \min(r_0,\bar{\rho}(a,b)):   a,b \in Y, \alpha(a)=x, \alpha(b)=y \} \]

Since the classes of $\be$ are compact, the infimum in this definition is attained, and thus $d_X(x,y)=0$ implies $x=y$.
It is clear that $d_X(x,x)=0$ and that $d_X$ is symmetric.  Let us consider the triangle inequality $d(x,z) \leq d(x,y) + d(y,z)$.
If $d(x,y) \geq r_0$ or $d(y,z) \geq r_0$ the inequality is clear; so we may assume $d(x,y), d(y,z) <r_0$.  Let $(a,b)$
attain the minimum in the definition of $d_X(x,y)$, and likewise $b',c$ for $y,z$.  Then $\bar{\rho}(a,b)< r_0$ so $aEb$, and likewise $b'Ec$.
We also have $b \be b'$ since $\alpha(b)=y=\alpha(b')$.  Thus $b'=n b$ for some $n \in N_b =N_a=N_c$.   So 
$\bar{\rho}(b',c) = \bar{\rho}(b, c')$ where $n c' = c$.   It follows that $d_X(x,z) \leq \bar{\rho}(a,c') \leq \bar{\rho}(a,b)+ \bar{\rho}(b,c') =
d_X(x,y)+d_X(y,z)$.    Hence $d_X$ is a metric on $X$.  It induces the same topology on $X$ as the quotient topology inherited from $Y$ via
$X=Y/\be$.

The pullback of an $r_0/2$- ball $a$ of $X$ to $Y$ is the union over a $\be$-class $b$ of the $r_0/2$-$\rho$-balls around elements of $b$.
Since $b$ is compact, finitely many of these $r_0/2$-balls cover $b$, so that the union  is contained in a finite union of $ r_0$-$\rho$-balls.  
Hence $R$ is commensurable with $\alpha^* R_X$, where $R_X =\{(x,x') \in X^2:   d_X(x,x') \leq r_0/2 \}$.  As all   the clauses of
\thmref{1} except for (3) are invariant under rescaling, we may replace $d_X$ by $2 d_X/r_0$, so as to obtain (3).  

  The connected components of $X$ are the images of the classes of $E$, and are 
 at distance at least $r_0$ from each other.  On each connected component, we have a transitive
 action of a Lie group $\mathsf{L} \cong H_a/N_a$, respecting the metric.  
 The induced metric on the space of connected components is locally finite by Claim (3).

 Let us now address the definability issues  (4).  Recall the map $h: \Omega \to \bar{Y}$ from the first lines of the proof.  
The equivalence relation $h(x)=h(y)$
is $\bigwedge$-definable,  by the pure probability  formula $\rho(x,y)=0$, and more generally $\bar{\rho}(hx,hy) = \rho(x,y)$.  Thus $(\bar{Y},Y,\bar{\rho})$
can be viewed as interpretable in $(\Om,R)$.  

Let $P \subset Y^k$ be an automorphism-invariant closed relation on $Y$.   
Pick   $a  \in Y$, and consider a bounded $\bar{\rho}$-ball $B$ around $a$, say of radius $r_1 \leq 1$.  Then $B$ is compact;
so the restriction $P|B = P \meet B^k$ is $a$-definable; say via a formula $\psi_P(x,a)$ of continuous logic.  Since $G$ acts transitively on $Y$,  $\psi_P(x,y)$ defines $P$ on $Y$,
provided $P$ is local, i.e. $P(x_1,\cdots,x_k)$ implies $\rho(x_i,x_j) \leq r_1$.  In particular this is the case for $P=\be$, since the classes of $\be$ are compact, and for the 
relations in the statement of (4).   
Pulling back to $\Om^*$ we see that $P$ is $\bigwedge$-definable (though not necessarily by a formula of pure probability logic.).  
  
\eprf

\begin{example}  \rm Let $A_n$ be the interval $[-n,n]$, with Lebesgue measure, and let $R(x,y)$ be the relation $|x-y| \leq 1$.   This sequence of $2$-approximate equivalence relations is not strictly approximately homogeneous;
the measure of $\{x: R(a,x) \}$ is $2$ towards the middle, but approaches $1$ near the endpoints.  However, even in full continuous logic with probability quantifiers, it is a.e. approximately homogeneous;
a formula whose quantifiers look out to distance $d$ will take the same value on all points except for the intervals $[-n,d-n]$ and $[n-d,n]$, whose measure $2k/n$ approaches $0$ with $n$.

\end{example}

\ssec{Definability and asymptotic structure}   \label{definability}  The  definability statements (4) in \thmref{1} are crucial for deducing asymptotic consequences from the structure of the limit.  We plan to return to this in the future, and for now only sketch a basic statement
 by way of illustration.  Let $(A_n,R_n)$ be a sequence of increasingly homogeneous locally finite graphs.     Fix $r>0$ and let $B_n$  be a ball of radius $r$ in $A_n$, and $B_\infty$ a ball of radius $r$ in an ultraproduct $A$.  
  By assumption the $B_n$ look increasingly similar as $n \to \infty$.

 Let $E$ be the formula defining the 'same connected component' equivalence relation on $B_\infty$.  
 Then $E$ has a finite number $f$ of classes on $B_\infty$.  For   almost all $n$, the same formula $E$ defines a partition
  of $B_n$ into $f$ classes.   One can further find a modified graph structure $R'_n$ on $B_n$, uniformly commensurable
  to $R_n$.   Let $C_n$ be one of the $f$ classes of $B_n$.     After refinement of the indices $n$,  $(C_n,R'_n)$
    converge to a bounded region on a Riemannian homogeneous space $H_r$,  with 
   graph structure 'distance $\leq 1$' relative to the Riemannian metric $d_{Riemann}$.  One can also find asymptotic versions of the metric $d_{Riemann}$ on the $C_n$.

\ssec{Homogeneity}   \label{grouptriple} \rm   
Consider an approximate equivalence relation $(\Om, R)$ with  transitive automorphism group $G$.   Let
$H$ be the stabilizer of a point $a \in \Om$.   Let $X=\{g \in G:  (a,ga) \in R\}$.  It is easy to see that $X=X \inv$ is an 
approximate subgroup of $G$,  and $HX=X$.   

Conversely,   a triple $(G,X,H)$ with $X$ an approximate subgroup of $G$, $H$ a subgroup with $HX=X$ gives an approximate equivalence relation  $ \R_{G,X,H}$ on $G/H$ with transitive automorphism group.  
Namely, $( g_1H, g_2H) \in \R_{G,X,H}$ iff $g_1 \inv g_2 \in X$.   When $H=1$ we write $\R_{G,X}$ or just $\R_X$.

\begin{problem}[Rigidity]\label{finite-axiomatizability?} 
Let $G$ be a semisimple Lie group $G$, with a maximal compact subgroup $K$.  Let $R(x,y)$ be the relation: $y \inv x \in K$.
Show that the pure probability logic theory of $(G,R)$ determines $G$ uniquely.   Further, does there exist  a single sentence $\si$ 
 of probability logic such that any increasingly homomgeneous sequence of approximate equivalence relations $X_i$
 whose $\si(X_i)$ converges to $\si(G)$ must in fact converge to $G$?   This would say that a resident of $X_i$
(for large enough $i$) can reasonably guess the limit $G$ that the sequence is tending to.   

\ 

\rm Assume the ultraproduct $M$ of \thmref{1}      has an associated  Riemannian symmetric space whose Lie group $L$ is simple, with    a finitely presented lattice $\Lambda$, with generators $\lambda_1,\ldots,\lambda_k$.  Describe  the element $\lambda_i$ of $L$ as a limit of  `rough' symmetries of the approximating graphs, and how to recognize the relations on the $\lambda_i$
from rough relations among these.    When this can be done and  $\Lambda$ is simple,  it should
becomes possible to  recognize it by looking at a single sufficiently good approximation.

  \end{problem}

\begin{problem} \rm \label{mixed}
\thmref{1}  describes completely the connected case, and the locally finite case.   But the    
mixed case is not fully described.   What are the possible homogeneous extensions $(\Omega,R)$ of a 
 homogeneous Riemannian space  $X$ by a homogeneous, locally finite graph $\Xi$?  
 
 We have copies $X_a$ of $X$, for $a \in \Xi$.  Say $G=Aut(X)$ is a  centerless semisimple Lie group.
 In this case,  for two points $a,b$ of $\Xi$ connected by an edge, there exists a unique isometry $X_a \to X_b$ at finite $d_R$-distance.   This gives
a homomorphism $\pi_1(X,a) \to Aut(X_a)$.   Does this fully describe the structure of $(\Omega,R)$, up to commensurability?
\end{problem}

\begin{problem}  \label{non-amenable} \rm  Theorem 5.16 of \cite{H-beyo}  describes the structure of approximate subgroups without an amenability
assumption.  Generalize this and \thmref{1} to homogeneous approximate equivalence relations; see \ref{grouptriple}.   Definability may be challenging as
the definability statement in (\cite{H-beyo}, Theorem 5.16) is weaker than in the amenable case; 
and it refers a priori to the group $G$, a  second-order object from the point of view of $(\Om,R)$.

\end{problem}

\begin{problem} \rm \label{nested}   

  In the model-theoretic limit, the pieces of a (well-chosen) partition of any structure into $n$ pieces tends to a union of {\em homogeneous} structures. 
 Given an arbitrary  approximate equivalence relation, can one partition into pieces that resemble various Riemannian homogeneous spaces?

Here is a more precise formulation:

Let $(M_i)_{i \in \Nn}$ be any family of locally finite, $k$-approximate equivalence relations;  we make no homogeneity assumption.  
Let $M$ be a (sufficiently saturated) ultraproduct;  so  $Aut(M)$ acts transitively on each complete type $P$.
The restriction to $P$ is then a homogeneous $k$-approximate equivalence relation, and a structure theorem should apply
(either using an answer to Problem \ref{non-amenable}, or disintegrating the measure over almost all $P$.).    What does this imply for the given $M_i$?

One might expect a statement of the following type: take any function $\alpha: \Nn \to \Nn$   growing to $\infty$, as slowly as you like.
Then one can find a subsequence $I$ of $\Nn$, and for $i \in I$ a partition of $M_i$ into $\alpha(i)$ definable sets $(D_{i,k}: k < \alpha(i))$, refining
the partitions $D_{i',k'}$ for $i'<i$; such that each branch resembles,
more and more closely, a $k$-approximate equivalence relation commensurable to a Riemannian homogeneous space.  By a branch we mean a choice
of a definable set $D_{i,k(i)}$ for each $i$, so that $D_{i,k(i)}$ implies $D_{i',k(i')}$ if $i'< i \in I$; it gives a sequence of graphs $D_{i,k(i)}(M_i)$.

See \cite{benjamini-hutchcroft} for the $1$-dimensional case.

\end{problem}

\end{section} 
\begin{section}{strong approximation:  from groups to graphs} \label{grtogr}

This section again addresses the structure of approximate equivalence relations; but 
here we assume  definability over finite fields, or more generally the existence of a dimension theory similar to the one available for pseudo-finite fields.    We will see that approximate equivalence relations in this setting are close to actual equivalence relations.

We first recall a key statement of   \cite{HP2}  
in the case of groups , \thmref{hp} below, that forms the model of the graph-theoretic generalization we aim for.  This was the main ingredient in model-theoretic
proofs of {\em strong approximation} over prime fields, for instance of the fact that if $H$ is a Zariski dense 
subgroup of $SL_n(\Zz)$ then $H$ maps surjectively to almost every  $SL_n(\Ff_p)$.

\begin{thm}[strong approximation lemma: groups] \label{hp}Let $F=\Ff_p$,  $p$ nonstandard.  Let
$G$ be a definable group over $F$,  and let  $(X_i:i \in I)$ be a family\footnote{not necessarily uniformly definable}  of definable subsets of $G$.  
Then:  (A)  there exists a definable  $H$ such that:  \begin{enumerate}
\item $H$ is a subgroup of $G$
\item  $H \subset <\union X_i> $.
\item $X_i / H$ is finite.
\end{enumerate}

\, 
If read for standard primes $p$, each $X_i$ should be definable uniformly in $p$; and in (3), each
$X_i / H$ has finite cardinality independent of $p$.
Moreover, if $G$ is an algebraic group: \\

\noindent (B):    there exists a homomorphism $h: \widetilde{H} \to G$ of algebraic groups, with finite kernel 
such that $H = h(\tH(F))$.

\end{thm}

\ssec{}\label{unip}   Here is how strong approximation follows from \thmref{hp}.   Let $G$ be a linear algebraic group. 
Applied to the family of one-dimensional unipotent subgroups $X_i$ of {\em an arbitrary subgroup $\G$ of $G$}, \thmref{hp} shows that $\G$ contains a {\em definable} normal subgroup $H$ (generated by the unipotent elements of $\G$)   such that $\G/H$ has no unipotent elements.   By Jordan's theorem it follows that $\G/H$ is 
Abelian-by-bounded.  In particular, the image of a Zariski dense subgroup of $SL_n(\Zz)$ in 
$SL_n(\Ff_p)$ has bounded index, and admits the description of   in \thmref{hp} (B).  
See \cite{HP2}, Prop. 4.3 and 7.3.   These results were previously proved by other means by Weisfeiler, Nori, Gabber.

\ssec{ Generalization to graphs}
The graph version applies to definable graphs over $F$.  For instance, the graph may be $(V,X)$ with $V$ a  variety and $X$ a subvariety of $V^2$, or more generally the projection of the rational points of a variety $W$ mapping onto $V^2$.  
We try to study them {\em up to graphs of bounded valency}.    
 Essentially we show that  after a  partition of the vertices into a (bounded) finite  number of pieces, and after   fibering each piece over a graph of bounded valency in each piece, $X$ generates an equivalence relation on each piece in a bounded number of steps.    Between any two distinct pieces,
 the induced graph has finite valency in at least one direction, following  an interesting partial ordering of the set of pieces.

While  \thmref{appr-graphs} (B) is formulated for pseudo-finite fields,  part (A) and 
 \propref{appr-graphs-inf}    are valid for structures with a finite, definable S1-rank in the sense of \cite{pac}.  We refer to such structures and their theories as S1-structures and theories.  
 
The reader who wishes can read both parts for pseudo-finite fields.

For a graph $(G,X)$, $X \subset G^2$, we let $\simg$  be the equivalence relation  of connectedness, i.e. $\simg$ is the smallest equivalence relation on $G$ containing $X$.  

 \begin{prop}  \label{appr-graphs-inf} Let $(G,X)$ be a definable graph in an S1-theory.  
 On each type $P$ of $G$ there exists a canonical $\bigwedge$-definable
 equivalence relation $E_P$ such that $E_P \subset \simg$, in fact $E_P \subset X[m]$ for some $m$, and such that the  induced graph on $P/E_P$ is locally finite.

   The definition of $E_P$ depends only on $P$ and on $\simg$.

 \end{prop}

\begin{rem}  \propref{appr-graphs-inf}  is valid more generally if $X= \union_i X_i$ is only $\bigvee$-definable; i.e. the  edges $X$ are given as a countable  union of definable sets   in a countably saturated model.   In this case local finiteness means that any finite union  $X' = \union_{i \in F} X_i$,  
  and any $\bar{a} \in P/E_P$, there are only finitely many $\bar{b} \in P/E_P$ with $(\bar{a},\bar{b}) \in X' / E_P$.  \end{rem}
 
 Here is a more detailed version, that allows passage to the finite.

\begin{thm}  [strong approximation lemma: graphs]   \label{appr-graphs} \, \\
\begin{itemize}
\item[(A)]  Let $(G,X)$ be an $S1$-structure  with $X \subset G^2$ symmetric, reflexive.       Then there exists  
 a definable   partition $G=D_1 \union \cdots \union D_n$     a definable function $f$ on $G$, and $m \in \Nn$
   such that:\begin{enumerate}
 
\item If $f(a)=f(b)$ then $d_X(a,b) \leq m$.  
\item   
For $a \in G$,  $f(X(a))$ is finite. 
\item For
  $i \leq j \leq n$,  let $\bar{X}_{i,j} = \{(f_i(x),f_j(y)):  x \in D_i, y \in D_j, (x,y) \in X \}$.
   
 Then for any $\bar{b} \in f(D_j)$ there are only finitely many $\bar{a} \in f(D_i)$
 with $(\bar{a},\bar{b}) \in \bar{X}_{ij}$.  
 \ \\

\end{enumerate}
\item[(B)] Let $F$ be an ultraproduct of finite fields.   Let $(G,X)$ be an $F$-definable graph:    
   $G$ is a definable  subset of  $V(F)$, $V$ an $F$-variety, and   $X   \subset G ^2$ is definable.  Then (A) holds with an {\em algebraic} $f$  and $D_i$:   we can find a quasi-finite morphism of $F$-varieties ${\rho}: \tG \to V$ with ${\rho}(\tG(F)) =  G(F)$,   and regular functions $\phi$ on $\tG$,  
  such that  for $u \in \tG$,
  \[ f({\rho}(u))= \phi(u)  \]

\end{itemize}
\end{thm}       
 
\

   {\em Preliminaries to proof.} 
We work in a sufficiently saturated and homogeneous model.  Two elements with the same type are then conjugate by some automorphism.  
 We write $d(b/a)$ for the dimension of the smallest $a$-definable Zariski closed set including $b$.
 More generally in the setting of finite S1-rank, $d(b/a)$ is the least S1-rank of a formula $\phi(x,a)$ true of $b$.  .
 We write $\Ind{a}{b}{c}$ if $d(a/b,c) = d(a/c)$.   In the case of pseudo-finite fields, dependence and independence are determined by the quantifier-free type; we have $\Ind{a}{b}{c}$ iff the algebraic locuse of $a/\acl(b,c)$ (namely the variety $V$ of smallest
dimension  defined over $\acl(b,c)$ with $a \in V$)  is already defined over $b$.  
 
 {\em Canonical bases}    Given tuples $b,c$ there exists a unique smallest algebraically closed $\bar{c} \subset \acl(c)$ with 
 $\Ind{b}{c}{\bar{c}}$.   We write:  $\bar{c} = \cb(b/c)$.  If $\acl(c) \subset \acl(c')$ and   $\Ind{b}{c'}{c}$ then $\cb(b/c')=\cb(b/c)$.  
 In the case of pseudo-finite fields,  $\cb(b/c)$ is the field of definition of the locus of $b$ over $c$. 
 
In simple theories of finite SU-rank, one can instead use canonical bases in the sense of simple theories,  (\cite{HKP}).
 
As a consequence of the existence of canonical bases, we see that if   $\Ind{b}{c}{d}$ and  $\Ind{b}{c}{d'}$ where $d,d' \subset \acl(c)$, and  $\acl(d'')=\acl(d) \meet \acl(d')$,
then $\Ind{b}{c}{d''}$.

{\em Elimination of imaginaries:}   $\Ff$ has a unique extension $\Ff_n$ of degree $n$; it is Galois with Galois group $\Zz/ n \Zz$;
when $F$ is finite, $Aut(F_n/F)$ has a canonical element $\phi_n$; we view this as part of the structure of each finite field $F$; correspondingly
we obtain an element $\phi_n$ of $Aut(\Ff_n/\Ff)$, that we view as definable.  This amounts to naming a   certain algebraic  imaginary   element $e_n$
coding $\phi_n$, for each $n$, on top of the field structure.   With this understood, $Th(\Ff)$ admits elimination of imaginaries; see \cite{ChH}.

\begin{proof}[Proof of   \propref{appr-graphs-inf}]

For $u,v \in G$,  write $u \sim v$ if $u,v$ lie in the same connected component of the graph $(G,X)$.  This is a  
countable union of definable relations, $\bigvee_{j \in \Nn} d_X(u,v) \leq j$.

 Fix for a moment an element $b \in G$, with type $P$.  Let $W=W_b$ be the connected component of $b$ within the graph $(G,X)$.  
 Choose $a \in W$   with  $d(b/a)$ as large as possible; let $  m_P = d(b/a)$.     Let $\d_P =d_X(a,b)$ (for definiteness, choose
 $a$ so that $\d_P$ is minimal, subject  to $d(b/a)=m_P$.)  
 Let $e =  \cb(b/a)$.   
 
\claim{1}   For any  $c \in W$, $e \in \acl(c)$.   In particular, $e \in \acl(b)$.

We first prove the claim under the assumption that  $\Ind{c}{b}{a}$.   We  have $d(b/c) \geq d(b/ac)=d(b/a)$.   
By maximality of $d(b/a)$, 
 $d(b/c)=d(b/a)= d(b/ac)$.   So  $ \Ind{b}{a}{c}$.  Hence $e=\cb(b/a)=\cb(b/ac)  \in \acl(c)$. 
 
Now for a general $c$, take $c'$ with  $tp(c/\acl(a))=tp(c'/\acl(a))$ and   $\Ind{c'}{b}{a}$.   Since $a \sim c$ we have also $a \sim c'$,
 so $c' \in W$.    By the special case just proved, we have 
$e \in \acl(c')$.     
 Since $e \in \acl(a)$, we have $tp(c/e)=tp(c'/e)$.  Thus  
  $e \in \acl(c)$.

 \claim{2}  If $tp(b' / \acl(e))=tp(b/\acl(e))$  and $\Ind{b}{b'}{e}$ then 
 then $b' \in W$;   in fact $d_X(b,b') \leq 2\d_P$.  
 
     Recall $e = \cb(b/a)$ so $\Ind{b}{a}{e}$.  
Using the  
 independence theorem over $\acl(e)$ (computed in $M^{eq}$), we can find $a'$ with 
 $tp(a',b / \acl(e) ) = tp(a',b'/  \acl(e) ) = tp(a,b /  \acl(e) )$.    By the definition of $\d_P$ and the choice of $a$  we have $d_X(b',a) = d_X(b,a) = \d_P$.   Since $tp(a',b)=tp(a',b') = tp(a,b)$ we also have $d_X(a',b')=d_X(a',b )=d_X(a,b)=\d_P$.   
 So $d_X(b,b') \leq 2\d_P$.

  \claim{3}  If $tp(b' / \acl(e))=tp(b/\acl(e))$   
 then $b' \in W$;   in fact $d_X(b,b') \leq 4\d_P$.  
 
   Indeed let $tp(b''/\acl(e))=tp(b/acl(e))$
 with $\Ind{b''}{b,b'}{e}$.  By Claim 2
  we have  $d_X(b,b'') \leq 2\d_P$ and  $d_X(b',b'') \leq 2\d_P$; so $d_X(b,b') \leq 4\d_P$.

\claim{4}       $\acl(e)  = \meet_{c \in W}  \acl(c)$.

We already saw one direction in Claim 1.  Conversely in Claim 2 we saw that if $tp(b' / \acl(e))=tp(b/\acl(e))$
and $\Ind{b}{b'}{e}$ 
then $b' \in W$; if   $d \in \meet_{c \in W}  \acl(c)$, then $d \in \acl(b) \meet \acl(b')$ so $d \in \acl(e)$.

\claim{5}  Let $tp(e'/b)=tp(e/b)$; then $\acl(e')=\acl(e)$.  More generally if $Aut(M/b')$ leaves $W$ invariant, 
and $tp(e'/b')=tp(e/b')$, then $\acl(e')=\acl(e)$.

 Indeed let $\si$ be an automorphism fixing $b'$ and with $\si(e)=e'$; then $\si(W)=W$, and by Claim 4 we have
 $\acl(e') = \si(\acl(e)) = \meet_{c \in \si(W)}  \acl(c) = \meet_{c \in W} \acl(c) = \acl(e)$.   
 
 \bigskip 
 
 \claim{6}  There exists a definable function $f_P$ such that $\acl(e) = \acl(f_P(b))$.
 
 Let $\be$ be a code for the finite set $\overline{E}$ of conjugates of $e/b$.  Then $\be \in \dcl(b)$.   
 So $\be=f_P(b)$ for some definable function $f_P$.  
 By Claim 5, each   element  of $\overline{E}$ is equi-algebraic with $e$; so we have also $\acl(e)=\acl(\be)$.  
 
  Thus far, the element $e$ played a role only via $\acl(e)$;   so we may replace $e$ by $\be$ with no loss.

  Define an equivalence relation $E_P$ on $P$, setting $x E_P y$ iff $f_P(x)=f_P(y)$ and
  $tp(x / \acl(f_P(x))) = tp(y / \acl(f_P(x)))$.      It is easy to see that $E_P$ is $\bigwedge$-definable.  
  Both equations together can be written as:  
    $tp(x / \acl(f_P(x))) = tp(y / \acl(f_P(x)))$.   (Since this implies $f_P(x)=f_P(y)$.)  As this equation refers only to $\acl(f_P(x))= \meet_{u \sim x} \acl(u)$, it does not depend on the specific choice
  of the   $0$-definable function $f_P$.

 \claim{7}     If $c,d \in P$ and $c \sim d$ then $f_P(c) \in \acl(f_P(d))$.
 
 \prf  Let $d'$ be independent from $c,d$ over $f_P(c)$, with $tp(d'/\acl(f_P(d))=tp(d/\acl(f_P(d))$.  In particular $d' \in P$, 
 $f_P(d')=f_P(d)$, and 
  $d'E_Pd$.  By Claim 3, $d'  \in W_d$.  So $c \sim d'$.  By Claim 1, $f_P(c) \in \acl(d')$.  Since $\Ind{c}{d'}{f_P(d)}$, we have $f_P(c) \in \acl(f_P(d))$.
 \eprf
 
  Local finitenss follows from Claim 7 and compactness.  In particular  each  $\sim$-class on $P$ is a countable union of $E_P$-classes.

 This ends the proof of  \propref{appr-graphs-inf}.
 
 \end{proof}
 
 \begin{proof}[Proof of \thmref{appr-graphs}]
 
We continue with the proof of \thmref{appr-graphs}.

 The equivalence relation $E_P$   is an intersection of $0$-definable equivalence relations; so for one of these equivalence relations $E_0$, for 
 some 0-definable set $D$ with $b \in D$, we have:  for all $b',b'' \in D$,  with $\d=\d_P$ we have: 
  \[ (*) \hspace{30pt}  b'E_0 b'' \hbox{ implies }  d_X(b',b'') \leq 4 \delta \]
  
    Let $f(b)=b/E_0$.
 A conjugate of $f(b)$ has the form $f(b')$; a conjugate over $\acl(f_P(b))$ has the form $f(b')$ where $tp(b/ \acl(f_P(b)) ) = tp (b'/ \acl(f_P(b)))$;
 and so since $E_0$ refines this equivalence relation, we have $b E_0 b'$, so $f(b)=f(b')$.  Thus  $f(b) \in \acl(f_P(b))$.

On the other hand, we saw that if $c,b$ are connected then $e \in \acl(c)$ and so $\be \in \acl(c)$.      In particular,
if $d_X(b,c) \leq 5 \delta_P$ then $\be \in \acl(c)$.  Thus we can also choose $D$ so that  for all $ b' \in D, c' \in G$, with $\d=\d_P$ we have:
\[ (**) \ \ \hbox{ if }   d_X(b',c') \leq 5 \delta  \hbox{  then }  f(b') \in \acl(c')\]
(here, $ f(b') \in \acl(c')$ can be replaced by a single formula, obtainable by compactness, that guarantees it.)

Let $F$ be the set of all triples $(D,f,\delta)$ with the above two properties.    We saw that any $b$ lies in $D$ for some $(D,f,\delta) \in F$.  By compactness there exists finitely many such
pairs $(D_i,f_i,\delta_i)$ such that $\union_i D_i = G$.  We may refine them so as to be disjoint, and then take the union to obtain
a single function $f$ defined on $G$, into the disjoint union of the ranges of $f_i$; so a fiber of $f$
is a  fiber  of some $f_i$ on $D_i$, and thus (*) holds with $\delta=\d_i$; and (**)   holds with $D=G$, $\d=\d_i$.
    \[ (*_i) \hspace{46pt}  \hbox{ if } b,b' \in D_i, f(b)=f(b_i)  \hbox{ then } \ \ \ d_X(b',b'') \leq 4\d_i\]  
\[ (**_i)  \hspace{20pt}  \ \ \hbox{ if }  b' \in D_i, c' \in G, d_X(b',c') \leq 5\d_i  \hbox{  then }  f(b') \in \acl(c')\]
    
Define $H= \{(x,y) \in G: f(x)=f(y) \}$.  By (*), since $f(x)=f(y)$ implies $x,y$ lie in the same $D_i$,   (1) holds  with $m=\max \d_i$.  Similarly (2) follows from the $(**)_i$ and compactness (using $1 \leq \d_i$ for each $i$.)  

 Towards (3), re-order the $D_i$ so that $\d_i \geq \d_j$ for $i \leq j$.

 Fix $i \leq j$, and let  $b'  \in D_i$, $c' \in D_j$.  
If $d_X(b',c') \leq 1$  
then $f(b') \in \acl(c'')$ for any $c'' \in D_j$ with $f(c'')=f(c')$; the reason is that $d_X(c',c'') \leq 4\d_j$ by 
$(*)_j$, 
 so $d_X(b',c'') \leq 4\d_j+1$ and by $(**)_i$,  since $4\d_j+1 \leq 5\d_i$, we have $f(b') \in \acl(c'')$.    From this it follows that $f(b') \in \acl(f(c'))$.  
 
  Thus (3).

Let us now prove part B.   We work over a subfield $F_0$ so that $F$ has elimination of imaginaries, and model completeness in the form:  every definable set is a projection of the set of rational points of a variety under a finite morphism.   See \remref{appr-graphs-remark} (\ref{fieldconstants})
for a better statement with control of $F_0$.  

Thus the equivalence relation $(x,y) \in H$ can be written as $f(x)=f(y)$ for some definable function $f$ into $F^k$.  For each complete type $P$ of $G$,
  one can find  a variety $G_P$, a morphism $h'_P:  G_P \to \Aa^k$ with finite fibers, and regular functions $\phi'_P$ on $G_P$
  such that for any $a \in P$, for some $b$ we have $h_P'(b)=a$ and $f(a)= \phi'_P(b)$.   This is a form of the 
 algebraic boundedness  of pseudo-finite fields (\cite{vddries}).  It follows from Ax's theorem applied to the graph of $f$.  (This is a finite projection of
 the rational points of a  quantifier-free type of ACF, which can be enlarged to be a locally closed subset of some variety; but this is then itself a variety.)
 
 By Ax's theorem, any  definable
subset of $G_P$ has the form $g(\tG_P)$ for some variety $\tG_P$ and morphism $g$ with finite fibers. In particular
this is true for $\{y \in G_P:  \phi'_P(y) = f(h'_P(y)) \}$.  Define $\phi_P$ on $\tG_P$ by $\phi_P=\phi'_P \circ g$, and ${\rho}_P = h'_P \circ g$.   Then $f({\rho}_P(x)) $ is given by 
regular functions, and  $P \subseteq \rho_P(\tG_P)$. 
   By compactness, we can find a finite number of ${\rho}_i: \tG_i \to \Aa^k$  such that $f({\rho}_i(x))$ is given by regular functions on
each  $\tG_i$, and $\union_i {\rho}_i(\tG_i) = G$.   Let $(\tG,{\rho})$ be the disjoint union of $(\tG_i,{\rho}_i)$.    
Then ${\rho}: \tG \to \Aa^k$ has finite fibers and   image $G$,  and there exists a tuple of regular functions $\phi$ on $\tG$, 
such that $\phi(x) = f({\rho}(x))$ for $x \in \tG$; the 'moreover'   follows.

\end{proof}

\begin{rems} \label{appr-graphs-remark}   \rm \begin{enumerate}

\item The 'moreover' is made explicit as follows:   say $G \subset F^k$.   Then there exists  a     morphism  of varieties $h: \tG \to \Aa^k$ with finite fibers,  and a tuple of regular functions $\phi$ on $\tG$, 
such that $\phi(x)=\phi(x')$ if $h(x)=h(x')$,  
and 
\[ H = \{(h(x),h(y)):  x,y \in \tG(F), \phi(x) = \phi(y) \} \]
In particular each $H$-class has the form $h(\phi \inv(c) \meet \tG(F))$ for some $c$.  

\item  One can add that piecewise on $G$ (i.e. on each of finitely many definable pieces $G_\nu$), the rational invariants
suffice to determine a class of $H$ up to finitely many possibilities.  In other words
 there are rational functions $\psi$ on $G_\nu$, such that the equivalence relation $H'$ defined by
$\psi(x)=\psi(y)$ contains $H | G_\nu$, and each $H'$-class is a finite (bounded) union of $H$-classes.  The $\psi(x)$ will list the 
coefficients of the minimal polynomials satisfied by the $\phi(y)$ for $h(y)=x$.

\item  The theorem can be stated for the family of finite fields, in place of a single pseudo-finite field.  Then one adds uniformity to the 
definability assumption and ot the conclusion.  In particular   $m$ is bounded, that
the valencies in (2) are uniformly bounded, and the complexity of $h,\tG, \phi$ is bounded independently of $p$.

\item  
Note in particular that   each $i \leq n$, $X \meet D_i^2$ induces a finite valency graph on $f(D_i)$.  Globally on $f(G)$, there is a multi-layered structure with layers $f(D_i)$; the graph has finite valency in the downwards direction, regarding arrows going from $f(D_j)$ to $f(D_{\leq j})$,  but not necessarily upwards.    

This kind of tree-like structure  is intrinsic and cannot be avoided; 
see \exref{nofirstexpansionbound}.

 \item 
Part A of \thmref{hp}
admits an analogue for stable or even simple theories, near a given regular type $p$,  using semi-regular $p$-weight in place of dimension.  Is there a similar generalization of   \thmref{appr-graphs}?  (It seems likely.)

\item  For simple Robinson theories, even of $SU$-rank one,   it was shown in   \cite{simple}  that $E/E_P$ may be a connected compact space.  
 But pseudo-finiteness was not considered there, and 
it would be interesting to see what happens when this  hypothesis is added.

 \item  Part B  of \thmref{appr-graphs}  is valid for structures interpretable in PF, not necessarily with the full induced structure.  (Thus the transitivity assumption is easy to attain, given a group action.)

\item \label{fieldconstants}  In Part B, we needed to  adjoin constants  so as to have model completeness as in Ax.  If in the statement we 
replace $V(F)$ by $V(F_n)$, where $F$ is a (periodic) difference variety and $F_n$ is the degree $n$ extension of $F$,  or alternatively Artin symbols,   this becomes unnecessary; see \cite{pac}.  In this formulation
  the theorem is valid over $F_0^{alg}$ where $F_0$ is the field of $0$-definable elements of $F$.    The reason for going to $F_0^{alg}$ rather than $F_0$
  is the  need to  name  the algebraic imaginary parameters $\phi_n$, 
 coding a generating automorphism of the Galois group of the degree $n$ extension of $F$, in order to have elimination of imaginaries; these 
 are definable over $F_0^{alg}$.

\item   Part A of  \thmref{hp}   can easily be recovered from part A of  \thmref{appr-graphs}.    Part B of \thmref{appr-graphs} partially  generalizes
Part B of  \thmref{hp}, but to obtain a group structure on the the covering variety appears to require additional argument.
 
\end{enumerate}

\end{rems}

\ssec{Definability of connected components}

We showed in general that connectedness over $\Ff_p$ of a definable graph can be definably reduced to a locally finite graph.
Here we impose  stronger assumptions that allow ruling out the locally finite graph, and thus showing definability of the connected components.
 This will be the graph analogue of \secref{unip}.    In particular under these assumptions, 
  if the graph generated by $X$ on $\Ff_p^n$ fails to be connected for infinitely many $p$,
  then this is so for an algebraic reason.

Recall that a definable group is {\em connected} if it has no proper definable subgroups of finite index. If $F$ is a pseudo-finite field
and $G$ is a simply connected algebraic group, then $G(F)$ is connected (\cite{HP2}.)  

\def\gG{\mathcal{G}}
 Let $T$ be an S1-theory.  Let $\gG$ be an ind-definable group with a transitive action on a definable set $(V,X)$.   Let $I$ be an index set,
 and for $i \in I$, let $(J_i(v): v \in V)$ be a definable family of definable subgroups of $\gG$.   Assume
conjugation invariance:  $J_i(gx) = g  J_i(x) g\inv$.   Let $\sim_V(a)$ denote the connected component of $a \in V$.  

\begin{cor} \label{inddefgrp}  Assume: in any model of $T$,
  whenever $a \in V$,    $J_i(a)$ is a connected definable group; 
   $ J_i(a) \cdot a \subseteq \sim_V(a)$; and if $(a,b) \in X$ then $j \cdot a = b$  for some $i \in I$ and some $j \in J_i(a)$, or dually.
 Then $\simv$ is definable. \end{cor}
 
\prf    
By assumption, $\gG$ acts on $(V,X)$ by automorphisms; so as a graph, $(V,X)$ has a unique type $P$.  Let $E=E_P$ be the
equivalence relation of \propref{appr-graphs}.   Then the $\gG$-action preserves $E$; and $V/E$ is a locally finite graph. 
Now for $a \in V$, and $i \in I$, by assumption $J_i(a)  \cdot a \subseteq \sim_V(a)$; by compactness, $J_i(a) \cdot a$ is within a finite radius
$X/E$-ball $B$ around $a$; now $B$ is finite by local finiteness; so a finite index definable subgroup of $J_i(a)$ fixes it pointwise.
But $J_i(a)$ is connected, so $J_i(a)$ acts trivially on $B$; in particular $J_i(a)$ fixes the $E$-class of $a$.  
Let $X'$ be the directed graph obtained by drawing a directed edge from  
$a$  to $j \cdot a$, for any $j \in J_i(a)$.  
We have just shown that the forward neighbours of $a$ are in the same $E$-class as $a$.  Applying this to other elements, we see that the
backwards-neighbours of $a$ are also in the same $E$-class.  But by assumption, any edge of $X$ is in $X'$ or the reverse graph;
so all $X$-neighbours of $a$ are in the same $E$-class.  Conversely we have $E \subset \simv$;
so the $E$-classes coincide with the $\simv$ classes.  Hence both are definable.
\eprf

For example we may take $\gG$ to be the   ind-definable  group of polynomial automorphisms 
of $V=\Aa^n$.   By a {\em transvection} (not necessarily linear) we mean a map $(x,y) \mapsto (x,y+f(x))$, where $(x,y)$ is a partition of the coordinates
of $\Aa^n$, and $f$ is a polynomial.  Note that any transvection $t$ forms part of a unipotent group $(t^ \alpha: \alpha \in F)$, where
$t^\alpha(x,y) = (x,y+\alpha f(x))$; in the theory $T$ of pseudo-finite fields of characteristic zero, this group is connected.      

\begin{cor}

  Let  $j_1,j_2$ be  definable, conjugation invariant maps from  $V$ to $\gG$, so that $j_i(a)$ is a transvection;    draw an edge from $a$ to $j_i(a) \cdot a$.  Let $\G$ denote the associated graph (of valency at most $4$.)    Then the connected component $\simv(a)$ in the finite field $\Ff_p$
is definable uniformly in $a$ and in $p$. \qed \end{cor}

This is a special case of \corref{inddefgrp}.  The proof shows also that  the graph $(V,X')$, obtained by connecting $a$ to each
$j_i^t(a)$, $i \in I$, $t \in F$, has  the same connected components as $(V,X)$, and   each connected component has finite $X'$-diameter.

\begin{example}  \rm  Fix   a linear algebraic  group $G$ (say over $\Zz$).   Consider the action of $G$ by conjugation on $G^n$.
Define a graph structure $\G$ on $G^n$, of valency at most $n$, by letting $a$  be adjacent to $a_i \inv  a a_i $, where $a_i$ is one of the $n$ components
of $a = (a_1,\ldots,a_n)$ and $k \in \Zz$.   Let $p$ be a prime,   let $a_1,\ldots,a_n \in G(\Ff_p)$ be unipotent, $a=(a_1,\ldots,a_n)$, and let $C(p;a)$
be the component of $a$ in $\G(\Ff_p)$.     Then 
$C(p;a)$ is definable in $\Ff_p$ uniformly in $a$ and $p$.  

If one `speeds up' the graph by declaring $a_i^{-k}  a a_i^k$ to be adjacent to $a$, then  here exists a bound $\beta$ such that for all primes $p$ and all $a$ as above, 
$C(p;a)$ has diameter $\leq \beta$.   
   
   \end{example}

Each component of $\G$ is contained in a  conjugacy classes of $n$-tuples of unipotent elements,   in the subgroup of $G$ that they generate.
The subgroup $G(a)$ generated by an $n$-tuple $a$ is itself a definable function; this is a consequence of \thmref{hp}.   From \thmref{appr-graphs}
we learn that the component of $a$ is a definable subset of the $G(a)$-conjugacy class of $a$.  In particular when $G=SL_r$ and
$a$ is a generating $n$-tuple of $SL_r(\Ff_p)$, I do not know if the conjugacy class of $a$ is  connected in $\G$; if it is not, the theorem shows that it
 is due to an algebraic invariant.     In particular if the conjugacy class in $G(\Ff_p)^n$ is connected away from a density zero set of primes, then it is connected for all but finitely many primes.

 See   \cite{BGS}, \cite{BGS1},   \cite{avnigarion}   \cite{cantat-loray} and references there for deep  work on more specific instances or classes of definable graphs.

 \ssec{Example of unbounded first expansion radius} \label{unbounded-expansion}

We conclude with a simple example showing that the "layered" or "tree-like" structure in the statement of  \thmref{appr-graphs} is unavoidable.
 
 Let $(G,X)$ be a (symmetric, reflexive) graph, definable in a structure of finite S1-rank  $d(G)$.   
 
 For $a \in G$, let
 $\xi(a)$ be the smallest $n$ such that  for all $n' >n$    we have $d(B_{n'}(a)) = d(B_n(a))$;
 where $B_n(a)$ is the $n$-ball of $(G,d_X)$.   

Here is an example where $\xi$ is unbounded.   It can be understood in ACF or in PF, and shows that the finite partition in \thmref{appr-graphs} is unavoidable.

 \begin{example}  \label{nofirstexpansionbound}   Let $G=\Aa^2$, and let $f: G \to G$ be an endomorphism of infinite order.
 Let $C$ be an irreducible curve on $G$, also of infinite order under $f$, i.e. not $f$-preperiodic. 
 Let $X$ be the union of the graph of $f$ and $C \times C$.  Then $\xi$ is unbounded. 
   Hence, there is no definable  equivalence relation $H$ on $G$ such that:
 \begin{enumerate} 
\item For some $m$, if $(a,b) \in H$ then $d_X(a,b) \leq m$.  
\item $X / H$ has finite  valency.
\end{enumerate}  \rm
Indeed, $\xi(g)=n$ if $g \in f^{-n}(C) \m \union_{k<n} f^{-k}(C)$.  Thus $\xi$ is unbounded.  Suppose an $H$ with (1,2) exists.  Let 
$a \in G$ be generic;  it suffices that $a \notin \union_{n \in \Zz} f^n C$.   
Since all $d_X$-balls around $a$ are finite,  and by (1),  the $H$-class of $a$ is finite.     Hence all $H$-classes must be finite except
possibly on some finite union of curves $f^{-n}(C)$.   On the other hand, the curve $C$ itself meets only finitely many $H$-classes, since the complete graph on $C$  has finite valency
modulo $H$.   
Let $n$ be the greatest integer such that $f^{-n} (C) / H$ is finite.  Then $f^{-n-1}(C) /H$ is infinite, yet $f$ induces a finite valency
graph on  the product $f^{-n-1}(C) /H \times f^{-n} (C)/H$, which is impossible.   \end{example}

 \ssec{Larsen-Pink}    One may ask about a graph-theoretic version of the Larsen-Pink inequalties.    We obtain such an analogue in a basic case,
 leaving  the more general case as a question.

  Consider an enrichment of the
 theory of fields, with a reasonable dimension function $\delta$, as in \cite{HW}. 
   For any partial type $S$, write $d(S)$ for the dimension of the 
 Zariski closure of $S$, also $d(c)=d(S)$ if $S=tp(c/M)$. Write $\d(S)=\alpha$ if $\alpha = \inf \d(S')$ as $S'$ varies over definable sets containing $S$.

  Let $P,Q$ and   $R \subset P \times Q$ be complete types over some base field $M$, such that  $\d(P),\d(Q), \d(R)$ exist.    
  Assume:   \noindent \\
$\diamond$ \ \ \ \ for  $(a,b) \in R$ we have 
$M(a)^{alg} \meet M(b)^{alg} = M$.

\begin{lem}\label{lp-graphs-1}
 Assume $d(P)=d(Q)$ and $d(R) = d(P)+1$.   Then  
 \[d(P) \d(R)  \leq \d(P) d(R)  \]
 \end{lem}
 
  \prf   let $R^t = \{(b,a): (a,b) \in R\}$.  Consider $(a_0,\ldots,a_k)$ with $(a_i,a_{i+1}) \in R \union R^t$,
 such that $\Ind{a_0}{a_{i+1}}{M(a_i)}$.  Then $d(a_0a_{i+1}) \leq d(a_0a_ia_{i+1}) = d(a_0a_i)+1$.
 Thus $d(a_0a_i)  \leq   d(P)+i$.  Once  $d(a_0a_i)=d(a_0a_{i+1})$ we have     $\Ind{a_0}{a_{i}}{M(a_{i+1})}$ and so by considering the canonical base we have
 $\Ind{a_0}{M(a_i)^{alg} } {M(a_i)^{alg} \meet M(a_{i+1})^{alg}}$, i.e. by assumption,  $\Ind{a_0}{a_i} {M}$; so $d(a_0a_i) = 2 d(P)$.
 Thus equality holds only at $i=d(P)$; and for $i< d(P)$ we have $a_i \in M(a_0a_{i+1})^{alg}$. 
 
 So far we used only algebraic independence.  But now choose $(a_0,\ldots,a_k)$ with $(a_i,a_{i+1}) \in R \union R^t$ such that
 $\Ind{a_0}{a_{i+1}}{M(a_i)}$ holds in the sense of $\d$ (and a fortiori algebraically.)  Using  $a_i \in M(a_0a_{i+1})^{alg}$, we have
 $\d(a_i/ M(a_0a_{i+1})=0$, and it follows that  for $i \leq d(P)$ we have $\d(a_0a_{i}) \geq   \d(a_0a_{i-1}) +\d(R_b)$ where $R_b = \{x: (x,b) \in R\}$.  So   
 $ 2 \d(a_0) \geq  \d(a_0a_i) \geq \d(a_0) + i\d(R_b)$.   For $i=d(P)$ we obtain   $\d(a_0) \geq d(P) \d(R_b)$.
 Thus $\d(R_b) \leq \d(P)/ d(P)$ and  also $\d(R) = \d(P) + \d(R_b) \leq \d(P) (1+ 1/d(P)) $ so $\d(R) d(P) \leq \d(P) d(R)$.     \eprf

  \begin{question} What can be said without the condition $d(R)=d(P)+1$?  \end{question}
  
  Let us now see what assumption $\diamond$ amounts to for graphs arising from a definable subset $X$ of a  definable group $G$.   First we analyze the condition $M(a)^{alg} \meet M(b)^{alg} = M$ in this case, showing that it means in essence that $X$ generates $G$.  
 The argument of this paragraph is valid in theories of finite Morley rank,  in particular the case that concerns us, ACF.   
Let $G$ be a connected definable group, $X \subset G$ a complete type over $M$.    Let $(b,c)$ be generic in $ G \times X$, and let $a=bc$.

\claim{}  $\acl(M(a)) \meet \acl(M(b)) = M$ iff $X$ is not  contained in an $M$-definable coset  of  a proper  subgroup of $G$. 

 \rm Indeed if $X \subset Hm$
then $aHm= bHm \in \acl(M(a)) \meet \acl(M(b)) $ but by genericity of $b$ in $G$, $bHm \notin \acl(M)$.    Conversely, assume $\acl(M(a)) \meet \acl(M(b)) \neq M$.   Let $tp(b'/  \acl(M(a))= tp(b/ \acl(M(a))$, with $b',b$ independent over $M(a)$.  Since $b=a c \inv$, there exists $d \inv$
with $(b,c,d)$ generic in $G \times X \times X$, and $b' = a d \inv = b c d \inv$.    If $e \in \acl(M(a)) \meet \acl(M(b)) \m M$,  we still have $e \in \acl(M(b'))$.
Continuing this way, for arbitrarily large $n$ we can find $(b,c_1,\ldots,c_{2n}$ generic in $G \times X^{2n}$, with 
$\acl(M(b c_1 c_2 \inv \cdots c_{2n-1} c_{2n} \inv)) \meet \acl (b) \neq M$.  For large enough $n$, $c_1 c_2 \inv \cdots c_{2n-1} c_{2n}\inv $ is  a generic of 
element of a definable group $H \leq G$; now $H = G$ is impossible since for generic $c \in G$, we do have $\acl(Mb c) \meet \acl(Mb) = M$.  Thus
$c_1 c_2 \inv \cdots c_{2n-1} c_{2n} $ is the generic of a proper definable subgroup, and the result follows.

   \begin{rem}  Let $G$ be a definable group, $X$ a definable subset.  
   Let $P=Q = G$ and let $R = \{(x,y) \in G^2:  x y \inv \in X \}$.  Then the  inequality \[d(P) \d(R)  \leq \d(P) d(R)\]    implies  the   Larsen-Pink   inequality 
   \[ \d(X)d(G) \leq \d(G) d(X).\]
     \rm       To see this, renormalize $\delta$ so that $\delta(G)=d(G)$.   Then the first inequality becomes $\d(R) \leq d(R)$;  equivalently 
 $\d(G)+\d(X) \leq d(G) + d(X)$.   Using   $\d(G)=d(G)$ additively now,  we obtain $\d(X) \leq   d(X)$ and hence $\d(X)d(G) \leq \d(G) d(X)$.  \end{rem}

\section{ The Galois group of a NIP measure}  \label{NIPgr}

Let $\mu(x)$ be a definable measure in a NIP theory,   let $\phi(x,u)$ be a formula, and $q(y)$ a type (over $0$.)

Define an equivalence relation on $q$:  $b E_{\mu,\phi} b'$ iff for $\mu$-almost all $x$, $\phi(x,b) = \phi(x,b')$.  This is then 
a $\bigwedge$-definable equivalence relation on $q$.  In a NIP theory, (\cite{NIP1}), it is co-bounded and thus the quotient $q/E$ is
compact in the logic topology. 
The group $G$ of automorphisms of $q/E$ is likewise a compact group.  
Since $q$ is a complete type, $G$ acts transitively on $q/E$.   In this section we will study this group $G=G_{\mu,\phi,q}$. 

 We could also consider finitely
many formulas $\phi_1(x,u_1),\ldots,\phi_k(x,u_n)$, each with a type $q_i(u_i)$; but by standard tricks one can find a single type $q=q(u_1,\ldots,u_n)$ projecting to $q_i$
at the $i$'th position, and a formula $\phi$, such that each  $\phi_i(x,c_i)$  $(c_i \models q_i)$ is equal to some $\phi(x,c)$ (with $c \models q$).   Thus the sets of such formulas $\phi(x,c)$ with $c \models q$ forms
a  directed system; and the projective
limit of the groups $G=G_{\mu,\phi,q}$  can be viewed as the Galois group of the space of weakly random  types.   The individual groups 
$G_{\mu,\phi,q}$ are the fundamental building blocks.   

   It turns out that a totally disconnected part may appear  both as a quotient ($G/G^0$) and as a normal subgroup $K/G^{00}$ of $G^0/G^{00}$.
  We show in \thmref{nipsymmetry} that the ``archimedean core" $G^{0}/K$ is a finite dimensional real Lie group.

As a corollary we obtain  \thmref{anand-nip-groups},  addressing a basic question of Anand Pillay regarding local finite-dimensionality of    
groups of connected components $G/G^{00}$ for a definable group $G$ in a NIP theory.   
  We define the local connected component $G^{00}_{\phi,q}$
 attributable to $\phi(x,u)$ ranging over a given type $q(u)$, and show again that it  is a finite-dimensional Lie group, up to profinite group extensions
 above and below.

We first present some basic material characterizing the target,  profinite-by-Lie-by-profinite groups.    
In this section, by 'Lie group' we will   mean:  finite-dimensional Lie group, with finitely many connected components.
Only compact Lie groups will be considered.  
All topological groups are taken to be Hausdorff.   

 \ssec{Profinite extensions of compact Lie groups} \label{profinite}
 
 Let $L$ be a connected, finite dimensional   Lie group, with identity element $1$, and let $\pi_1(L):=\pi_1(L,1)$.   Then 
the universal cover $\tilde{L}$ admits a group structure, induced for instance by multiplication of paths;
there is an exact sequence $1 \to \pi_1(L) \to \tL \to L \to 1$.    The connected group $\tilde{L}$ acts trivially
by conjugation on the finitely generated group $\pi_1(L)$; so this is a central extension.  
Any  subgroup $N$ of $\pi_1(L)$ of finite index is normal in $\tL$, and we can form the quotient
$\tL/N$.  Taking the inverse limit over all $N$ we obtain a group $\widehat{L}$;   
by construction  it is an extension of the Lie group $L$ by the  finitely generated profinite group $\widehat{\pi_1(L)}$,
profinite completion of $\pi_1(L)$.    $\widehat{L}$ admits a surjective continuous  map onto any connected, pointed, finite covering group of L (unique with $1 \mapsto 1$), and
at the limit, a surjective continuous  map onto any  connected  central extension $E$ of $L$ by a profinite group.  We call   $\widehat{L}$ the {\em universal profinite covering group} of $L$.

\ssec{ Groups of finite archimedean rank}

\<{defn}\label{p-by-lie-def}   \rm A connected compact topological group $K$ is {\em profinite-by-Lie}, if  it  has a closed normal profinite subgroup 
$K_0$  such that $K /K_0$ is a  Lie group.

We will say that   compact topological group $G$ is profinite-by-Lie-by-profinite, or
  {\em has finite archimedean rank},  if the connected component $G^0$ of $G$ is  profinite-by-Lie.
\>{defn}

We will see that $G$ is actually profinite-by-Lie-by-finite in this case, though $G/G^0$ need not be finite.

\<{rems}  \rm \<{enumerate} 
\item  If $K$ is a compact group, then $K/K^0$ is profinite.
\item  
Let $G$ be a topological group with a  profinite normal subgroup $K$, such that  $G/K$ is a Lie group with finitely many connected components. 
 Then $K$ is determined up to finite index
  by these conditions.  Indeed if $L$ is another such group, then $KL/K$ is a normal profinite subgroup of $G/K$; as  Lie groups do not contain infinite normal profinite groups, $KL/K$ must be finite, and by symmetry $K,L$ are commensurable.

\item 
If  $K_0$ is profinite and normal in a connected group $G^0$ then it is central in $G^0$.  Indeed the connected group $G^0$ acts by conjugation on $K_0$,
so each orbit is connected; but the totally disconnected group $K_0$ has no connected subsets other than points.

\item Let $A$ be a compact abelian group; we have $A \cong Hom(\widehat{A}, U_1)$ where $U_1$
is the unit circle in $\Cc$, and  $\widehat{A}$ is the character group.  Then 
  $A$ has finite archimedean rank iff $\widehat{A}$ does as a discrete group, i.e. $\Qq \tensor \widehat{A}$ is a finite-dimensional $\Qq$-vector space.
  
  \item  A compact group $G$ has finite archimedean rank iff the center $Z$ of $G$ has finite archimedean rank,
  and $G/Z$ has a finite-index subgroup isomorphic to the product
  of  a profinite group with a  compact real Lie group.  (See \propref{plp} (\ref{plp3}).
  \>{enumerate}
\>{rems}

  \begin{example}  \label{plp-ex} \rm
Take a   finite dimensional Lie group $L$; a central extension $\bar{L}$ of $L$ by a profinite group $A$; another profinite group $P$, and a central extension $\bar{P}$ of $P$ by $A$;
  and form $G = \bar{P} \times^A \bar{L}$, the quotient of $ \bar{P} \times   \bar{L}$ by the antidiagonal subgroup of $A^2$.   Then $G$ is profinite-by-Lie-by-profinite, and compact if $L$ is. 
  
  In fact $G$ is also profinite-by-Lie in this case, as shown by the image of $\bar{P} \times A$ in $G$.  Along with \propref{plp} (\ref{plp3}), 
  this explains the remark in brackets in \defref{p-by-lie-def}.
    \end{example}

\<{prop} \label{plp}  Let $G$ be a compact topological group. \begin{enumerate}
\item if $G$ is connected profinite-by-Lie, then $G$ is a quotient of the  universal profinite cover of a compact connected Lie group.
\item If $G$ is   Lie-by-profinite,  then some subgroup $G'$ of finite index in $G$  splits
as a direct product of a connected Lie group with a profinite group, both normal in $G$.  
 
\item \label{plp3} If $G$ is of finite archimedean rank, then 
$G$ has an open finite index subgroup isomorphic to $  \bar{P} \times^A \bar{L}$, with 
  $\bar{P},A,\bar{L}$ as in Example \ref{plp-ex}.    Hence $G$ is profinite-by-Lie-by-finite.  
\item  $G$ has finite archimedean rank iff there exists a bound $r$ and a family $(N_i: i \in I)$ of closed normal subgroups of $G$, 
closed under finite intersections,
such that $\meet_{i \in I} N_i = (1)$ and each $G/N_i$ is a Lie group of dimension $\leq r$.
 
\item   Let $H$ be a closed subgroup of $G$.  If $G$ has finite archimedean rank, then so does $H$.
  \end{enumerate}
  \>{prop}

\prf (1)  This follows from the discussion in \secref{profinite} above.

(2)   We have an exact sequence 
 \[ 1 \to L \to G \to P \to 1  \]
with $L$ a Lie group, $P$ profinite.  As 
$L$ has no small subgroups,    $G$ has an open neighborhood $O$ containing no nontrivial subgroup of $L$.
By Peter-Weyl there exists a continuous homomorphism $r: G \to U_n \leq GL_n(\Cc)$ with kernel $R$ contained in $O$ and hence meeting $L$ trivially.  
(See e.g. \cite{tao}, 1.4.14 for this consequence of Peter-Weyl.)
On the other hand $G/RL \cong P/Im(RL)$ is a profinite group, that embeds in a subquotient
of $GL_n(\Cc)$; so it finite.   Thus   $RL$ has finite index in $G$, and  splits as a semi-direct product of $R$ with $L$.  Conjugation induces a homomorphism
\[ R \to Aut(L) \to Aut( Lie(L)) \subset GL_n(\Cc)\]
where $Lie(L)$ is the Lie algebra of $L$.  
The composition $R \to GL_n(\Cc)$ has finite image as 
R is profinite; so some open subgroup $R'$
of $R$ acts trivially on $Lie(L)$, hence on the connected component $ L^0$; so a finite index normal subgroup  $R''$ of $R$ acts trivially on $L$.  We may take $R''$ normal in $G$.  Thus $R''L$ is a direct product as stated.  Note that $R \cong RL/L$ is isomorphic to a subgroup of $P$,
and hence is  profinite.

(3)  
We have two exact sequences:
 
\[ 1 \to {G^0} \to G \to P \to 1  \]
\[ 1 \to Q \to {G^0} \to L \to 1  \]
with $P$,$Q$ profinite, $L$ Lie, ${G^0}$ connected.  $Q$ is contained in the center of ${G^0}$
(since a connected compact group has Abelian $\pi_1$, and acts trivially by conjugation on it.)

Since $L$ has no small subgroups, there
 exists a neighborhood $U$ of $1$ in $G$ such that any subgroup  of $G$ contained in  $U \meet {G^0}$ is contained
 in $Q$.     By Peter-Weyl   there exists  a homomorphism
$r: G \to GL_n(\Cc)$ with kernel $R \subset U$.   So $R \meet G_0 \leq Q$.  $G/(RG^0)$ is a quotient of the profinite group $G/G^0$
and also of the group $G/R$ that has finitely many connected components; so  $G/(R{G^0})$ is finite, i.e.
$RG^0$  has finite index in $G$.   We may thus assume $R{G^0}=G$, and likewise replace 
    $Q$ by $R \meet {G^0}$, and $L$ by the new ${G^0}/Q$; 
now $Q$ is normalized by $R$.   

Since  the closed subgroup   $R \meet {G^0}$ is commensurable with $Q$, it is profinite.   Since $R/G^0$ is also profinite,
$R$ is   profinite.

% ****%double line breaks checked to here.  ****

Since ${G^0}$ is connected, it acts trivially on the totally disconnected group $R$; i.e. ${G^0}, R$ commute.  Thus 
  $R$ acts trivially on ${G^0}$.  (Thanks to the referee for this short argument.)    In particular $Q$ is central in $R{G^0} =G$.  So $G =  R \times_Q {G^0}$.
 
 (4)  
  
  Assume such a bound exists.
   $\dim(G/N_i)$ is bounded independently of $i$.  It follows that for some $N_0$, for all $N_i \leq N_0$, 
   $\dim(G/N_0)=\dim(G/N_i)$; so  $[N_0: N_i] < \infty$.  Since $\meet \{N_i: N_i \leq N_0\} = 1$,  $N_0$ is profinite; while $G/N_0$ is
   a compact Lie group.  
   
Conversely, assume $G$ has finite archimedean rank.  By (3), $G$ is profinite-by-Lie-by-finite, so it has an open profinite-by-Lie normal subgroup $G_1$
of finite index;  and clearly $G_1$ has   a family $(N_i':  i \in I)$, with $\dim(G_1/N_i') \leq r'$ and $\meet_i N_i'=1$.    Let $N_i$ be the intersection of the conjugates of $N_i'$.
Each $N_i'$ has at most $[G:G_1]$ conjugates, and thus $G/N_i$ is a Lie group of dimension $\leq r' [G:G_1]$; clearly $\meet_i N_i =(1)$.  

(Thus $G$ is profinite-by-(not necessarily connected Lie).)

(5)   Follows from (4), since a closed subgroup of a Lie group is a Lie group of lower (or equal) dimension.
 \eprf 

We actually need the following lemma only when $L$ is compact, but state it more generally with a view to a later generalization (see \ref{open7} (\ref{locallycompact}).)

   \begin{lem} \label{manifold} Let $L$ be a Lie group acting transitively and faithfully on a connected manifold $Y$, with compact stabilizer $S$ of a point $y \in Y$.   Then   $\dim(L) \leq  \dim(Y)^2$.   \end{lem}
   \prf   
   Since $S$ is compact, fixing $y \in Y$ one can find an $S$-invariant inner product $b$ on the tangent space $T_y$, and this propagates to an $L$-
     invariant Riemannian metric on $Y$.   Using the exponential map
     one sees that the pointwise stabilizer $S'$ of the  tangent space $T_y$ fixes a neighborhood of $y$.
     The set of points $y'$ such that $y'$ and each vector in  $T_{y'}$  is fixed by $S'$ is closed and open ,
     hence equals $Y$.   So $S'=1$, hence the homomorphism $L \to End T_y$ is injective.
  Thus $\dim(S) \leq Aut(T_y,b) =  {{ \dim Y} \choose {2}}$, and hence  $\dim(L) \leq \dim(Y) + {{ \dim Y} \choose {2}} = \frac{\dim(Y)^2+\dim(Y)}{2}   \leq  \dim(Y)^2$.
\eprf

Let $Y$ be a compact $C^1$ differentiable manifold.  Let $d_Y$ be a metric on $Y$.  We say that $d_Y$
is {\em compatible with the manifold structure} if for any chart $c: W \to Y$, where $W$ is an open subset of
$\Rr^n$ and $d_W$ is the metric induced from $\Rr^n$, the map $c: (W,d_W) \to (Y,d_Y)$ is bi-Lipshcitz.
Any two compatible metrics on $Y$ are bi-Lipschitz equivalent, and hence the packing dimension is
well-defined and does not depend on the choice of a compatible metric.   In fact we the packing dimension equals
the dimension of $Y$ as a manifold, as can be seen by going to charts and using the standard Euclidean metric.

\begin{lem} \label{manifold2}  Let $G$ be a compact Lie group, $H$ a finite-dimensional Euclidean space, $\rho: G \to Aut(H)$ a faithful unitary representation.   
Let $X=Gv$ be an orbit of $G$, with metric $d_X$   induced from $H$.  
 Then  $X$ is a smooth submanifold of $H$, and $d_X$ is Lipschitz-compatible with the manifold structure on $X$.
\end{lem}

 \prf   
Let   $F$ be the graph of $\rho: G \to Aut(H)$.    Then 
$F $ is a closed subgroup of $G \times Aut(H)$, which  is analytic by the closed subgroup theorem (Von Neumann 1929, Cartan 1930).  
 Since the projection $F \to G$ is analytic, and an isomorphism,  the inverse map  $G \to F$ is too; composing
 with  the projection $F \to Aut(H)$ we see that $\rho$ is a real-analytic, and in particular  $C^1$, isomorphism
 between $G$ and an analytic subgroup $A$ of $Aut(H)$.  
 
 $G_v = \{g \in G: gv = v\}$.   We give $G/G_v$   the quotient manifold structure (\cite{bourbaki} 5.9.5). 
Likewise let $A_v = \{a \in A: av = v \}$, and give $A/A_v$ the the quotient manifold structure.  We obtain induced analytic
isomorphisms $G/G_v \to A/A_v$.  
 
 Consider the map $\alpha:  A \to X$, $a \mapsto a(v)$, and the derivative $\alpha'$ of $\alpha$ at $1$,
 from the Lie algebra  
of $A$  into $T_v H = H$.   
 The kernel is precisely the Lie  algebra $T_1A_v$ of $A_v$.    Hence
 the induced map $A / A_v \to H$ induces an injective linear map on  tangent spaces, and
 so is  an immersion  by \cite{bourbaki} 5.7.1.   Since $A/A_v \cong G/G_v$ is compact and $A/A_v \to H$ is injective,
 it follows that the image $X$ is a differentiable submanifold of $H$, and $G/G_v \cong A/A_v \to X$ is an isomorphism
 of $C^1$-manifolds.

For the final point,  is easy to see more generally that if $Y$ is a  smooth submanifold of $\Rr^n$,   the metric induced by a Euclidean structure on $\Rr^n$
is compatible with the manifold structure on $Y$.   For instance near a point $a \in Y$, let $H'$ be the tangent space
to $Y$ at $a$ embedded as a linear subspace of $H$, and consider the orthogonal projection $\pi:Y \to H'$; at an $\e$-neighborhood of $a$ the distance to $H$ is $O(\e^2)$, and so $\pi$ is bi-Lipshcitz.

\eprf

\begin{lem}\label{lipschitz}   If $f: Y \to X$ is a Lipschitz map between metric spaces,  then the packing dimension of $Y$
 is at least equal to the packing dimension of $X$. \end{lem}  \prf   In fact let   $\g_X(\e)$ be the maximal number of disjoint $\e$-
 balls in $X$,   likewise $\g_Y(\e)$, and let $l$ be the Lipschitz constant of $f$;  then $\g_Y(\e) \geq \g_X(l \e)$.  
  If $a_1,\ldots, a_k$ are points of $X$ such that the $l \e$-balls around the $a_i$ are pairwise disjoint,
   let $b_i \in X$ be such that $f(b_i)=a_i$.  Then the $\e$-balls around the $b_i$ are pairwise disjoint:  
  if $c \in Y$ and $d_Y(b_i-c), d_Y(b_j-c) < \e$, and $d = f(c)$, 
then $d_X(a_i,d), d_X(a_j ,d) <l \e$, a contradiction to the disjointness of the $l \e$-balls around $a_i,a_j$. \eprf

\begin{lem} \label{manifold3}  Let $G$ be a compact, connected Lie group, $H$ a Hilbert space , $\rho: G \to Aut(H)$ a unitary representation.   
Let $X$ be an orbit of $G$, and assume $G$ acts faithfully on $X$.  
Let $\delta$ be the packing dimension of  $X$ as a metric space with the metric induced from $H$.   Then  $\dim(G) \leq  \delta^2$.
\end{lem}

\prf   

  By Peter-Weyl,  there exist finite-dimensional $G$-invariant subspaces $H_i' (i \in \Nn)$ whose direct sum is dense in $H$.
 Let $H_n = \oplus _{i \leq n} H_i'$, and $\pi_n: H \to H_n$ the orthogonal projection.  Since $\pi_n$ 
 is Lipschitz (in fact with constant $1$),  $\pi_n(X)$  still has packing dimension $\leq \delta$.    
 Also, if $g \in G$ fixes $\pi_n(x)$ for each $n$ and each point $x \in X$,  then $g$ fixes $X$ and hence $H$, so   $g=1$.
 Let $K_n = \{g \in G: g|\pi_n(X) = Id\}$.    Then $\meet K_n = (1)$.  By the descending chain condition on closed  subgroups of $G$ we see that $K_n=1$ for some $n$, so the action of $G$ on $\pi_n(X)$ is faithful.  Replacing $H$ by the subspace of $H_n$  generated by $\pi_n(X)$, and using \lemref{lipschitz}, we are reduced to the case that $H$ is finite dimensional.  
 
 In this case, by \lemref{manifold2}, $\delta$ equals the dimension $\dim(X)$ of $X$ as a manifold.  
 
 By \lemref{manifold}, $\dim(G) \leq  \dim(X)^2 \leq   \delta^2$.

\eprf

\ssec{Interaction of NIP formulas with definable measures} \label{NIPformulas}

We consider a   sort $X$ carrying a definable measure, and a NIP formula $\phi(x,u)$ relating another sort $U$ to $X$. 
More generally we allow  $U$ to be a $\bigwedge$-definable set (without parameters); and the measure $\mu$ may just be defined on the
 algebra $B(\phi)$ of subsets  of $X$ generated by  the sets $\phi(x,b)$, $b \in U$.
 
We will find a compact Hausdorff quotient $\sU$ of $U$ through which the interaction is mediated.  We are interested in the fibers of the natural map from 
$\sU$ to the space of types over $\emptyset$.  These are principal homogeneous spaces for a compact group $\sG$.   We show here
that $\sG$ is made up of totally disconnected groups and a single finite-dimensional Lie group.

    We begin by recalling a number of basic objects associated with $(X,U,\phi,\mu)$:
  a   compact space $\sX$ with a probability measure, a compact metric space $\sU$, a compact topological group ${\sG}$ acting on $\sX$ and $\sU$.

Let $M=(X,U,\phi,\mu)^M$ be an $\aleph_1$-saturated model.

We assume $Aut(M)$ acts  transitively on $U(M)$ (i.e. we consider one type of $U$ at a time.)  

We can define a pseudo-metric on $U$: 

\[d(a,a') = \mu(\phi(a,x) \triangle \phi(a',x)) \]
 
With associated equivalence relation:
 
\[ a E_{\phi,\mu} a' \iff \mu(\phi(a,x) \triangle \phi(a',x))=0\]

So the quotient $\sU=\sU_{\phi,\mu} = U/E_{\phi,\mu}$ becomes a metric space.

$E= E_{\phi,\mu}$ is a  $\bigwedge$-definable equivalence relation.   Clearly $a E_{\phi,\mu} a'$ iff  for all  weakly random $p$ over $M$, 
 $\phi(a,x) \in p \iff \phi(a',x) \in p$.    Since $\phi$ is NIP, the number of weakly random $\phi$ types
 has cardinality bounded independently of $M$ (see \secref{weaklyrandom}).   Hence $E$ is co-bounded.  Thus we also have a logic topology on $\sU$.  The identity
 map is continuous from the (compact) logic topology to the (Hausdorff) topology induced by the metric, so they coincide.

 By compactness,  for any $\e>0$ there is a finite bound
 on the size of a set $A \subset U(M)$ such that $\mu(\phi(a,x) \triangle \phi(a',x)) \geq \e$ for all $a \neq a' \in A$.  A maximal finite set of this kind can be found in
 $U(M_0)$ for any countable elementary submodel $M_0$ of $M$.  It follows that $\sU$ is separable, and  does not depend on $M$.   Being a compact metric space, $\sU$ is  second countable.
  
  Let $\sX$ be the space of weakly random $\phi$-types.
    
   We have a map 
 $i: U \to L^2(X)$, mapping $a$ to the characteristic function of $\phi(a,x)$.  We have $|| i(a)-i(b) ||_2^2 = d(a,b)$.
  The definable measure $\mu$ induces a regular Borel measure on $\sX$, and $i$ extends to an embedding
  $i: \sU \to L^2(\sX)$.  

  Let ${\sG}$ be the group of permutations of $\sU$ induced by automorphisms of $M$.  These are the permutations preserving all
  images in $\sU$ of $\bigwedge$-definable subsets of $U^n$.  In particular $d$ is preserved, so ${\sG}$ consists of isometries of $\sU$.
  We give $\sG$ the topology of pointwise convergence (equivalently here, the uniform topology.)
  It follows that ${\sG}$ too is second countable.   Since ${\sG}$ is closed in the group of all isometries of $\sU$, it is compact.  
  We will refer to ${\sG}=: \sG_{\mu,\phi}$ as the compact symmtery group of $\mu$ relative to $\phi,U$.
  
  The action of ${\sG}$ extends canonically to an action by automorphism on   $\sX$.  
   As $\sG$ preserves
  the relations giving the measure of any finite Boolean combination of $\phi(a,x)$, this action is measure-preserving.  It thus induces
  a unitary action of $\sG$ on $L^2(\sX)$.  Note also that $g \mapsto (gu,v)$ is continuous.  
  If an automorphism of $M$ fixes $\sU$, it fixes $\sX$ too, since it fixes the image of the relation $\phi$ on $\sU \times \sX$.   
 
  Conversely if $g$ fixes $\sX$, it is clear from the definition of $E$ that $g$ fixes $\sU$.

 \begin{thm} \label{nipsymmetry} Let $\phi \subset X \times U $ be a NIP formula, and assume $U$ is a complete Shelah strong type,
 i.e. a $\bigwedge$-definable set carrying no nontrivial  definable equivalence relations with finitely many classes.
 Let $\mu$ a  definable Keisler measure on $X$.

  Then  $\sG_{\mu,\phi}$ has finite archimedean rank.
 \end{thm}
 
 \prf
 Let $\sG= {\sG_{\mu,\phi}}$, and let $\sG^0$ be the connected component of $\sG$.  By the completeness assumption on $U$,
 $\sG^0$ acts transitively on $\sU$.  
 
By Peter-Weyl there is a cofinal system ${\mathcal{N}}$ of closed normal subgroups $N$ of $\sG^0$ such that $\sG^0/N$ is a connected Lie group.   By `cofinal' we mean that any neighborhood of the identity element of $\sG^0$ contains some $N \in \mathcal{N}$;
equivalently,  we have $\meet {\mathcal{N}} = (1)$.   
Since $G$ is second-countable, we may take ${\mathcal{N}}$ to be countable.  

For $N \in {\mathcal{N}}$, we can factor out $(\sU,\sX, i,\mu)$ by the action of $N$, to obtain 
$(\sU/N,\sX/N, i,\mu_N)$ where $\mu_N$ is the pushforward measure on $\sX/N$.  We will also use the orthogonal projection
$\pi_N:  L^2(\sX) \to L^2(\sX)^N$, where $L^2(\sX)^N$ is the   subspace of $L^2(\sX)$ consisting of $N$-invariant functions.
We have a canonical identification of $L^2(\sX)^N$ with $L^2(\sX/N)$.
By computing $\int fg $ for $g \in  L^2(\sX/N)$, it is easy to see that 
$\pi_N(f)$ is the integral of $f$ with respect to Haar measure on $N$, i.e. $\pi_N(f) (x) = \int f(n(x)) d_N(n)$.
Hence, for continuous $f$ on $\sX$, we have $f = \lim_{N} \pi_N(f)$, i.e. for any $\e>0$, for some open neighborhood $W$
of $1$ in $G$, whenever $N \subset W$, we have $||\pi_N(f) - f ||_\infty <\e$.

 The action of $G/N$ on $ L^2(\sX/N)$ is faithful:  if $g \in G \m N$, then by Urysson's lemma there existss a continuous function $f$ on $G/N$
 vanishing at $1_{G/N}$ but not on $g$; clearly $^g f \neq f$.

Define $i_N= \pi_{N} \circ i:  \sU \to L^2(\sX / N)$.  
  
  By \propref{ls4.6},  
   the image $i(\sU)$ of $\sU$ in $L^2(\sX)$ has finite packing dimension.
  By \lemref{lipschitz}, the packing dimension of $i_N(\sU)$ is at most $\delta$, for each $N \in \mathcal{N}$.  
  Since $\sG^0$ acts transitively on $U$,     the image of $\sU$ in $L^2(X/N)$ forms a single orbit of ${\sG^0}/N$-orbit.
   By \lemref{manifold3}, $\dim(\sG^0/N) \leq \delta^2$.  So   $\dim ({\sG^0}/N)$ is bounded, independently of $N \in \mathcal{N}$.   
   By \propref{plp} (4),   ${\sG^0}$ has finite archimedean rank, hence (by definition) so does $\sG$.    \eprf

\begin{rem} \rm The theorem is valid more generally for continuous logic; the measure $\mu$ is better interpreted as an expectation operator in this case, but
 by Riesz, induces a measure $\mu$ on $\sX$; 
 we still have an embedding of $U$ into $L^1(X,\mu)$ mapping $a$ to  the real-valued function $\phi(a,x)$; we use the induced norm from $L^1$.
 \end{rem}

\ssec{The local connected component  of a NIP group}

In this paragraph we address Pillay's question, discussed in the introduction.

  Let $G$ be a definable group.  Let $\phi(x,u)$ be a NIP formula, with $x$ ranging over $G$ and $u$ over a definable set $U$.
  
  Assume $G$ carries a definable, left translation invariant measure, 
  at least on the Boolean algebra of subsets generated by left translates of formulas $\phi(x,b)$.

  Recall the definition of a weakly random type, \secref{weaklyrandom}.  
Let $G^{00}_\phi$ be the  intersection of all the stabilizers of weakly random $\phi$-types.  Equivalently:
 \[  g \in G^{00}_\phi \iff   \mu( \phi(gx,b) \  \triangle\, \phi(x,b) )= 0   \hbox{ for all }    b  \]

More generally, if $q(u)$ is a partial type ,let 
\[ G^{00}_{\phi,q} = \{g \in G:    \mu( \phi(gx,b) \  \triangle\, \phi(x,b) )= 0   \hbox{ for all }    b \models q \}.\]

  Thus $G^{00}_{\phi}$   is a $\bigwedge$-definable subgroup of $G$ of bounded index, with compact quotient $K$. 
  
  When    $\phi(x,y) := (xy\inv \in P)$, $P$ being a definable subset of $G$,  we also write $G^{00}_P$.

We can view $G/G^{00}_\phi$ as the maximal quotient of $G/G^{00}$ attributable to $\phi$.  
 $G/G^{00}_{\phi,q}$  the maximal compact quotient of $G$ attributable to instances $\phi(x,b)$ with $b \models q$.

 \<{thm}   \label{anand-nip-groups}   
  
 For a definable set $P$, the  compact group $G/G^{00}_P$  
 has finite archimedean rank.   
 
 More generally, if $q$ is a complete Shelah strong type, then   $G/G^{00}_{\phi,q}$  has finite archimedean rank.

   \>{thm}

\prf  This reduces to \thmref{nipsymmetry} by a standard transposition from definable groups to semi-direct factors of automorphism groups.  Namely introduce a new sort $X$, with a regular action of $G$ on $X$ `on the right', 
and no additional structure beyond the original   structure on the original sorts.  We can view $X$ as another copy of $G$; the new language includes
the old and the map $X^2 \to G$, $(x,y) \mapsto x \inv y$.   
 Define $\phi$ on $X^2$:  $\phi(x,y) \iff P(x \inv y)$.  
The left-invariant definable Keisler measure $\mu$ induces
a definable measure $\mu_X$ on $X$.    

Let $\hat{G}$ be the compact symmetry group of $\mu_X$ relative to $\phi$, as defined in \secref{NIPformulas}.
 Now $G$ acts by automorphisms on $X$:  $(g,x) \mapsto g x$, fixing the old sorts including
$G$ itself.  Indeed this identifies $G$ with the automorphism group of the new structure over the old.  
 It is easy to check that  under this identification,  $G/G^{00}_P$ becomes a closed subgroup of $\hat{G}$.   Hence by \thmref{nipsymmetry}  and \propref{plp} (5),  $G/G^{00}_P$ has finite archimedean rank.  
 
 The more general statement is deduced in the same way from \thmref{nipsymmetry}, letting $U = X \times q$, and
 $\phi'(x, x',u) = \phi(x \inv y, u)$.  Note that $G$ acting on the right is transitive on $X$ while fixing $q$ pointwise, while for homogenous $M \models T$, $Aut(M)$ is transitive on $q$, so altogether the automorphism group of the new structure is transitive on $U$. 
  \eprf

 \<{example}   Take $G$ to be a saturated elementary extension of $(\Zz_p,+,D)$ where  $D = \{ x \in \Zz_p: v_p(x) \in 2 \Zz \}$.
Then  $G^0=G^{00} $.  The measure theoretic stabilizer of $D$ is $G^0$, and $G/G^0 \cong \Zz_p$. \>{example}

 \<{example}   \rm  Here we show that the Lie group can interact with the totally disconnected group in \thmref{anand-nip-groups}, so that  $G^0/G^{00}$ 
 is not a product of a totally disconnected group with a Lie group.  
 Note that $\Zz (1,1)$ is a discrete subgroup of the topological group $\Zz_p \times \Rr$; it intersects $\Zz_p \times (-1,1)$ trivially.   Let $G=(\Zz_p \times \Rr) / \Zz (1,1)$.  The pathwise connected component of the identity in $G$ is a dense subgroup, the image of $\Rr$, and so $G$ is connected. 
 The image of $\Zz_p$ is a compact normal subgroup of $G$, with $G/ \Zz_p  \cong  T:=  \Rr / \Zz $, so $G$ is profinite-by-Lie.
 
  It is easy to see that the subgroup generated by $(1,1)$ is relatively $p$-divisible in $(\Zz_p \times \Rr)$, and  $G$ has no nonzero $p$-torsion elements.  
 Thus the exact sequence $0 \to \Zz_p \to G \to \Rr/\Zz \to 0 $ does not split; even as a pure group, $G$ does not contain a copy of $T=\Rr/\Zz$;  and this will not be fixed by moving to a subgroup of finite index, or a finite quotient.
 
 As a set, $G$ can be identified with $\Zz_p \times [0,1)$.   The group law is then defined by 
  \[(z,t) + (z',t') = (z+z'+c(t,t') , m(t,t')) \]
  where $c(t,t') \in  \{0,1\} $, and  $0 \leq m(t,t') < 1$ in $\Rr$.      Thus the group $G$ is definable in the ring  $\Qq_p \times \Rr$.  
  (Equivalently, in the model-theoretic disjoint union of the valued field $\Qq_p$ and the ordered field $\Rr$.)
 The theory of this structure is NIP.

 The product measure of the two Haar measures on $\Zz_p$ and on $[0,1)$ (identified with $\Rr/\Zz$)  is invariant for this multiplication too.
So it gives a generically stable measure.    Let $D = \{(x,y): v(x) \in 2 \Zz, y \in [0,1/2) \}$ say; a definable set of positive measure.

Let $M^*$ be a saturated elementary extension.  Let $S$ be the $\mu$-stabilizer of D.   It's a  $\bigwedge$-definable subgroup of bounded index, and $G(M^*)/S \cong G$.  $G$ is connected and profinite-by-Lie, but cannot be split as a product.

\>{example}

 \ssec{Open questions} \label{open7}
 \begin{enumerate}
 \item   Is the transitivity assumption on $U$ necessary in \thmref{nipsymmetry}?    The original question was for {\em all} $\phi(x,b)$.   
  
 The question reduces to $U$ comprising finitely many types $U_1,\ldots,U_k$, provided the bound on $\dim(G)$ does not depend on $k$. The issue
 is that the action on each $U_i$ may have a kernel $N_i$.  For definable $N_i$,  the Baldwin-Saxl lemma would imply that the intersection of 
 $k_0$ of the $N_i$ is already trivial, where $k_0$ is bounded in terms of the VC dimension.   
 Since $N_i$ is only $\bigwedge$-definable,  Baldwin-Saxl does not directly apply and some substitute is needed.   In any case we obtain
 that the compact group $\sG$ embeds into a product of Lie groups of bounded dimension; in case $\mu$ takes only finitely many values, the $\mu$-stabilizer
 of each formula $\phi(x,b)$ is definable, and so the usual Baldwin-Saxl applies and shows transitivity is not needed.
 
 \item  \label{locallycompact}  A local logic version, allowing for non-compact definable groups over $\Rr$ and $\Qq_p$, would be very interesting.

 \item  In \cite{macpherson-tent}, Macpherson and Tent study profinite definable groups $G$ in a NIP structure $M$, along with a formula $\phi(x,u)$
such that any open subgroup has the form $\phi(x,b)$ for some $b$.    \footnote{Actually their formulation is slightly stronger condition,   that the $\phi(x,b)$ for $b$ from $M$ define exactly the set of open subgroups.  However it is probably sufficient for the main theorem of \cite{macpherson-tent} to assume that the family contains all the open subgroups, as well as other groups that may have infinite index.}
In this situation the Haar measure yields a natural definable  measure on $\phi$-types, and 
the  `fullness' assumption implies that $G^{00} = G^{00}_\phi$ and profiniteness of $G(M)$ implies that  $G/G^{00} \cong G(M)$.  
They show that $G(M)$ is a finite product of finite dimensional $p$-adic analytic groups.  

This can be seen as a profinite / adelic analogue of Pillay's conjecture on the archimedean part of $G^{00}_\phi$, though under a much stronger hypothesis. 

   It would be very interesting to put these results on a common footing.  In particular, does \thmref{nipsymmetry}
  have an adelic analogue?   If $G$ is a pro-$p$-group, is it $p$-adic analytic?   Is $G$ in general interpretable in the adeles, or rather in the 
  disjoint union of finitely many $p$-minimal $p$-adic fields $\Qq_p$ and the o-minimal field $\Rr$ (all possibly enriched analytically)?

\end{enumerate}

\end{section}

 \appendix
\<{section}{Stability for invariant relations}
 
We develop   the basic results of stability, presented here   in Theorems \ref{stable-f} and \ref{stable-i}.  We view them as
a reduction, modulo a certain ideal, of binary relations to unary ones; thus a kind of measurability result for binary relations
for the product measure.  The theory is primarily due to Shelah, and for the most part we follow standard presentations.  
   Shelah understood the significance of having the theorem over an arbitrary base structure and not 
just over an elementary submodel, and introduced imaginary elements and the algebraic closure as the precise obstructions to this.   In \cite{lazy}, the theory was
extended beyond the first order setting.    In \cite{kim-pillay},   the main theorem was proved for arbitrary invariant stable relations {\em over a model}.  A little later, for simple theories, the ``bounded closure" with its compact automorphism
group was recognized by these authors as the obstacle to existence of 3-amalgmation.   This was the first use of Lascar's
compact Lascar group;  in the case of finite S1 rank, $3$-amalgamation was known to hold with  ordinary  algebraic closure 
and the associated profinite group ( in \cite{pac}.)     

See \cite{clpz} for a good presentation of the compact and general Lascar types; we will use it below.
In \cite{byu}, the theory was beautifully developed for continuous real-valued relations;   \ref{stable-i} is a (less elegant) generalization
for more general $\inft$-definable stable relations.  

Here we treat arbitrary automorphism-
invariant stable relations, over any base set.  We show that the fundamental theorems of stability theory hold, with strong Lascar types
as the natural obstacles to both uniqueness and existence.   
 
For  $\inft$-definable relations, generalizing slightly the 
continuous real-valued case, 
 we show that compact Lascar types or Kim-Pillay
types suffice.

We begin with describing a local setting, allowing notably to discuss stable independence over an ``imaginary" element of the form $a/E$, where $E$ is a $\bigvee$-definable
 equivalence relation.  We will need it in order to treat approximate equivalence relations canonically, in particular preserving any group actions on them. 
 This generalizes the usual setting if one takes the metric $d$ to be bounded.

 We will sometimes assume the language is countable;  the generalization to the general case (by considering countable sublanguages) is routine.  
 $A$ will denote a countable base set; we will  sometimes use a countable elementary submodel $M$ containing $A$.
{\em The unqualified words definable, $\bigvee$-definable, $\bigwedge$-definable always mean: without parameters.}

 When $R \subset X \times Y$, and $a \in X$, we let $R(a) = \{b: (a,b) \in R \}$.
 Define $R^t \subset Y \times X$, $R^t = \{(b,a): (a,b) \in R \}$.  When the context leaves no room for doubt, for $b \in Y$ we will write $R(b)$ for $R^t(b)$.

\ssec{Local structures}   \label{local-structures}

Let $\Uu$ be a structure with a metric $d: \Uu^2 \to \Nn$, such that for any $n$, $\{(x,y): d(x,y) \leq n \}$ is  $\bigvee$-definable.  Recall that a {\em definable} subset of $\Uu^n$ is the interpretation of a formula 
(without parameters); a $\bigvee$-definable set is a union of definable sets.  
  \footnote{  If many sorts are allowed, we still assume the domain of  $d$ is the set of all pairs,
belonging to  the union of all sorts.   There are natural generalizations to bigger semigroups than $\Nn$,
 both in the direction of continuous metrics and of uncountable languages, but we restrict here to the main case.}
A typical way to obtain such a structure is to begin with an arbitrary binary relation $R_0$ on another structure $\Uu_0$.
Let $\tE$ be the equivalence relation generated by $R_0$.  Then any $\tE$-class is naturally a local structure; the metric distance $d(x,y)$
is the length of a shortest chain $x=x_0,\ldots,x_n=y$ with $R(x_i,x_{i+1})$ or $R(x_{i+1},x_i)$ for each $i$.  Here the  small distance relations are 0-definable;
if one takes a family of relations instead, it is only $\bigvee$-definable.

  A {\em graph} is a  set $\Omega$ with a symmetric binary relation $R$.     We let 
 $R(a):=\{b: R(a,b) \}$. 
 We define the associated {\em metric}   $d_R(x,y)=\min n. \ (\exists x=x_0,\ldots,x_n=y, R(x_i,x_{i+1}))$.   
  The graph is connected if $d_R(x,y)$ is always defined.  Note  $R(a) \union \{a\}$ is the $d_R$ $1$-ball around $a$.  
  We define $R^k(a)$ to be the $d_R$ $k$-ball around $a$, i.e. $R^k$ is the  composition of $R$ with itself $k$ times.

 A relation $R(x_1,\ldots,x_n)$ is {\em local} if for each $i,j \leq n$,  
for some $m$,  $R(x_1,\ldots,x_n)$ implies  $d(x_i,x_j) < m$.    (For unary relations, this poses no constraint.)  
We will be  concerned only with local   relations.   
 We will say, when only local relations are allowed, that the structure is local.  (This is closely related to the Gaifman graph,   used in finite model theory, and to Gaifman's theorem on this subject.)   
 
 Note that one cannot freely introduce dummy variables;
if we wish to involve an additional variable $y$, it must be added along with a formula that ensures $d(x_i,y)<m$ for some $m$.
Geometrically this means we allow bounded products of the form $X \times_\delta Y = \{(x,y): x \in X, y \in Y, \psi(x,y)\}$
where $\psi$ implies $d(x,y) < \delta$ for some $\delta$.  
Formulas are formed by such controlled addition of dummy variables, conjunction, disjunction or difference of formulas with the same set of variables, projections  distance-bounded universal quantifiers, of the form:  $(\forall x)(d(x,y) \leq l \to \phi(x,y))$.

  If $E$ is a $\bigvee$-definable equivalence relation in a saturated structure,   then each $E$-class can be presented as a local structure;  the local structures setting will enable us  to speak about independence over an $E$-class (viewed as a (generalized) imaginary element of the base.)  We can present $E$ as having the form $d(x,y) < \infty$, where $d$ is a metric such that
  $d(x,y) \leq n$ is definable, for each $n$.  Then we can take the basic relations to be the $d$-bounded ones (this does not depend on the choice
  of $d$.)    

Any relation $R$ has local traces $R_{|l}$, the intersection of $R$ with distance- $\leq l$ between any pair of variables.  Note that $R$ can be recovered from the $R_{|l}$,
in the specific structure at hand; so that the automorphism group of the local structure obtained in this way is identical to the original one.

If a local structure $\Uu$ has a constant symbol, or more generally a nonempty bounded definable set $D$,   then it can be viewed as a $\bigvee$-definable structure; it is the the union of the definable sets of points at distance $\leq n$ from $D$,
each of these being $0$-definable.  
In general however,  the automorphism group here need not respect any specific  inductive presentation.

The metric can   be extended to imaginary sorts; first to $\Uu^n$ via:  $d((x_1,\ldots,x_n), (y_1,\ldots,y_m)) = \max(\max_{i} \min_j d(x_i,y_j) ,\max_{j} \min_i d(x_i,y_j) )$   ; then to a quotient by a bounded
equivalence relation, with quotient map $\pi: \Uu^n \to \Uu^n/E$, with distance defined by   $d(u,v) = \inf \{d(x,y): \pi(x)=u,\pi(y)=v\}$.

We assume $\Uu$ is saturated as a local structure, i.e.  
 any $d$-ball is saturated; equivalently any small family of definable sets has nonempty intersection, provided the family includes a bounded set,
and that any finite subset   has nonempty intersection.   Local saturation can be achieved by taking an ultrapower using bounded functions only. 

 A remark on ultraproducts:
  if $(N_i,d_i)$ are a family of local structures for the same language, and $(N,d)$ is an ultraproduct in the usual sense,
one has an equivalence relation:  $d(x,y) \leq n$ for some standard $n$; each equivalence class is a local structure, and Los's theorem holds.
thus an ultraproduct here requires a choice of an ultrafilter along with a component, rather than just an ultrafilter.

\ssec{Locally compact Lascar types}
An $\inft$-definable relation $E= \meet_i E_i$ is a {\em cobounded equivalence relation} if
in any elementary extension $N$ of the given structure, $\meet_i E_i(N)$ is an equivalence relation $E(N)$,
and $N/E(N)$ has cardinality bounded independently of $N$.

Call a sort $S$ {\em separated} if it carries an $\inft$-definable cobounded  local equivalence relation.   

If $S$ is separated, 
  let $\eqlc^S$ be the   intersection of all    $\inft$-definable  cobounded  local equivalence relations on $S$.  Then $\eqlc^S$ is the unique smallest such relation.  It may change if we add parameters to the language. 
  If the identity of $S$ is clear we write simply $\eqlc$.

Let $\pi=\pi^{lc}_S : S \to S/\eqlc$ be the quotient map.  On $S/\eqlc$ 
 we define a topology:  $Y$ is closed iff $\pi \inv Y $ is locally $\inft$-definable. 
 
\<{lem} The quotient by $\eqlc$ is a locally compact space (and $\si$-compact) space.     \>{lem}
\prf  (Cf. \cite{simple}, \cite{clpz} for the bounded case, of Kim-Pillay spaces.)    Let $a \in S$, and let $B_n$ be the ball of radius $n$ around $a$, in $S$.
Then $\pi(B_n)$ is compact (so $S/\eqlc$ is $\si$-compact.)  
Since $\eqlc$ is local, say $d(x,y) < m$ for $(x,y) \in S^2$ with $x \eqlc y$.    
Then $\pi(a) \notin \pi(S \m B_m)$.  But $\pi(B_m) \union \pi(S \m B_m) = \pi(S)$.
Thus the compact set $\pi(B_m)$ contains   a neighborhood of $\pi(a)$, namely 
  the complement of $\pi(S \m B_{m})$.   The proof of the Hausdorff property is similar.
  
  \eprf 
  \<{rem}  
The local algebraic closure $\acl(\emptyset)$ in a given sort $S$ can be defined as 
the union of the locally finite definable sets.  The automorphism group of $\Uu $ has a quotient group acting faithfully
 on  $\acl(\emptyset)$, referred to as locally profinite Galois group of $S$.   It is 
 a totally disconnected locally compact group.   The stabilizer of a nonempty subset of  $\acl(\emptyset)$ is a compact group (fixing one point implies
leaving invariant balls of various radii.)     

 One can similarly define the {\em local compact closure} to be the union of $S/\eqlc$, over all sorts
$S$ such that $\eqlc$ is defined.    
\>{rem}

Now  consider the more general setting of   $Aut(\Uu)$- invariant equivalence relations (we will simply say: invariant to mean $Aut(\Uu)$-invariant.)  Assume $S$ has an 
$Aut(\Uu)$- invariant  cobounded local equivalence relation.  Then it has a smallest one; it is denoted $\eqlas$.   
This equivalence relation is generated by $\union_m \theta_m(a,b)$, where $\theta_m(a,b)$ holds iff $a,b$ begin an indiscernible sequence, and $d(a,b) \leq m$.
When $d$ has the property that any two elements are connected by a chain of elements of distance $1$, as is the case in the main examples, $\eqlas$ is generated by $\theta_2$.   At any rate,  $\eqlas$ is an $F_\si$ relation (a countable union of $\inft$-definable relations.)

\ssec{} \lbl{ }
 
Assume:

 ($\clubsuit$):  for some $m_0$, for all $n$, any $n$-ball is a finite union of $m_0$-balls.

\noindent{\bf Remark.}  
$\clubsuit$  is true in the setting of a measure, finite on balls.  More precisely assume $\mu$ is a  definable measure, each ball of radius $1$ has nonzero measure, and each ball of radius $\leq 3$ has finite measure.  Then by Rusza's trick, any ball of radius $3$ is a union of finitely many balls of radius $2$ (consider a maximal disjoint set of radius $1$- balls in
the radius $3$ ball; then enlarging them to radius $2$ would cover the larger ball.)    Assume in addition that the metric space is ``geodesic" in the sense that any two points of length $n$ are joined
by a path of length $n$, where the successive distance is $1$ (as is the case for Gaifman graphs.)  Then it follows inductively that any ball of radius $n$ is a union of finitely many balls of radius $2$.

\noindent{\bf Remark.}  
We are interested only in types of elements at finite distance from elements of $\Uu$.   
In the presence of $\clubsuit$, any such type has bounded distance $\leq m_0$ from some element of $\Uu$.  It follows
that if $\Xx$ is an   $Aut(\Uu)$- invariant  closed set of types over $\Uu$, then $\Xx$ contains a compact subset $X$ with $Aut(\Uu) X = \Xx$
(namely the types of distance $\leq m_0$ from a given point.)

\ssec{Ideals of definable sets}

We will work with   saturated (local) structures $\Uu$.
  {\em Invariance} refers to the action of $Aut(\Uu)$, or $\Aut(\Uu/A)$ for a small substructure $A$.  A set {\em divides} if for some $l$ it has an arbitrarily large set of $l$-wise disjoint conjugates
(i.e. any $l$ have empty intersection).

We will consider ideals of $\Uu$-definable sets (of some sort $S$).  Say $I$ is {\em definably generated} if it is generated by a definable family of definable sets.  
Say $I$ is $\bigvee$-definable if it is generated by some bounded  family of definably generated ideals.    
 
Equivalently, for any formula definable 
$D \subset S \times S'$, 
$\{b \in S':  S(b) \in I \}$ is $\bigvee$-definable.  If $I$ is $Aut(\Uu/A)$-invariant, then $\{b \in S':  S(b) \in I \}$ is in fact $\bigvee$-definable over $A$.

Dually, $I$ determines a partial type over $\Uu$, generated by the complements of the definable sets in $I$.  Any extension of this partial type is
called {\em $I$-wide}.   
We say $a/A$ is $I$-wide if $a$ does not lie in any $A$-definable set lying in $I$.  Note that $tp(a/A)$ will then extend to an $I$- wide complete type over $\Uu$.

If $f: S \to S'$ is a 0-definable surjective map, and $I$ is a $\bigvee$-definable ideal, let $f_*I = \{D: f \inv D \in I \}$.  This is a 
$\bigvee$-definable ideal on $S'$, proper if $I$ is proper.  If $c/A$ is $I'$-wide, then $c=f(b)$ for some $I$-wide $b/A$.

If $I,I'$ are two ideals (on $S,S'$), we can define an ideal $I \tensor I'$ on $S \times S'$, generated by the sets $D \subset S \times S'$ such that for some $D_1 \in I$, for all $a \in S \m D_1$, $D(a) \in I'$.   So if $a/A$ is $I$-wide and $b/A(a)$ is $I'$-wide, then $(a,b)/A$ is $I \tensor I'$-wide.      Conversely,
if $(a,b)/A$ is $I \tensor I'$-wide, then $a/A$ is $I$-wide, and - assuming $I'$ is $\bigvee$-definable - $b/A(a)$ is $I'$-wide:  to see the last statement, if $b \in D(a) \in I'$, then since $I'$
is $\bigvee$-definable, there exists $\theta(x)$ true of $a$ such that $D(a') \in I'$ for all $a' \in \theta$; let $D' = \{(a',b'): b' \in D(a'), a' \in \theta \}$;
then $D' \in I \tensor I'$; and $(a,b) \in D'$.  

Inductively, we define $I^{\tensor n}$, $I^{\tensor (n+1)} = I^{\tensor n} \tensor I$.  We will say $b=(b_1,\ldots,b_n)$ is $I$-wide if it is $I^{\tensor n}$-wide.

\begin{example}  \rm
Let us mention here some canonical ideals,   relative to a given complete type $p$.  There is Shelah's  forking ideal $I_{sh} $, generated by
the set $Div(p)$ of  formulas that divide (over $\emptyset$).  Given any invariant measure $\mu$ (such that $p$ is wide), we have the ideal $I_{\mu}$ of all formulas
of $\mu$-measure zero.      If $\mu$ is definable, then $I_\mu$ is $\inft$-definable.  We have
$Div(p) \subseteq I_{sh} \subseteq I_\mu$, for any   invariant measure $\mu$.
\end{example}
 
 If $I$ is an ideal on $S'$,
 let $SDiv(I)$ be the family of generically $I$-dividing subsets of $S$; i.e. the family of sets $Q(b)$, $b \in S'$, $Q$ an $A$-definable subset of $S \times S'$, 
such that for some $n$, for any $I^{\tensor n}$-wide $(b_1,\ldots,b_n)  $
with $tp(b/A)=tp(b_i/A)$, $ \meet _{i=1}^n Q(b_i) = \emptyset$.  Note that if $I \subset J$ then $I ^{\tensor n} \subset J^{\tensor n}$, so $SDiv(I) \subset SDiv(J)$.
   
Let $\widehat{I}$  be the ideal generated by $SDiv(I)$.  We have $SDiv(I) \subseteq Div$ and so $\widehat{I} \subseteq I_{sh}$.
If  $I$ is $\bigvee$-definable over $A$, so are $SDiv(I)$ and  $\widehat{I}$.

\<{defn}  \lbl{strong-aa}
Let $R \subset P \times P'$ be an invariant relation over $A$, and let $I$ be a $\bigvee$-definable ideal on $P$.   
Say $R$ holds  {\em $I$-almost always} if for any $c \in P'$, for any $b \in P$ with $b/A(c)$ $I$-wide, 
we have $R(b,c)$.  Say $R$ holds {\em $I$-almost always in the strong sense} on $P \times P'$
if the transpose $R^t=\{(y,x): (x,y) \in R\}$ holds  $\widehat{I}$ almost always.   \>{defn}

 Explicitly,    $R$ holds {\em $I$-almost always in the strong sense} on $P \times P'$ if whenever $(b,c) \in P \times P' \m R$,
there exists an $A$- definable local $Q \subset P \times P'$ and $n \in \Nn$ such that $(b,c) \in Q$, and 
for any $I^{\tensor n}$-wide $n$-tuple $(b_1,\ldots,b_n)$, $P \meet \meet _{i=1}^n Q(b_i) = \emptyset$. 

If   $R \subset S \times S'$ is an invariant relation,  $I$   a $\bigvee$-definable ideal on $S$,  and $P \subseteq S, P' \subseteq S'$   invariant sets,
we will also say that  $R$ holds {\em $I$-almost always in the strong sense on $P \times P'$} if $R \meet (P \times P')$ does.  

\<{lem} \lbl{aa} Assume $R$ holds $I$-almost always in the strong sense on $S \times S'$.  Then: \<{enumerate}
\item  $R$  holds $I$-almost always. 
\item  If $tp(c/A(b))$ does not divide over $A$, and $tp(b/A)$ is $I$-wide, then $R(b,c)$.
\>{enumerate}  
\>{lem}

\prf  (1)  Suppose not; let $Q,n$ be as in \defref{strong-aa}.   Let $b_1=b$.  Inductively find $b_k$ such that $Q(b_k,c)$ and $b_k$ is wide over $A(c,b_1,\ldots,b_{k-1})$;
this is possible since $Q(c)$ is wide.  But then $c \in \meet_{i=1}^n R(b_i)$, a contradiction.

(2).  Supose $\neg R(b,c)$.   Let $Q$ be a definable set as in \defref{strong-aa}, 
so that  for any $I$-wide $(b_1,\ldots,b_n) \in S^n$, $\meet _{i=1}^n Q(b_i) = \emptyset$. 
As $tp(b/A)$ is $I$-wide, one can find $b_i \models tp(b/A)$ for $i \in \Nn$, such 
that $tp(b_n/A(b_1,\ldots,b_{n-1}))$ is wide.   Then any subsequence of length $n$ of this infinite sequence is $I^{\tensor n}$-wide,
so the intersection of $Q(b_i)$ over any such subsequence is empty.  It follows that $tp(c/A(b))$ divides over $A$.  
\eprf

 \ssec{Stable invariant local relations}

  \<{defn} Two definable relations $P(x,y),Q(x,y)$ are {\em stably separated} if   there is no   sequence of pairs $(a_i,b_i): i \in \Nn$ with $P(a_i,b_j)$ and $Q(a_j,b_i)$
 for $i<j \in \Nn$.  \>{defn}
 
 Let $R \subset S \times S'$ be an  $Aut(\Uu/A)$-invariant relation.   
 \<{defn}   
  $R$ is  {\em stable} if whenever $(a,b) \in R$ and $(c,d) \in (S \times S') \m R$, then there exist $A$-definable sets $Q,Q'$
   such that $Q(a,b), Q'(c,d)$ and $Q,Q'$ are stably separated.  \>{defn}

\<{rem} $R$ is stable iff  there is no indiscernible sequence $(x_i,y_i)$ such that for $i \neq j$, 
   $R(x_i,y_j)$ iff $i<j$.  \>{rem}

 \prf   If no such indiscernible sequence exists, then whenever $(a,b) \in R$ and $(c,d) \in (S \times S') \m R$, $tp(a,b)$ and $tp(c,d)$ 
 must be stably separated; by compactness, for some definable $P$ approximating $tp(a,b)$ and $Q$ approximating $tp(c,d)$,
 $P,Q$ are stably separated.   Conversely if $(a_i,b_i)$ is an indiscernible sequence as in the remark, then $tp(a_1,b_2)$
 is not stably separated from $tp(a_2,b_1)$ though $R(a_2,b_1)$ and $\neg R(a_1,b_2)$.
 \eprf

 \<{thm}  \lbl{stable-f} Let $\Uu$ be a   local structure, with $\clubsuit$.  Let $\ff$ be a   family of invariant stable local relations on $S \times S'$.  
 Let $E_\ff$ be the intersection of all  
   
  co-bounded  invariant local equivalence relations on $S$, 
 such that each class is a Boolean combination of a bounded
 number of sets $R(b) \subset S$, $R \in \ff$.   Then for each complete type $\bP$ in $S$,  there exists a proper, $\bigvee$-definable  ideal $I(\bP)$ on $S$, satisfying:
 
 (*) If $R \in \ff$, $P \subset \bP$ is an $E_\ff$-class, and $Q$ is an $E_{\ff^t}$-class on $S'$, then either $R$ holds  almost always in the strong sense for $I(\bP)$ on $P \times Q$, or $\neg R$ does.

 Also, symmetry holds:   if  
 for $P,Q$ as above, if $\bQ$ is a complete type with $Q \subset \bQ$, 
  then on $P \times Q$, $R$ holds almost always for $I(\bP)$ iff $R^t$ holds almost always for $I^t(\bQ)$.
  
 \>{thm}

\<{rem} \lbl{loc0}  \rm $E_\ff$ has a distinguished class $S^-$, such that  for any  $R \in \ff$, $\neg R$ holds almost always on $S^- \times S'$ in the strong sense for $I_S$. 
Away from this class, $E_\ff$ is a local relation.    (See proof,   above \lemref{symm}.)

\>{rem}

We remark that there also   exists a canonical proper $\bigvee$-definable ideal $I_S$, such that the dichotomy (*) and symmetry hold $I_S$-almost always.
However it may trivialize certain types on $S$.

Though the proofs go through for any $\ff$, we will
assume below that $\ff=\{R\}$ to simplify notation.  (In fact the theorem reduces easily to the case that $\ff$ is finite; and then, - replacing $S$ by $S \times \ff$, and considering
the relation $\widehat{R} ((x,R),y) \iff R(x,y)$ -  to the case that $\ff$ has a single element $R$.)

We will use the space $S_D(\Uu)$  of all bounded global types on a sort $D$, i.e. types containing a formula implying $d(x,a) \leq n$ for some $a, n$.  If $x$ is a variable of sort $D$, we will also write $S_x(\Uu)$.
  Let $(d_p x) R = \{b: R(x,b) \in p \}$.  If $(d_p x) R= (d_{p'}x) R$,   we say $p,p'$ define the same $R$-type.   We do not define a topology on the set  of global $R$-types.

\begin{lem}   \lbl{s0}  Let $M$ be a countable model.  Let ${R'}(x,y),{R}(x,y)$ be   definable relations (of which at least one is local.)
Assume ${{R}'}(x,y)$ and ${R}(x,y)$ are stably separated.  Then for any type $p$ over $M$ there exists a finite Boolean combination $Y$ of sets ${R}(x,c_i)$ with $c_i \in M$,
such that $d_p y{{R}'} \implies Y$ while $Y, d_py {R}$ are disjoint.  \end{lem}

\prf      Let $c \models p |M$.  Define $a_n,b_n,c_n \in M$ recursively .  Given $c_1,\ldots,c_{n}$, the equivalence relation:  $\bigwedge_{i \leq n} {R}(x,c_i) \iff {R}(x',c_i) $ has at most $2^{2n}$ classes; if none of these classes meets both $d_p y {R}$ and $d_p y {R}'$,
then some union $Y$ of these classes contains $d_p y {R}$ and is disjoint from $d_p y {R}'$, and the lemma is proved.  Otherwise, choose 
 $a_n,b_n$  such that $d_py{R}(a_n)$, $d_py{{R}'}(b_n)$,
while $a_n$,$b_n$ lie in the same sets ${R}(x,c_i)$, $i \leq n$.   Then, find $c_{n+1}$ such that ${{R}'}(d,c_{n+1}) \iff {{R}'}(d,c)$, where $d \in \{a_i,b_i: i \leq n\}$.
 
For $n<k$ we have ${{R}'}(b_n,c_k)$.  
Applying Ramsey with respect to the question ${R}$ and refining the sequence $(a_n,b_n,c_n)$,
we may assume that $ {R}(b_n,c_k)$ for all $n>k$ or for no $n>k$; but the former is impossible since ${{R}'},{R}$ are stably separated.  So $ \neg {R}(b_n,c_k)$ for all $n>k$

Since $a_n,b_n$ have the same ${R}$-type over the smaller $c_i$, it follows that $\neg {R}(a_n,c_k)$ for $n>k$.  
But for $n<k$ we have    ${{R}'}(a_n,c_k)$;
so the sequence $(a_n,c_n)$ contradicts the stable separation of ${{R}'},{R}$.   

\eprf

\<{cor}  \lbl{noperfect-stsep}  Assume $L$ is countable.  Let ${{R}'},{R}$ be stably separated local definable relations on ${S} \times {S'}$.  There does not exist an uncountable set $W \subset S_x(M)$ such that for $p \neq p' \in W$, for some $b \in M$, ${{R}'}(x,b) \in p$ while ${R}(x,b) \in p'$.  \>{cor}
 \prf  Let $Y_p$ be an $M$-definable set such that $d_p{{R}'} \to Y_p \to \neg d_p {R}$ (\lemref{s0}).  There are only countably many choices for $Y_p$, so  
 there will be $p,p' \in W$ with $Y_p=Y_{p'}$.  Now if ${{R}'}(x,b)  \in p$ then $b \in Y_p=Y_{p'}$ so $\neg {R}(x,b) \in p'$.   \eprf 

It follows that there is no map $f$ from the full binary tree $2^{<\omega}$ into ${S'}$, such that for each branch $\eta \in 2^\omega$, 
$\bigwedge {{R}'}(x,f(\eta | n+1): \eta(n)=0) \wedge \bigwedge {R}(x,f(\eta | n+1): \eta(n)=1)$ is consistent.    By compactness, for some finite $n$, no such map exists for the height-$n$ tree $2^n$.   
We define the {\em rank} of a partial type $W$ to be the maximum $m$ such that there exists $f: 2^m \to {S'}$, with 
$W \wedge \bigwedge {{R}'}(x,f(\eta | n+1): \eta(n)=0) \wedge \bigwedge {R}(x,f(\eta | n+1): \eta(n)=1)$ consistent for each $\eta \in 2^m$.

  Let $R$ be a stable invariant relation on ${S} \times {S'}$.   
  
\<{lem}\lbl{st5}  Let $p,{q}$ be types over $\Uu$.  Assume:  for
any stably separated local definable $\phi,\psi$, for some ${e}={e}_{\phi,\psi}$ we have:  ${e} \subset p,{q}$ and  $rk_{\phi,\psi}(p) = rk_{\phi,\psi}({e})=rk_{\phi,\psi}({q})$.
  Then $p |R = {q}|R$.  \>{lem}
 
 \prf Let $c \models p$ and ${d} \models {q}$.  Suppose $p |R \neq {q}|R$.  Then for some $b \in \Uu$,
 $tp(b,c)$ implies $R$ but $tp(b,{d})$ implies $\neg R$.  As $R$ is stable, $tp(b,c)$ and $tp(b,d)$ are stably separated; hence by compactness, some $\phi(x,y) \in tp(b,c)$ and $\psi(x,y) \in tp(b,d)$ are
 stably separated.  Let ${e}={e}_{\phi,\psi}$, $l=rk_{\phi,\psi}({e})$.   Let $[\phi(x,b)]$ be the set of types extending $\phi(x,b)$.   It follows that either $rk_{\phi,\psi}({e} \meet [\phi(b,x)]) <  l$ or $ rk_{\phi,\psi}({e} \meet [\psi(b,x)]) <  l$.  
 But $rk_{\phi,\psi}(p)=rk_{\phi,\psi}({q})=l$,  a contradiction.  
 \eprf
 
\noindent {\bf Remark} (Uniqueness of finitely satisfiable extensions).    
Thus  if ${e}$ is a partial type, $e \subseteq p,q$,  and   $rk_{\phi,\psi}(p) = rk_{\phi,\psi}({e})=rk_{\phi,\psi}({q})$ for all stably separated $(\phi,\psi)$, then for all stable invariant relations $R$ we have 
$p|R = {q}|R$.  This hypothesis holds if ${e}$ is a type over  a model $M$, and $p,{q}$ extend ${e}$ and are finitely satisfiable in $M$.

 \noindent{\bf Remark}  (Determination by Ind-definable part).  
We can also deduce that for any global $p,p'$, if $p'$ contains all schemes
\[\{ \psi(x,b): \theta(b) \}\]
that are contained in $p$, then for any stable invariant relations $R$,
 $p|R = p'|R$. 
For this, for each stably separated pair $(\phi,\psi)$, we look at the deepest $(\phi,\psi)$- binary tree {\em contained} in $p$ (rather than consistent with $p$.)  
\bigskip

For any partial type $Q $, we let $\widehat{Q}$ denote the set of types over $\Uu$ extending $Q$.  

\<{prop} \lbl{small}  Let $R$ be a stable local  invariant relation on $S \times S'$.  Assume $\clubsuit$.  
Let $X$ be a   nonempty closed  invariant subset of $\widehat{S}$.
Let  $X|R = \{(d_p x ) R: p \in X\}$.   

Then $1 \to 2 \to 3$:  \<{enumerate}  
\item  $X$ is minimal.
\item  for any stably separated $\phi,\psi$ defined over $A$, 
  $rk_{\phi,\psi}(p) $ is constant (does not depend on $p \in X$.)
 \item $X|R  $ has cardinality bounded independently of $\Uu$; in fact $|X|R| \leq \aleph_0^{|L|}$.
\>{enumerate}
Moreover,   a minimal nonempty closed invariant subset of $X$ exists.

\>{prop}

\prf   (1)   implies (2) since the set of elements of $X$ of $(\phi,\psi)$-rank $\geq n$ is a closed, invariant subset of $X$.

Now assume (2).  Fix $\phi,\psi$ stably separated, and say    $rk_{\phi,\psi}(p) =m$ for $p \in X$.  For each ball $B$ of the metric $d$, the intersection of $B$,$X$ and the complement of all definable sets of $(\phi,\psi)$-rank $\leq m$ is empty; by (local) compactness, $B \meet X$ is covered by finitely 
many definable sets of  $(\phi,\psi)$-rank $\leq m$.  Thus $X$  is covered by countably 
many such definable sets, say ${e}(\phi,\psi,l), l \in \Nn$.  Each $p$ now determines a function $\chi_p:  (\phi,\psi) \mapsto l$, where $l$ is least such that
$p \in {e}(\phi,\psi,l)$.  But in turn $p|R$ is determined by this function.  For if $p , p' \in X$ and $\chi_p=\chi_{p'}$, then by \lemref{st5}, $p|R = p'|R$.  
This proves (3).   

For the moreover, given a complete type $P$, let $ \widehat{P} $ be the set of types over $\Uu$ compatible with $P$.
Then $X$ meets some $  \widehat{P} $ nontrivially, so letting $Z = X \meet  \widehat{P} $ it suffices to show that   any nonempty closed invariant subset $Z$ of $ \widehat{P}  $ contains a minimal nonempty closed invariant subset.   
     Fix
   $b \in P$, and let $B$ be the ball defined by $d(x,b) \leq 2m_0$.  
  By $\clubsuit$, any type $p$ over $\Uu$ meets some $m_0$-ball; by saturation of $\Uu$, this $m_0$-ball contains a $P(\Uu)$-point $a$; so   $d(x,a) \leq 2m_0$ is
compatible with $p$.  By invariance, $d(x,b) \leq 2m_0$ is compatible with some $p' \in Z$.  Thus $\widehat{B} \meet Z \neq \emptyset$
(where $    \widehat{B}   $ is the set of all   types over $\Uu$ of elements of $B$.)   So if $Z_i$ is a descending chain of 
   nonempty closed  invariant subsets of $S_{{R}}^P(\Uu)$,   then $Z_i \meet \widehat{B}$ is
   nonempty, and as $\widehat{B}$ is compact, $\meet Z_i \meet \widehat{B}$ is
   nonempty, and in particular $\meet Z_i$ is nonempty.   Thus by Zorn's lemma a minimal element exists. 
  \eprf

 Let $S,S'$ be sorts, and $R \subset S\times S'$ be  invariant, stable.   

Let $Gen^R$ be the set of all restrictions $p|R$, where $p$ is a global type of $S$ and $p|R$ has a small orbit under $Aut(\Uu)$.  (The total number of orbits
is small, say by \lemref{st5}, so $Gen^R$ is small.)  When relativizing to a small set $A$, so $R$ is $Aut(\Uu/A)$-invariant, we write $Gen_A^R$.  
  
 Any type $P$ on $S$ extends to some element of ${Gen^R}$, 
by \propref{small}.     It follows that for any $\eqlas$-class $X$ on $S$  there exists an element $q_X$ of $Gen^R$ such that for any small $N$,
$q_X | N$ is realized in $X$.  Indeed {\em some} $\eqlas$-class of $P$ has this property; since all $\eqlas$-classes in $P$ are conjugate, all have it.
 
Similarly define ${^R Gen} = Gen^{R^t}$ on $S'$.

Define an equivalence relation $E_{\ff}$ on $S$ by:
$(a,b) \in E_{\ff}$ iff for all $p \in {^R Gen}$ and $R \in \ff$, $(d_py)R(a,y) \iff (d_py)(R_b,y)$; and dually define $E_{R^t}$ on $S'$.    $E_{\ff}$ is co-bounded since ${^R Gen}$ is bounded.  $E_{\ff}$ is local since $R$ is local:    if $a E_{\ff} b$ then for some $c$, $R(a,c)$ and $R(b,c)$; 
so $d(a,b) \leq d(a,c) + d(b,c)$.

We say that $q|R$ is consistent with an invariant set $Z$ if   any small subset $q_0$ of $q|R$ is realized by some element of $Z$.

\<{lem}[symmetry and uniqueness]  \lbl{symm} Any $E_{\ff}$-class  on $S$ is consistent with a unique $q \in {Gen^R}$.  If $q \in {Gen^R}, q' \in {^R Gen}$,  $a  \in S, a' \in S'$,
and $q$ is consistent with $E_{\ff}(a)$, and $q'$ with $E_{\ff^t}(a')$, then
$d_{q'} y R(a,y) \iff d_q x R(x,b)$.  
\>{lem}

\prf  We prove the symmetry statement first, following the standard route.   Suppose for contradiction that it fails for $q,q',a,a'$.  
Say $d_{q'} y R(a,y)$ holds but $dq _x R(x,b)$ fails.
Construct $a_n,a'_n$ so that  $a_n \models q| A(a'_i: i<n)$, $a_n E_{\ff} a$, and $a'_n \models q'| A(a_i: i<n)$, $a'_n E_{R^t} a'$.  Then since $a_n E_{\ff} a$,  $d_{q'} y R(a_n,y)$  holds, and similarly $d_q xR(x,a'_n)$ fails.  Thus if $i>n$ then $R(a_n,a'_i)$ holds but $R(a_i,a'_n)$ fails.
This contradicts the stability of $R$.

We have already shown that there exists $q' \in {^R Gen}$ consistent with the Lascar type $E_{\ff^t}(a')$.  Now if $q_1,q_2 \in {Gen^R}$ are both consistent with $E_{\ff}(a)$, then by symmetry we have $d_{q_1}x R(x,b) \iff d_{q'} y R(a,y) \iff  d_{q_2}x R(x,b) $.  Thus $q_1 = q_2$.
\eprf
 
Because of this lemma, if $\chi$ is an $E_{\ff}$-class and $q$ is the unique element of ${Gen^R}$ consistent with it,
we can write $(d_\chi x)R(x,y)$ for $(d_q x)R(x,y)$.

 Let $\chi $ be an $E_{\ff}$-class, consistent with $q$.  Let $M$ be a substructure such that for any two elements $q_1 \neq q_2 \in {Gen^R}$, there exists $b \in M$ with $R(x,b) \in q_1$ but $R(x,b) \notin q_2$, or vice versa.    Let $E^\ff_M$ be the equivalence relation:  $a E^\ff_M b$ iff
 for any $R \in \ff$ and   $b \in M$, $R(a,b) \iff R(a,b')$.  
 Then $\chi$ is a cobounded equivalence relation, each class is a bounded Boolean combination of sets $R^t(b)$, and $E^{\ff}_M$ refines $E_\ff$.
 Indeed by construction a unique element $q \in {Gen^R}$ will be consistent with a given $E^\ff_M$-class $\chi$.  So for any $q' \in {^R Gen}$,
 let $d$ be such that $tp(d/M)$ is consistent with $q'$; then 
 for $a \in \chi$, $R(a,y) \in q'$ iff $R(x,d) \in q$.

Since all $E_{\ff}$ classes of a complete type $P$ over $A$ are $Aut(\Uu/A)$- conjugate, it follows from uniqueness that all elements $q$ of $  {Gen^R}$ consistent with $P$ are $Aut(\Uu/A)$- conjugate.  

  We {\em choose} a minimal  nonempty closed $Aut(\Uu/A)$-invariant  set $X=X_P$ of 
global types extending $P$, as in \lemref{small}.  By this lemma, for any $\phi,\psi$, $\beta_p(\phi,\psi) = rk_{\phi,\psi}(p) $ does not depend on 
the choice of $p \in X$.    Let $I(P)=I(X_P) $ be the ideal generated by  all definable sets $D$ such that for some $\phi,\psi$,
$ rk_{\phi,\psi}(D) < \beta_p(\phi,\psi)$.

\<{lem}[dividing] \lbl{divi} Let $q'$ be a global type of elements of  $S'$,
Assume $q'|R^t \in {^R Gen}$, $P$ is an $E_{\ff}$-class, and $R(a,y) \in q'$ for $a \in P(\Uu)$.

For $i \in \omega_1$, let $b_i  \models q' | A(b_j: j<i)$.
Then for any $a \in P(\Uu)$, for cofinally many $\a \in \omega_1$ we have $R(a,b_\a)$.  
 \prf

Re-define $b_i$ (without changing the type of the sequence) as follows:  let $M_i \prec \Uu$ be a small model containing $a_j$ for $j<i$, and let
$b_i \models q' |M_i$.  Let $M = \union _{i < \omega_1} M_i$.  For any pair $(\phi,\psi)$, for some $i< \omega_1$, we have
$rk_{\phi,\psi} (tp(a/M_i)) =  rk_{\phi,\psi} (tp(a/M)) $.  Since $\omega_1$ has uncountable cofinality, 
for some $\a < \omega_1$, for any $\phi,\psi$ , 
$rk_{\phi,\psi} (tp(a/M_\a)) =  rk_{\phi,\psi} (tp(a/M)) $.   Since $M_\a \prec \Uu$, there exists a global type $q$ 
extending $tp(a/M_\a)$ such that $rk_{\phi,\psi} (tp(a/M_\a)) =  rk_{\phi,\psi} (q) $.  By \lemref{st5}, $q|R$ is uniquely determined.
On the other hand since $q'|R^t \in Gen_A^R(S')$, it is clear that $q' |R^t  \in Gen_M^R(S')$.  Since $R(a,y) \in q'$, by \lemref{symm},
$R(x,b) \in q$ if $tp(b/M_\a)$ is consistent with $q'$.   Hence $R(x,b_i) \in q$ for $i \geq \a$. But we can also construct a global type
$q^+$  extending $tp(a/M_{\a+1})$ with  $rk_{\phi,\psi} (tp(a/M_{\a+1})) =  rk_{\phi,\psi} (q^{+}) $.   As $rk_{\phi,\psi} (tp(a/M_{\a+1})) = rk_{\phi,\psi} (tp(a/M_{\a}))$,
it follows that $q=q^+$; as $R(x,b_\a) \in q$ we have $R(x,b_\a) \in q^+$, i.e. $R(a,b_\a)$.  
 
\eprf

\>{lem}   
 
 It follows from \lemref{divi} (as well as from \lemref{symm}, as we saw before)  that $(d_py)R(x,y)$ is a bounded (but infinitary) Boolean combination of instances of $R(x,b)$; namely $(d_py)R(a,y)$
 iff $R(a,b_j)$ holds for cofinally many $j$, where $(b_j)$ is a sufficiently long sequence as in the lemma.
 
\prf[Proof of \thmref{stable-f}]   We will use the equivalence relation $E_\ff$ and the ideals $I(P)$ defined above
\lemref{divi}.   We have to show
 
 (*) If $R \in \ff$, $P \subset \bP$ is an $E_\ff$-class, and $Q$ is an $E_{\ff^t}$-class on $S'$, then either $R$ holds  almost always in the strong sense for $I(P)$ on $P \times Q$, or $\neg R$ does.

Pick $p \in X(\bP)$, and $p' \in X(\bQ)$ (with respect to $^tR$.)  By definition of $E_\ff$, for any $a \in P$, $p'(y)$ implies $R(a,y)$, or 
else for any $a \in P$, $p'(y)$ implies $ R(a,y)$.  Without loss of generality the latter holds.  Now suppose $\neg R(c,b)$ holds with $c \in P, b \in Q$.
 As $p'(y)$ implies $R(a,y)$ and $E_{\ff^t}(a,c)$, $p'(y)$ also implies $R(c,y)$.  
Let $r=tp(c,b/A)$.  We have to show that the condition in \defref{strong-aa} holds, i.e. that 
 for some $n$, and some $D \in r$,  $\union D(x,y_j) \union \neg I_{\ff^t}^{\tensor n} (y_1,\ldots,y_n)$ is inconsistent.  
 Otherwise, there exists a sequence $c,b_1,b_2,\ldots$ with $b_k/ A(b_1,\ldots,b_{k-1})$ wide for $I_{\ff^t}$ for each $k$, 
 and $r(c,b_i)$ holds for each $i$.   
 Let $\si$ be an automorphism taking $(c,b)$ to $(c,b_1)$.  Then $q'=\si(p')$ is a global type, $q' |R^t \in GEN$, consistent with  
 $E_{\ff^{t}}$-class of   $\si(b_1)$, and $q'(y)$ implies $R(c,y)$ (since $\si(c)=c$.)  By   \lemref{divi}, $R(c,b_i)$ holds for some $i$.  But $r$ is a complete type, and cannot be consistent with both $\neg R(c,b)$ and $R(c,b_i)$.  This shows that $\union D(x,y_j) \union \neg I_{\ff^t}^{\tensor n} (y_1,\ldots,y_n)$ is indeed inconsistent.

We saw  that $(d_py)R(x,y)$ is a bounded Boolean combination of instances of $R(x,b)$; hence any  $E_\ff$ - class can be expressed as Boolean combination of a bounded
 number of sets $R(b) \subset S$, $R \in \ff$.   Given this, the finest co-bounded equivalence relation with this property refines $E_\ff$, 
 and so also satisfies (*).

\eprf

\<{rem}  Let $p(x,y)$ be a  partial type.  Then there exists a unique smallest stable  invariant relation $P$ containing $p$.  (I.e. $p$ implies $P$.)  $P$ is $F_\si$.   Likewise for `equational' in place of stable. \>{rem}

\prf   We prove the stable case; the equational case is the same, with $a_0=a,b_0=b$ below.
For any invariant relation $P(x,y)$, let $P'(a,b)$ hold iff there exists an indiscernible sequence of pairs $(a_i,b_i)$ with $a_1=a,b_0=b$, and $P(a_0,b_1)$.
Clearly $P'$ is $\bigwedge$-definable if $P$ is; and $P$ is stable iff $P=P'$.  Also if $P = \bigvee_j P_j$ then $P'= \bigvee_j P_j'$; and the operation $P \mapsto P'$ is monotone.  
So let $P_0=p$, $P_{n+1}=P_n'$ and $P= \union_{n \in \Nn} P_n$. Then $P$ is $F_\si$ and stable, and contained in any stable invariant relation containing $p$.  \eprf

Presumably  $P$ is usually not $\inft$-definable.  (For instance when $p$ implies $\eqlas$ and $\eqlas$ is not $\bigwedge$-definable.)
Note that $\eqlas$ is itself a stable invariant relation.

\ssec{$\inft$-definable stable relations}

 We will discuss a stable, $\inft$-definable relation $R(x,y)$;  the results go through in the same way for a set $\ff$ of such relations.  We assume for simplicity that the language $L$ and
base set $A$ are countable, so $R = \meet_n R_n$ for some sequence $R_n$ of definable relations, with $R_1 \supset R_2 \supset \cdots$;  the general case reduces immediately to this.   We work with a universal domain $\Uu$.

 First we note that
the $p$-definition of $R$ is $\bigwedge$-definable, for any type $p$.

\<{lem}[definability]  \lbl{s1} Let   $p \in S_x(M)$.    Let $R = \meet_n R_n$, with $R_n$ definable.  
Then \begin{enumerate}
\item  $d_p R$ is $\inft$-definable over $M$;  it  is an intersection of    Boolean combinations of   sets $R_{n}(c)$
with $c \in M$. 
\item 
  In fact for any 
 $m \in \Nn$    there exists  $n=n(m)$ and  a finite Boolean combination $Y$ of sets $R_{n}(x,c_i), c_i \in M$,   such that $d_pR \to   Y \to d_pR_m$.
\end{enumerate}
\>{lem}

\prf  It suffices to prove (2).  By  stability, there is no sequence $d_n,e_k$ with $\neg R_m(d_n,e_k)$ for $k>n$ and $R(d_n,e_k)$ for $k<n$.  By compactness, for some $n_0$, there is no sequence with 
$\neg R_m(d_n,e_k)$  for $k>n$ and $R_{n_0}(d_n,e_k)$ for $k<n<n_0$.  Thus $\neg R_m, R_{n_0}$ are stably separated.   By \lemref{s0}, there exists a
finite Boolean combination $Y$ of sets $R_{n_0}(x,c_i), c_i \in M$,   such that $d_pR_{n_0} \to Y \to d_pR_m$.  
 
\eprf

 \<{lem} \lbl{s2}    Any  $E_\ff$-class of elements of $P$   is  $\inft$-definable with parameters, on any complete type $P$.  It is cut out   by certain sets  of the form $(d_qy)R(x,y)$. 
 \>{lem} 
 
 \prf  
   Let $P$ be a complete type of ${S}$.

We can   find $a \in P$ such that $Q(a)= \{q \in {^R Gen}:  a \in (d_qy)R(x,y) \}$ is maximal, i.e. not properly contained in any $Q(a')$ (with $a' \in P$).    This uses Zorn's lemma, and the fact that
$(d_qy)R(x,y)$ is $\bigwedge$-definable, so if $(d_q y)R(a_i,y)$ and $tp(a_i/M)$ approaches  $tp(a/M)$
in the space of types over $M$, 
 then $(d_q y)R(a,y)$.

Let $Q=Q(a)$.  
Now $aE_\ff b$ iff for each $q \in Q$, $(d_q y ) R(b,y)$.  So the $E_\ff$-class of $a$ is $\bigwedge$-definable.

Since all $E_\ff$- classes in $P$ are conjugate,  all $E_\ff$-classes in $P$ are $\inft$-definable.  As $P$ was arbitrary, the lemma follows.  
 \eprf

\<{cor}  If $a \eqlc b $  then $(a,b) \in E_\ff$. \>{cor}

\prf   In any case $a \eqlc b $ implies that $a,b$ have the same complete type; so it suffices to show this  for $a,b \in P$, where $P$ is a complete type.

Define:  $aEb$ iff $tp(a/c)=tp(b/c)$ for any $E_\ff$-class $c$ (i.e. there exists an automorphism fixing $c$ and taking $a$ to $b$.)   Clearly $E \subset E_\ff$.
Let   $\{C_i: i \in I\}$ list all the classes.  
then $aEb$ iff for each $i$, $(\exists c)(\exists d)(c,d \in C_i \& ac \equiv bd)$.  Since each $C_i$ is $\inft$-definable by \lemref{s2}, $E$ is $\inft$-definable.  
Since the number of classes $C_i$ is bounded,  and elements with the same type over some representative $c_i \in C_i$ also have the same type
over $C_i$, it is clear that $E$ is cobounded.  Hence $\eqlc \subset E$, so $\eqlc \subset E_\ff$.
\eprf

  From this and \thmref{stable-f} we obtain:

\<{thm}[locally compact equivalence relation theorem]  \lbl{stable-i}  Let $\ff$ be a nonempty family of  $\bigwedge$-definable stable local relations on $S \times S'$.  
Assume $S'$ is a complete type.  
There exists a proper  $\bigvee$-definable
 ideal $I'$ of definable subsets of $S'$, such 
that   if   $R \in \ff$, and 
   $P,Q$ are classes of $\eqlc$ on $S,S'$ respectively,
  then $R$ holds almost always on $P \times Q$ in the strong sense for $I'$, or $\neg R$ does.   
   Symmetry holds as in \thmref{stable-f}.    Also, the analogue of  \remref{loc0} is valid.
      \>{thm}

     In particular,
  fix $a$ and assume   $tp(a/A)$ forms a single $\eqlc$-class; then for $b$ such that 
  $tp(a/Ab)$ or $tp(b/Aa)$ does not divide over
$A$,   the truth value of $R(a,b)$ depends only on $tp(b)$.

Note that in the case of a definable measure, the measure $0$ ideal is $\bigwedge$-definable and so in general properly contains    the ideal $I'$ we found here;
they coincide only when both are definable.

 \begin{cor}\label{2continuity}  Let $R=\meet_n R_n$ be a $\bigwedge$-definable stable local relations on $S \times S'$.  
Assume $S'$ is a complete type.  Let $P,Q$ be classes of $\eqlc$ on $S,S'$, and assume $R$ holds almost always on $P \times Q$,
as in \thmref{stable-i}.  Then for each $n$ there exists a neighborhood $U$ of $(P,Q)$ such that if $(P',Q') \in U$ then $R_n$ holds almost always on $P' \times Q'$.
\end{cor}  

\>{section}

\begin{section}{Over a model}
 
 The entire thrust of this paper is to give a {\em lightface} account of higher measure amalgamation, choosing no constants.

Here we record  the much better understood situation {\em over a model}  in a similar language.   The idea is not to study the correlations in detail, but simply to take  an elementary submodel $M_0$ as if it were completely known, and describe the situation almost everywhere `above $M_0$', relying on  the fact that anything that may happen with positive probability has already happened in $M_0$.

  \thmref{overamodel} is a model-theoretic version of the hypergraph Szemeredi (or quasirandomness) lemma. 
 The methods are essentially  those of  Theorem 5 of Towsner \cite{towsner}, and  the results of 
 Tao cited there.      The results are valid only over a model, and in addition, only 'almost everywhere';  they are blind to phenomena occurring on  measure zero sets of $n$-types, 
  and so cannot give a meaningful stationarity lemma valid for all types (or even for almost all $n$-tuples of $1$-types, as opposed to almost all $n$-types.)

 We assume here that $L$ is a countable language, $T$ a complete theory, $X,Y$ definable sets carrying
 definable measures $\mu_X,\mu_Y$.   Form the multiple integral measures, and assume  
 Fubini holds for the product measures on $X \times Y^n$, for each $n$.  
 Let $\phi(x,y)$ be a definable relation on $X \times Y$.

 \begin{lem} \label{wrfs}  
   Let $M$ be a countable model, and $\phi(x,y)$ a formula.  Let 
 \[ B(\phi) =\{tp(a/M):   \mu_Y \phi(a,y) > 0 ,  \bigwedge_{m \in M} \neg \phi(a,m) \} \]
 Then $B(\phi)$ has measure zero.
 \end{lem}
 
 \prf  Note that $B(\phi)$ is a Borel  subset of $S_x(M)$; in fact the intersection of an open set with a closed set.
 
  Fix $\phi$, and let  
  \[ B_\e  =\{tp(a/M):   \mu_Y \phi(a,y) \geq \e ,  \bigwedge_{m \in M} \neg \phi(a,m) \} \]
 So as $\e$ descends to $0$,  $B(\phi)$ is the increasing union of the sets $B_\e$, and   it suffices to  show that
 each closed set $B_\e$ has measure zero, or just that 
  $\mu_X(B_\e)< \e$.  Fix $\e>0$, and 
  let \[ X_\e=\{tp(a/M):  \mu_Y \phi(x,y) \geq \e \} \]
Let $n=n(\e)$ be large, so that $\mu_X(X_\e) (1-\e)^n < \e$, and set 
\[ W = \{(x,y_1,\ldots,y_n) \in X \times Y ^n:  \bigwedge_{i=1}^n \neg \phi(x,y_i ) \} \]
Let $\mu = \mu_X \tensor \mu_Y \tensor \cdots \tensor \mu_Y$.   
Clearly, $\mu  (W) \leq \mu_X(X) (1-\e)^n < \e $.   
Let \[Y_\e=\{y \in Y^n:    \mu_X \{x:   (x,y) \in W \} \geq \e \} \]
By Fubini, $Y_\e$ cannot have full measure in $Y^n$.  
 
So $Y':=Y \m Y_\e$ is not a null set.   
  
Since $\mu_X$ is a definable measure and  $M$ is a model, we have $Y'(M) \neq \emptyset$.  Thus 
 for some $m_1,\ldots,m_n \in M$, $\mu_X \{x:   (x,m_1,\ldots,m_n) \in W \}< \e$.
But $B_\e \subset \{x:   (x,m_1,\ldots,m_n) \in W \}$; so $\mu_X(B_\e)< \e$.  Letting $\e \to 0$ we see that $\mu_X( B(\phi)) = 0$.
 \eprf
 
 \begin{cor}\label{wrfs-cor}  For almost all types $tp(a/M)$ in $X$,    any weakly random type (in $Y$) over $Ma$   is finitely satisfiable in $M$.
 \end{cor}
 
 \prf   By \lemref{wrfs},  $B:=\union_\phi B(\phi)$ has measure zero (here $\phi$ ranges over all formulas $\phi(x,y)$ over $M$.)  Assume $tp(a/M) \notin B$.  Let  
 $tp(b/Ma)$ be weakly random.  Then for any formula $\phi(x,y) \in tp(a,b/M)$, we have  $\mu_Y \phi(a,y) >0$,
 by weak randomness.  Hence by definition of $B$, $\phi(a,m)$ holds for some $m \in M$.   \eprf 
 
 \begin{rem}\rm  An independent family $(E_a)$ of  finite equivalence relations may, in general, be definable; then the effect of
 $E_a$ cannot be accounted for before one is aware of the parameter $a$, and 
 one  cannot expect 
$3$-amalgamation to hold over $M \union \{a\}$, but only at best over $M \union bdd(a)$.  Thus $4$-amalgamation 
cannot hold over $M$, in general, if we attempt to amalgamate extensions that are not algebraically closed.
   
   \lemref{wrfs}  shows
nevertheless that for {\em almost all types}, amalgamation {\em is} possible; for $a$ realizing a  random type over $M$, the finitely many classes of $E_a$ will already be represented in $M$.
 \end{rem}

   \def\Lt{\mathcal{L}}
   
 It will be useful to state a (tautological) measure theoretic lemma on compatibility of conditional expectation with random fibers.  
 
 \begin{lem} \label{E-fiber}  Let $X \to Y \to Z$ be Borel   maps between standard Borel spaces, let $\mu_X$ be a   Borel probability
 on $X$, with pushforwards $\mu_Y$ on $Y$ and $\mu_Z$ on $Z$.     
 For $z \in Z$, let $Y_z$ be the fiber above $z$ and let $X_z$ be the fiber above the composed map $X \to Z$.
 Assume $\mu$ `disintegrates' as an integral over $Z$ of a Borel family of measures $\mu_z$ on the fibers $X_z$ (so $\mu_X= \int_Z \mu_z$.)  
 Let $\phi:X \to \Rr$ be a bounded Borel function, with expectation $E (\phi) $ on $Y$.    For an $L^1$-function $\psi$ on $X_z$, let $E_z(\psi)$
 denote the expectation on $Y_z$ with respect to $\nu_z$.   Let $\mu_{Y,z}$ be  the pushforward of $\mu_z$ to $Y_z$.   Then $\mu_Y= \int_{z \in Z} \mu_{Y,z}$; and 
  for $\mu_Z$-almost all points $z \in Z$,  
 we have 
 \[  E(\phi) | Y_z =_{ Y_z -a.e.}  E_z( \phi| X_z) \]
 
\end{lem}
   
   \prf  When $Z=\{0,1\}$, $Y=Y_0 \dot \union Y_1$, and (if $\mu_Z(Y_0)>0$, and pulling back $Y_0$ to $X_0 \subset X$) the statement is that
   $E(\phi)| Y_0 = E (\phi| X_0)$, which is clear.  The general case   follows by   approximation.  \footnote{Or in   Radon-Nikodym style:    using separability of $L^2$, or countable generation of the algebra,  it suffices to show for a Borel function $\psi$ on  $Y$
that   $ \int_{Y_z} \psi(y) E(\phi)  = \int_{Y_z}  \psi(y)  E_z( \phi| X_z)$, for almost all $z$.  
This in turn is equivalent to showing for any bounded Borel $\theta$ on $Z$ that 
$\int_Z \theta(z) \int_{Y_z} \psi(y) E(\phi) = \int_Z \theta(z)  \int_{Y_z}  \psi(y)  E_z( \phi| X_z)$.
  Now
$E(\theta(z)\psi(y) \phi ) = \theta(z) \psi(y) E(\phi) $, so the left hand side is just
$\int_X \theta(z) \psi(y) \phi(x)$.   Similarly $E_z(\phi|X_z) \theta(z) \psi(y)= E_z( \phi|X_z) \theta(z) \psi(y))$,
  so the right hand side is $\int_{z \in Z}  \int_{X_z} \phi|X_z) \theta(z) \psi(y)=  \int_X \theta(z) \psi(y) \phi(x)$ too.}
   \eprf

 \ssec{}  Let $\Lt$ be  a continuous logic language, and   $\mathcal{T}$  be a stable theory of $\Lt$.   We assume $\mathcal{T}$ 
  eliminates quantifiers and imaginaries.  
 
 Assume given further a  
 piecewise, partial interpretation $\St$ of $\Tt$  in $T$:   namely
 a family $\mathcal{F}$ of maps  from   sorts of $L$ to various sorts of $\Tt$, such that 
 \begin{enumerate}
 \item  If $f: X \to    Y$ and $g: X' \to Y'$ are in $\mathcal{F}$ then so is $f \times g: X \times X' \to Y \times Y'$; 
 \item  If  $f: X \to    Y$ lies in $\mathcal{F}$, and $g: Y \to Y'$ is a ($\bigwedge$)-definable map of $\Tt$, then
 $g \circ f$ lies in $\mathcal{F}$;
 \item  The   pullback of
 any $\Lt$-$\bigwedge$-definable subset of $\Tt$ under any $f \in \mathcal{F}$ is $L$-$\bigwedge$-definable.
 \end{enumerate}
 
 By {\em partial} we mean that the maps $f$ need not cover between them a full model of   $\mathcal{T}$, but perhaps only a substructure.

 When $N \models T$ is sufficiently saturated and homogeneous,  and $A \subset N$ is countable, we 
 have $dcl(A) = Fix Aut(N/A)$.  
 We denote by $\St(A)$ the   definable closure of $A$ within $\St$:  
 \[ \St(A) = \{f(a):  f \in \mathcal{F}, a  \in dcl(A) \} \]
 Also write $\St_M(A) := \St(M \union A)$.
   
When given a tuple $(a_1,\ldots,a_n)$, write $a[n]:= (a_1,\ldots,a_n)$, and $a[n-i]:= (a_1,\ldots,a_{i-1},a_{i+1},\ldots,a_n)$.

 We will use the fact that in a stable theory if $a$ forks with $c$, then some $\phi(x,c) \in tp(a/c)$ causes forking;
 i.e. $\phi(a,c)$ takes value $>0$, and for any $c'$, if $\phi(a,c')$ takes value $>0$  then $a,c'$ are not independent.
 (The forking is due to some $\psi(x,y)$; let $\phi_1(x,y,e)$ (with $e$ a paramter in  $bdd(0)$)  be the absolute value of the difference between $\psi(x,y)$
 and the $tp(a/bdd(0))$-definition of $\psi$; then quantify out $e$, taking an appropriate supremum.)

 \begin{lem} \label{ind-constants} Let $M \subset N \models T$, $(a_1,\ldots,a_n,c) \in N$, and assume $tp(c/M(a_1,\ldots,a_n))$ is 
 finitely satisfiable in $M$.  Then $\St_M(a[n])$ is independent from $\union_i \St_M(a[n-i],c) $ over  $\union_i \St_M(a[n-i]) $ .
 \end{lem} 
  
\prf   Let $b \in  \St_M(a[n])$, and let $\phi(y,u)$ be a formula of $\Lt$, where $u=(u_1,\ldots u_n)$.   Also let
  $d_i \in \St_M(a[n-i],c) $, suppose  $\phi(b,d_1,\ldots,d_n) >0 $, and that $\phi(b,u) >0 $    causes forking
over $\union_i \St_M(a[n-i])$.   
We may write $d_i = f_i( a[n-i],c)$, where $f_i$ is a $\bigwedge$-definable function.  Since
  $tp(c/M(a_1,\ldots,a_n))$ is  finitely satisfiable in $M$, there exists $c' \in M$ with
  $\phi( b, (f_i(a[n-i], c'))_i) > 0$; let $d_i'=  f_i(a[n-i], c'))$.   Then $tp_\phi(b / d_1',\cdots,d_n')$ forks
  over $\union_i \St_M(a[n-i])$.    But this is a contradiction since $d_i' \in \St_M(a[n-i])$.  
 \eprf   
  
 Compare  \cite{towsner}, Lemma 4.     Note that Fubini is not required here.
  
  \begin{defn}
 Let   $Q_k$ be a collection of $k$-types, closed under restrictions and permutations of variables, $Q= \union_n Q_k$.  
 Consider a downward-closed family $S$ of subsets of  $\{y_1,\ldots,y_N\}$, all containing some base set $s_0$ and with $|s \m s_0| \leq k$ for $s \in S$..

 We say that $Q$ is a $(\leq k,\infty)$-amalgamation family if for any such family $S$,   and any map $j:S \to Q$ compatible with restrictions, such that $j(u) \in Q_{|u|}$, the union $\union_{u \in S} j(u)$ is consistent, and in fact extends to an element of $Q_N$.      
  
 (This is equivalent to $l$-amalgamation for each $l \leq k+1$, over a base set, in the sense of \cite{cigha}.) 
      
       An $n$-tuple whose type is in $Q_n$ will be called  $Q$- independent.  
       \end{defn}
       
 Note that for $k>2$, the hypothesis is incompatible with the presence of a linear ordering.   However it holds in many simple theories; for instance pseudo-finite fields over a base $A$ such that
 definable and algebraic closure coincide over $A$.  (Such a base exists for `most' but not all completions of the theory of pseudo-finite fields.).  A more refined version taking algebraic closure of each node in the system is
 valid in all pseudo-finite fields, but we restrict ourselves to the simpler case.

 \begin{thm}
       [Stationarity] \label{stationarity}  Let $\mu$ be a definable measure on a sort $X$.  
  
 Let $Q$ be a $(\leq n-1,\infty)$ -amalgamation family.  
  
Let $\phi_j(y_1,\ldots,y_n,x)  $ be   formulas   with $y_j$ a dummy variable not mentioned in $\phi_j$, and let 
  $tp(a^1,\ldots,a^n) \in Q_n$.  Then the quantity
\[ \mu(\bigwedge_{j \leq n}  \phi_j(a,x)  )   \]
depends only on the $n$-tuple of $n-1$-subtypes of  $tp(a^1,\ldots,a^n) $ and not on the full $n$-type.
 \end{thm}   

  \prf  Write  $p \eqtn p'$ to mean that the two $n$-types agree on any restriction to $<n$ of the variables. 
  
  For the sake of readability  we take $n=3$ as a representative case, and write $(a,b,c)$ for $(a^1,a^2,a^3)$.  
    Suppose $(a,b,c)$ and $(a',b',c')$ are $Q$-independent elements, and $tp(a,b,c) \eqthree tp(a',b',c')$.  We have to prove that for formulas $\phi_j$ as above, 
 $\mu( \bigwedge_j  \phi_j(a,b,c,x) ) = \mu(\bigwedge_j  \phi_j(a',b',c;x) )$. 
 
We may assume here that $tp(a,b,c,a',b',c')$ is $Q$-independent (using amalgamation over $\emptyset$.)
  
  We will construct   $M=\{a_1,a_2,\ldots\}$ such that  (i) $tp(a/M,b,c)$ is finitely satisfiable in $M$, (ii) likewise for $tp(a'/Mb'c')$, 
  and  (iii)  $tp(abc/M) \eqthree  \tp( a'b'c'/M)$.   During the induction, at stage $n$, we let 
  $\bar{a} = (a_1,\ldots,a_n)$.   
  We will be concerned with $tp(a,b,c/ \bar{a})$ and $tp(a',b',c'/ \bar{a})$ (but not especially with the type of $a,b,c$ over $a',b',c'$.)

  Assume $M_n=\{a_1,\ldots,a_{n}\}$ have been found, 
   with $tp(abc/M_n) \eqthree \tp( a'b'c'/M_n)$.  If $n$ is odd we work towards (i), if
  even towards (ii).  Say $n$ is odd.   Then we need to find $d$ such that  $d,a,b,c,a',b',c',a_1,\ldots,a_n$ is $Q$-independent, and :\\
  (i)  $\tp(d,b,c,a_1,\ldots,a_n)  = \tp(a,b,c,a_1,\ldots,a_n)$  \\
  (iii)  $\tp(a, b,c, / d, a_1 \ldots a_n)  \eqthree \tp(a',b',c' / d,a_1 \ldots a_n) $.  
  
  To meet (i), we extend $tp(a/b,c,\bar{a})$ to a type $p(x,a,b,c,\bar{a}$ over $a,b,c,\bar{a}$, so that $p \in  Q$;   this will be $tp(d/a,b,c,\bar{a})$.
 Next (moving the elements $a',b',c'$ if needed, recalling we are concerned only with their type over $\bar{a}$) we determine a type $tp(d,a',b',c',\bar{a}$, so that $tp(a',b',c' / \bar{a}, d)= \eqthree tp(a,b,c / \bar{a}, d)$; then (iii) is satisfied.  
This is possible using the induction hypothesis and $(2,\infty)$-amalgamation (or $(2,3)$-amalgamation over $\bar{a},d$.)

  Thus $M$ can be constructed satisfying (i,ii,iii).  
  Now the result follows from \lemref{ind-constants}. 
  \eprf

   \begin{rem} \rm   The notion of measure stationarity (the conclusion of \thmref{stationarity}) arose in early work of Elad Levi on the definable higher Szemeredi lemma.   
  Levi observed that it would suffice for a definable version of  Gowers' proof of higher-dimensional Szemeredi.  
  For the case of pseudo-finite fields,  
   stationarity was eventually proved in the stronger quantitative form, see   \cite{elad-alexis}.

  But for general theories this  remains interesting.  
  \end{rem}

    \begin{question}  \rm If $Y$ also admits a definable measure commuting with $\mu$, and Fubini is assumed, does   stationarity 
    \ref{stationarity} imply a precise formula similar to  \thmref{overamodel}, and valid on a set of full measure?   It seems plausible that  this can be proved by double counting and using Cauchy-Schwarz, as in \thmref{ind-prob} (2).     

Also, the proof should extend assuming higher amalgamation holds only for systems of algebraically closed substructures.
\end{question}

\begin{thm} \label{overamodel} Let $\mu$ be a definable measure, with Fubini.  
  Let $M$ be a model.  Then the measure spaces $S_{x_1,\ldots,x_n}(M)$
form an independent system.  \end{thm}
 
 Equivalently, the associated measure algebras $L_{x_1,\ldots,x_n}(M)$, with the standard embeddings among them,  form an independent system in the usual sense of stability.  (See Problem \ref{higher}.)

 By a L\"owenheim-Skolem argument, we may (and will)  assume the language as well as $M$ are countable.

 For readability we will write  omit the variable letter $x$, writing 
$\phi(123)$ for $\phi(x_1,x_2,x_3)$, $\phi(124) $ for $\phi'(x_1,x_2,x_4)$, $L(123)$ for the measure algebra of formulas in $x_1,x_2,x_3$ over $M$,
$L(12,23,13)$ for the join of the measure algebras $L(ij)$ ($1 \leq i,j \leq 3$), $S(12,23,13)$ for the corresponding (measured) Stone spaces.

Further we let $E(\phi; 12,13,23)$ denote the conditional expectation of $\phi$ relative to  $L(12,23,13)$.

Over larger structures  $M(b)$, $M(c)$ and $M(bc)$,  we have the measured Boolean algebras $L(1b)$ of formulas in $x_1$ over $M(b)$,
and likewise $L(1c)$, $L(1bc)$, and $L(1b,1c)$ genreated by $L(1b) \union L(1c)$;   and the Stone spaces $S(1b)= S_{x_1}(M(b))$ and similarly $S(1c)$ and $S(1bc)$.

  We view   formulas $\phi$ as $\{0,1\}$ valued (so conjunction is the same as multiplication),  or more generally valued in a bounded interval of $\Rr$ (so multiplication is still defined.)

 We use an integral symbol to denote expectation when it is absolute and not conditional; the 
  integral can always be understood to be over the largest space around, such   $S(1234)$ when the variables are among $x_1,x_2,x_3,x_4$.  Sometimes we will nevertheless indicate the intended space
by  a subscript, e.g. $\int_{1}$ for the integral over $S(1)$.

   \def\Ephi{ \widehat{\phi}}
 
\prf  The case $n=4$ is representative.   We want then to prove independence of $L(123)$ from $L(124,134,234)$ over $L(12,13,14)$.  

    It suffices to to prove that 
\[ \int \phi(123)\phi'(124) \phi''(134) \phi'''(234) = \int E(\phi(123); 12,13,23)  \phi'(124) \phi''(134) \phi'''(234)  \]
this will then extend to all  bounded $L^1$-functions on $S_{123}(M)$ in place of   $\phi(123)$; having replaced $\phi(123)$ by
$E(\phi(123); 12,13,23)$, we can continue and do the same with $\phi'$, etc. 
Let $\Ephi = E(\phi(123); 12,13,23)$.

It suffices to prove that for any random triple $b,c,d$ (i.e. $\tp(bcd/M)$ is random), we have
$\int_1 \phi(1bc)\phi'(1bd) \phi''(1cd)  \phi'''(bcd) = \int_1 \Ephi(1bc)  \phi'(1bd) \phi''(1cd))  \phi'''(bcd) $;  
or, taking out the constant factor $\phi'''(bcd)$, that
\[ \int_1 \phi(1bc)\phi'(1bd) \phi''(1cd)   = \int_1 \Ephi(1bc)  \phi'(1bd) \phi''(1cd))   \]

By \lemref{wrfs} and \corref{wrfs-cor}, $\tp(a/Mbcd)$ is finitely satisfiable in $M$.  Thus \lemref{ind-constants} applies,
and shows that 
\[ \int_1 \phi(1bc)\phi'(1bd) \phi''(1cd)    = \int_1 E(\phi(1bc); L(1b)L(1c) )  \phi'(1bd) \phi''(1cd))    \]

Thus the equality in the following claim finishes the proof.  

\claim{} Let $\Ephi(123) = E(\phi(123); 12,13,23)$.   Then for random $\tp(bc/M)$ we have:   \[\Ephi(1bc) = E(\phi(1bc); L(1b)L(1c) ) \] 

\prf   Note that 
$S(1bc)$ can be identified with the fiber above $\tp(bc/M)$ of $S(123) \to S(23)$;and likewise for $S(1b,1c,bc)$.   Thus the claim follows
from \lemref{E-fiber}, applied to the maps $S(123) \to S(12,13,23) \to S(12)$ and the   fibers above $\tp(bc)$.  
\eprf
 
\eprf

  Note that the type partition is canonical, once the model is chosen, and approximated by partitions into definable sets.
 
  Since the theorem is only valid  over a model, it loses sight of possible symmetries.  But if a definable group $G$ acts and is $\mu$-measure preserving, $G(M)$ acts on the type spaces and we do have equivariance.

 One can deduce a version of the Hoover-Kallenberg higher dimensional de Finetti theorem  in a similar (or actually easier) way:

\begin{prop}[Hoover-Kallenberg]
 \label{hk} Let $\mu$ be a definable measure.   Let $\Nn_1+\Nn_2$ be the disjoint sum of two copies of $\Nn$ ordered by $\Nn_1 < \Nn_2$, and let 
$(a_i:i \in  \Nn_1+\Nn_2)$ be an indiscernible sequence.   Let $M:=\{a_i: i \in \Nn_1\}$.  For $u \subset \Nn=\Nn_2$ let  $S_u$ be the space of types 
  in variable $x$ over $M \union \{a_i: i \in \  u\}$, with measure induced by $\mu$.   Then the $S_u$ form an independent system
  of measure spaces.   \end{prop}   
  
 \prf  For $i \in \Nn$, that $\tp(a_i / M \union \{a_j: j >i\})$ is finitely satisfiable in $M$.  Hence \lemref{ind-constants} applies.
 
 \eprf

 (Is Fubini needed?)

\begin{rem} \rm [NIP]   Assume NIP, and work over a model.  Then a statement much stronger than \thmref{overamodel} holds:  let $\mu(x)$ be a definable measure
(with no Fubini self-commutation assumptions.)  Let $B_n$ be the Boolean algebra of formulas in variables $x_1,\ldots,x_n$,
and let $B_{{1}\choose n}$ be the subalgebra generated by formulas $\psi(x_i)$ in a single $x_i$ variable.  
Then  for every formula $\phi \in B_n$  and $\e>0$ there exists $\phi' \in _{{1}\choose n}$ with  $\mu(\phi \symd \phi')< \epsilon$.  Equivalently,   the 
induced inclusion of $\si$-additive measure algebras, up to the null ideals,  is an isomorphism.     For $n=2$ this  is proved in 
  \cite{hps} 1.7(1), and under slightly different assumptions (essentially Fubini)  as Theorem 4.1(a,b) of \cite{lovasz-szegedy}.
 (It is curious that while the two teams of authors were entirely unaware of the parallel work in another field,    the Arxiv submissions are two days apart.)  The case of arbitrary $n$ follows immediately by induction from the case $n=2$.   
Once one knows that the measure algebra $M(X \times Y)$   is generated by the the $M(X)$ and $M(Y)$,  it follows that the measure algebra $M(X \times Y \times Z)$ is gen. by $M(X \times Y) \union M(Z)$ and hence by $M(X) \union M(Y) \union M(Z)$.      With Fubini assumed, a strengthening of this, both quantitative and qualitative, especially for distal theories, is  obtained in   \cite{chernikov-starchenko}.   

All of these sources allow arbitrary parameters.   Using definability of the measures, they thus hold over a model.

\begin{question} Does the 
above strong stationarity for NIP theories, identifying   $B_n$ with $B_{{1}\choose n}$, hold over $bdd(0)$?  \rm Possibly a statement of this type may follow by the method of \ref{stationarity}, noting that only the $n=2$ case is needed and that $(2,3)$-amalgamation is obtained in \thmref{ind-prob}.
\end{question}

For an extraordinary generalization to higher-arity NIP,   see \cite{chernikov-towsner},  Cor. 6.10 or Cor. 11.4.    
  
\end{rem}

   \end{section}

\begin{section} 
{An example from mixing}
This appendix to \S 3 is  intended to illustrate the use of expectation quantifiers and the various version of the independence theorem  \ref{ind-prob}.

We look at the convolution of two real-valued functions $f,g$ on a group.  This is well-studied in connection with mixing, see 
 \cite{LST} and \cite{GH}; I learned about this from a minicourse by Itay Glazer and Emmanuel Breuillard in Oxford, in spring 2024.  
 We will give a simple stability-theoretic proof of a special case of Theorem 1 of \cite{LST},
 namely for groups $G(\Ff_q)$ (modulo center) where $G$ is a simply connected algebraic group. (Using ACFA in place of PF we could also cover Rees and Suzuki group, i.e. the bounded rank case of \cite{LST}.)
 \corref{sc-measures}  covers e.g. nilpotent groups, and was written in response to 
 a  question of  Glazer's for vector groups; he  independently proved the vector group case  by analytic means.
 
Here the simplest version of stationarity (\thmref{uniqueness}) will suffice.   
But \thmref{overamodel} would not do; it is valid for almost all pairs of types, but in the proof of \propref{convolution-constant} it is essential to use the same type twice.

 \rm Let $G$ be a group carrying a left-invariant definable measure $\mu$.     $G$ may include    additional relations  (of discrete or of continuous logic.)   
 We will use pure probability logic quantifiers.  Formally, such  quantifiers do not directly distinguish  the graph of $\cdot$ from $\emptyset$; 
rather we first define, at the quantifier-free level, binary relations such as   $g(t,x)= g(t \inv x)$; and only then apply probability quantifiers to such relations.

 $G$ may for instance be a compact group made discrete - taken with the discrete metric as a CL structure, and  with the Haar measure serving to interpret expectation quantifiers.    Assuming the basic relations are measurable,   it follows that all formulas obtained by continuous connectives and expectation quantifiers are also measurable.  
 Or $G$ may be an amenable group with a finitely additive invariant measure.
 But the main example to have in mind  will be an   ultraproduct of finite groups with their normalized counting measure $\mu$.  In this case $\mu$ is both left and right invariant.    (This will be used in (\ref{sc-measures}) to see that the convolution is well-defined and continuous on $L^1$.)

We begin with the observation that the convolution $f*g$ is definable:  

\[  f*g(x) = E_t f(t) g(t \inv x)  \]

We say two real-valued functions $h,h'$ are equal a.e. if $E_x( |h(x)-h'(x)| ) = 0$.
Let $w$ be an element of $G$; define $h^w(x)=h(wx)$.  

For $h$ a definable function $G \to \Rr$
 
define the stabilizer of $h$ to be $Stab(h) := \{w \in G:   h =_{a.e.} h^w \}$.  Then $Stab(h)$ is an  $\bigwedge$-definable subgroup of $G$.

\begin{prop} \label{convolution-constant} Let $f,g$ be definable functions into $\Rr$, and let $h=f*g$ be their convolution.   Then the stabilizer of $h$ is 
a $\bigwedge$-definable subgroup of $G$ of bounded index.  In particular if $G$ admits no nontrivial definable homomorphisms into compact groups,
then $h$ is constant a.e.
\end{prop}

\prf   Let $w$ be an element of $G$, and $h=f*g$ the convolution.  By left invariance, substituting $w\inv t$ for $t$, we have
\[ h^w(x) =  E_t (f(t) g(t \inv w x))  = E_t (f(w \inv t) g(t \inv x) )  \]
By \propref{stablei}, this is a (real-valued) stable relation between $w$ and $x$.

Let $M_0$ be a countable model (taking the bounded closure of $0$ will also work.)
By \thmref{uniqueness}, provided $tp(a/M_0,b)$ does not divide over $M_0$ via a stable formula, 
 the value of $h^b(a)$ depends only on $tp(a/M_0)$ and $tp(b/M_0)$.

 Suppose $tp(b/M_0)=tp(c/M_0)$, yet $E_x |h^b(x) - h^c(x)| \geq \e >0$;  by the remarks  following \defref{ind-def}, there exists $a$  such that    $|h^b(a) - h^c(a)| \geq \e$ and 
$tp(a/M_0(b,c))$ does not divide, in particular $\ind{a}{b,c}$ holds; a contradiction.   Hence   $E_x |h^b(x) - h^c(x)|=0$, so
\[ h^b(x)=_{a.e.} h^{c} (x) \]
It follows that $b\inv c \in Stab(h)$.  We have shown that  any two elements with the same type over $M_0$
lie in the same coset of $Stab(h)$; hence $Stab(h)$ has bounded index.
\eprf

Let $G_n$ be a family of groups, endowed with left and right translation invariant  finitely additive measures.   We say $G_n$ is a quasi-random family, in the sense of Gowers, if 
for each $d \in \Nn$, for all sufficiently large $n$,  $G_n$ has no $d$-dimensional representations.
It follows from    \cite{conant-h-p} , Theorem 1.1 that  if $G$ is any  nonprincipal ultraproduct
of the $G_n$, then $G$  admits no nontrivial definable homomorphism into a compact Lie group,
and hence by Peter-Weyl no nontrivial definable  homomorphism into any compact group.
We use this below:

\begin{cor}  \label{convolution-constant-quasirandom} Let $( G_n)$ be a quasi-random family.  Let $b>0$ and let  $f,g: G_n \to [-b,b] \subset \Rr$ be functions with $||f||_1=||g||_1=1$.    Then  
 $||f* g - 1 ||_1 \to 0$ as $n \to \infty$.
 \end{cor}
 
 \prf  Set
 $h=f*g$.     By \ref{convolution-constant}, in any ultraproduct $G$, the stabilizer of $h$ is all of $G$; so $h$ is constant a.e., and since $||h||_1=1$ the constant value is $1$.
 \eprf

\begin{rem}  \rm In the case of bounded rank families of finite simple groups, a simpler  proof of \corref {convolution-constant-quasirandom} can be given.  
Let $G$ be a definably simple group in the theory of pseudo-finite fields.    It follows
from \cite{H-PAC}, 7.8 that any elementary extension of $G$ is also simple; 
hence $G$ cannot have a $\bigwedge$-definable subgroup of bounded index {\em in any expansion to a bigger language.}    Thus by \propref{convolution-constant},   $f*g$  is constant a.e.
 
\end{rem}

\begin{defn} A connected algebraic group $G$ over $\Qq^a$ is called
{\em simply connected} if  there is no surjective homomorphism $\widetilde{G} \to G$  of algebraic groups over $\Qq^a$ with nontrivial finite kernel.     \end{defn}  This definition is usually found in the setting of semi-simple groups, but we extend it to all connected algebraic groups. Commutative algebraic groups  whose geometric points form a divisible group are clearly simply connected; in particular in characteristic zero. vector groups such as $\Gg_a^n$ are simply connected.

 \begin{cor}  \label{convolution-sc}
     Let $G$ be a simply connected algebraic group.     
     Let $f,g: G(\Ff_p) \to \Rr$ be uniformly definable (hence uniformly bounded),  $||f||_1=||g||_1=1$.  Then 
     $||f*g  - 1 ||_1 \to 0$ as $p \to \infty$.  (And similarly for prime powers, and for $G(\Ff_p)/H_p$ where
     $H_p$ is some uniformly definable normal subgroup of $G(\Ff_p)$.)  
\end{cor}

\prf   Let $h=f*g$.   It suffices to show that $F \models E_x(|h-1|)=0$ for any  ultraproduct $F$ of the finite fields $\Ff_p$.    We have
$||h||_1 =1$ so it suffices to show that $h$ is $G(F)$-invariant a.e.  This in turn follows from \propref{convolution-constant}, once we show $G$ has no proper $\bigwedge-$ definable subgroups of bounded index.

Now Theorem 8.5 of \cite{H-PAC} shows that $G$(F) has no definable subgroups of finite index.   Theorem 6.3 there says that any $\bigwedge$-definable group $H$ is an intersection of definable groups $H_i$.
If $H$ has bounded index, then each $H_i$ has bounded index and hence by compactness, finite index; but then $H_i=G$ for each $i$ so $H=G$.  \eprf

\begin{rem}  \rm
Unlike \corref{convolution-constant-quasirandom}, where arbitrary bounded functions $f_n,g_n$ are allowed, in \corref{convolution-sc}  it is essential that they be uniformly definable over finite fields.

For instance on $\Gg_a$, if  $f_p(x \mod p)=10$ for $0\leq x < p/10$, and $f_p(u)=0$ for all other $u \in \Ff_p$, we see that the ultraproduct $f$ requires at least 10 self-convolutions to become uniform, and not $2$.  
This accounts for the additional model theoretic ingredient (8.5,6.3) quoted above; the 
homomorphism $n \mod p \mapsto exp(2\pi i /p)$ exists, and must be shown not to be uniformly definable.   
\end{rem}

To connect to convolution of pushforward measures, we will need one simple geometric lemma:

\begin{lem}  \label{dominant-fibers}  Let $f: Y \to X $ be a dominant morphism of  irreducible varieties over  an ultrapower $F= \lim_u \Ff_q$ of finite fields $\Ff_q$; let $d$ be the generic fiber dimension, and $d_X = \dim(X)$.   Let $F(x) =  q^{-d}  |f \inv(x)| $.   
Then there exists a  definable  (in $Th(F)$) real-valued function $f$ on $X$
and a proper subvariety $X_0$ of $X$, such that  $f - F$ tends to $0$ uniformly  
along $u$ on $X \m X_0$.   Also $||F-f||_1$ tends to $0$.
\end{lem}
\prf  Up to the last sentence, this is the main result of \cite{cdm}.  
We may take $f$ to vanish on $X_0$.   
The last sentence
follows since $Y_0:=f \inv(X_0)$ is a proper subvariety of $Y$ and hence has dimension $< \dim(Y)$.  So $q^{-\dim(X)} \sum_{x \in X_0} F(x)  = q^{\dim(Y)} | Y_0| =O(q^{-1/2})$.    
\eprf

      \begin{cor} \label{sc-measures} Let $G$ be a simply connected algebraic group defined over
      $\Zz[m \inv] $,  and let $\eta_i: Y_i \to G$ be a dominant morphism of  irreducible varieties ($i=1,2$). 
      Let $\nu^i_p$ be the normalized counting measure on $Y_i(\Ff_p)$ and let $\mu^i_p = (\eta_i)_* \nu_p^i$ be the pushforward of $\nu^i_p$ to $G$.   Also for a prime $p>m$ let $\mu_p$ be Haar on $G_p(\Ff_p)$.
      Then $||\mu_p - \mu^1_p * \mu^2_p ||_1 \to 0$ as $p \to \infty$.
     \end{cor}

\prf   We may write $\mu_p^i = F_i \mu_p$.  Let $f_i$ be as given in 
\lemref{dominant-fibers}; so $||\mu_p^i - f_i \mu_p||_1 \to 0$.  
Now the statement follows using continuity of convolution on the $L^1$-norm.   
\eprf 

To compare this to  \cite{LST}, Theorem 1,   set 
  $Y_i = G^{d_i}$ where $w_i = w_i(x_1,\ldots,x_{d_i})$, and let $f_i$ be the word map;
  then $f_1 * f_2$ is the word map associated to $w_1w_2$, since they have disjoint variables.

\end{section}


\begin{thebibliography}{10}
 

 
 
 
\bibitem{avnigarion} Nir Avni, Shelly Garion, connectivity of the product replacement graph of simple groups of bounded lie rank,
Journal of Algebra 320(2), 945--960,  2008


\bibitem{bkt}    Bakker, B, B Klingler, and J Tsimerman. “Tame Topology of Arithmetic Quotients and Algebraicity of Hodge Loci.” Journal of the American Mathematical Society. 33.4 (2020): 917--939.  


\bibitem{bft}  
Benjamini, Itai; Finucane, Hilary; Tessera, Romain,  On the scaling limit of finite vertex transitive graphs with large diameter. Combinatorica 37 (2017), no. 3, pp. 333--374;   arXiv:1203.5624.

\bibitem{benjamini-hutchcroft}  Itai Benjamini, Tom Hutchcroft, 
Large, lengthy graphs look locally like lines.  Bulletin of the LMS, Nov. 2020,  \url{https://doi.org/10.1112/blms.12436}

\bibitem{benjamini-schramm} 
Itai Benjamini, Oded Schramm,  Recurrence of Distributional Limits of Finite Planar Graphs,  Electronic Journal of Probability, Electron. J. Probab. 6, 1--13, (2001)
 
\bibitem{benyaacov-c-nip}   Ita\"i Ben Yaacov. Continuous and random Vapnik-Chervonenkis classes. Israel Journal of Mathematics, 2009, 173, pp.309--333
 
\bibitem{benyaacov-g}  Ita\"i Ben Yaacov, Model theoretic stability and definability of types, after A. Grothendieck, Bull. Symb. Log. 20 (2014), no. 4, 491--496, Bulletin of Symbolic Logic



\bibitem{BJM} Samuel Braunfeld, Colin Jahel, and Paolo Marimon, 
 When invariance implies exchangeability
(and applications to invariant Keisler measures), forthcoming.

\bibitem{cat}  Ita\"i Ben Yaacov,  Schroedinger’s cat, Israel Journal of Mathematics 153 (2006), 157-191

\bibitem{byu}  Ita\"i Ben Yaacov,  Alexander Usvyatsov,
Continuous first order logic and local stability  
 Trans. Amer. Math. Soc. 362 (2010), 5213-5259. 
 
 
 \bibitem{benynip}
 Ben Yaacov, Ita\"i,
Continuous and random Vapnik-Chervonenkis classes.  
Israel J. Math. 173 (2009), 309-333
  \bibitem{bybhu}  Itai Ben Yaacov, Alexander Berenstein, C. Ward Henson, and Alexander Usvyatsov, 
 Model Theory for Metric Structures .    In Model Theory with Applications to Algebra and Analysis, Vol. II, eds. Z. Chatzidakis, D. Macpherson, A. Pillay, and A.Wilkie, Lecture Notes series of the London Mathematical Society, No. 350, Cambridge University Press, 2008, 315--427
 
\bibitem{byb}   Ben Yaacov, Itay; Berenstein, Alexander Imaginaries in Hilbert spaces. Arch. Math. Logic 43 (2004), no. 4, 459--466.
 
 
 
 
 
 \bibitem{benykeisler} 
 Ben Yaacov, Ita\"i; Keisler, H. Jerome,
Randomizations of models as metric structures. 
Confluentes Math. 1 (2009), no. 2, 197--223. 
 
 \bibitem{benyaacov-unbounded}   Ben Yaacov, Ita\"i,  
continuous first order logic for unbounded metric structures, 
 Journal of Mathematical Logic Vol. 08, No. 02, pp. 197--223 (2008)  



\bibitem{bourbaki}  N. Bourbaki, \'El\'ements de Math\'ematique, Vari\'et\'es diff\'erentielles et analytiques, 
Fascicule de r\'esultats, Springer 1981.

\bibitem{BGS}  Bourgain, Gamburd and Sarnak,
Markoff triples and strong approximation,  C. R. Acad. Sci. Paris, Ser. I 354 (2016) 131--135.
 
\bibitem{BGS1}  Jean Bourgain, Alexander Gamburd, Peter Sarnak,  
Affine linear sieve, expanders, and sum-product,  
 Invent math (2010) 179: 559--644  


  \bibitem{bgt}  Emmanuel Breuillard, Ben Green, Terence Tao, The structure of approximate groups.   arXiv:1110.5008,
Publ. Math. Inst. Hautes Études Sci. 116 (2012), 115--221.

\bibitem{breuillard-lectures}  Breuillard, Emmanuel, 
Lectures on approximate groups and Hilbert's 5th problem,
IMA Vol. Math. Appl., 159
Springer, [Cham], 2016, 369--404.

 


\bibitem{cantat-loray}  S. Cantat, F. Loray, Dynamics on character varieties and Malgrange irreducibility of Painlev\'e VI equation, Ann. Inst. Fourier (Grenoble) 59 (2009), pp. 2927--2978.



\bibitem{clpz}  Casanovas, E.; Lascar, D.; Pillay, A.; Ziegler, M. Galois groups of first order theories. J. Math. Log. 1 (2001), no. 2, 305-319.  
 

 
 
 
 \bibitem{chang-keisler-c}   Chang, Chen-chung; Keisler, H. Jerome Continuous model theory. Annals of Mathematics Studies, No. 58 Princeton Univ. Press, Princeton, N.J. 1966 xii+166 pp.
 
   
\bibitem{cdm}  Z. Chatzidakis, L. van den Dries and A. Macintyre, Definable sets over finite fields, Journal f\"ur die reine und angewandte Mathematik 427 (1992), 107-135.


   \bibitem{ChH} Chatzidakis, Z.,
Hrushovski, E., Model Theory of difference fields, AMS Transactions v.
   351 (1999), No. 8, pp. 2997--3071
   
     
 \bibitem{CH}   Cherlin, Gregory; Hrushovski, Ehud Finite structures with few types. Annals of Mathematics Studies, 152. Princeton University Press, Princeton, NJ, 2003. vi+193


\bibitem{chernikov-simon} A. Chernikov, P. Simon:  Definably amenable NIP groups, JAMS  31, Number 3, July 2018, Pages 609--641 

\bibitem{chernikov-starchenko} 
Artem Chernikov, Sergei Starchenko,  Definable regularity lemmas for NIP hypergraphs, arXiv:1607.07701 

\bibitem{chernikov-towsner} Artem Chernikov, Henry Towsner,  Hypergraph regularity and higher arity VC-dimension,  arXiv:2010.00726,

\bibitem{elad-alexis}  Alexis Chevalier, Elad Levi,  
An Algebraic Hypergraph Regularity Lemma, 	arXiv:2204.01158.

\bibitem{chevalier-H}  Alexis Chevalier, E. Hrushovski, Piecewise Interpretable Hilbert Spaces, arXiv:2110.05142 [math.LO]

 
 
 


\bibitem{conant-h-p}  Compactifications of pseudofinite and pseudo-amenable groups,
Gabriel Conant, Ehud Hrushovski, Anand Pilla,  arXiv:2308.08440.



 \bibitem{vddries} L. van den Dries, \textit{Dimension of definable sets, algebraic boundedness and henselian fields}, Ann. Pure Appl. 
Logic \textbf{45} (1989), 189--209.


\bibitem{bert} {BERT}: Devlin, Jacob and
      Chang, Ming-Wei  and
      Lee, Kenton  and
      Toutanova, Kristina,  Pre-training of Deep Bidirectional Transform,  Proceedings of the 2019 Conference of the North {A}merican Chapter of the Association for Computational Linguistics: Human Language Technologies, Volume 1,  June 2019;  
arxiv.org/pdf/1810.04805

 \bibitem{GH}
Glazer, I. and Hendel, Y.I. (2021), On singularity properties of convolutions of algebraic morphisms - the general case. J. London Math. Soc., 103: 1453--1479. https://doi.org/10.1112/jlms.12414

\bibitem{gt}  An approximate logic for measures (with Henry Towsner), Israel Journal of Mathematics, Volume 199, Number 2 (2014), 867--913.

\bibitem{gowers-long}  Gowers, W. T. ; Long, J. 
Partial associativity and rough approximate groups.  Geom. Funct. Anal. 30 (2020), no. 6, 1583--1647.

\bibitem{gromov} Misha Gromov, {\em  Metric Structures for Riemannian and Non-Riemannian Spaces}. Birkhuaser, 1998.

\bibitem{gromov-ergo}  Misha Gromov,  Structures, Learning and Ergosystems: Chapters 1-4, 6.
 December 30, 2011.   See \url{https://www.ihes.fr/~gromov/category/ergosystems/}
 
 \bibitem{gromov-email} e-mails to the author, 9/15/10 and 9/17/10.
   
 
 
 
 
 
 
 
 
 
 
 
 
 
 
 
 
 
 

\bibitem{grothendieck52} Alexandre Grothendieck, Criteres de compacit\'e dans les espaces fonctionnels  generaux.  
American Journal of Mathematics
Vol. 74, No. 1 (Jan., 1952), pp. 168-186

\bibitem{HKP}

Hart, Bradd; Kim, Byunghan; Pillay, Anand,  Coordinatisation and canonical bases in simple theories
J. Symbolic Logic 65 (2000), no. 1, 293--309.

\bibitem{haussler-welzl}   Haussler, David; Welzl, Emo,  
$\epsilon$-nets and simplex range queries. 
Discrete Comput. Geom. 2 (1987), no. 2, 127-151.



\bibitem{H-PAC}   Hrushovski, E., Pseudo-finite fields and related structures, 
 Model theory and applications,  151--212, Quaderni Mat., 11, Aracne, Rome, 2002.
 
 

\bibitem{simple} Hrushovski,   Simplicity and the Lascar group, preprint 

\bibitem{cigha} 
Ehud Hrushovski,  Groupoids, imaginaries and internal covers, math.LO/0603413 , 
Turkish Journal of Mathematics , Volume 36, Issue 2, (2012) pp. 173--198.    
 

 \bibitem{nqf} Hrushovski, Ehud, Stable group theory and approximate subgroups. J. Amer. Math. Soc. 25 (2012), no. 1, pp. 189--243
 
 
\bibitem{ax+}  Ax's theorem with an additive character.    EMS Surveys in Mathematical Sciences
Volume 8, Number 1/2, 2021, pp. 179-216.  { \tt{   https://ems.press/journals/emss/issues/1652}} 
    


 
 
 

\bibitem{NIP1}  Ehud Hrushovski, Anand Pillay.  On NIP and invariant measures.    JEMS  13 (4), June 2011.  arXiv:0710.2330 math.LO.

 
    
\bibitem{HP2}
Hrushovski, E.  \&
Pillay, A., 
Definable subgroups of algebraic groups over finite fields.
J. Reine Angew. Math. 462 (1995), 69--91.

\bibitem{H-beyo} Ehud Hrushovski, Beyond the Lascar Group, ArXiv  [2011.12009] 
 
\bibitem{pac} Hrushovski, Ehud, 
Pseudo-finite fields and related structures. Model theory and applications, 151-212, 
Quad. Mat., 11, Aracne, Rome, 2002. 

\bibitem{HPP}
Hrushovski, Ehud; Peterzil, Ya'acov; Pillay, Anand Groups, measures, and the NIP.  J. Amer. Math. Soc.  21  (2008),  no. 2, 563--596.


\bibitem{hps}  Ehud Hrushovski, Anand Pillay, Pierre Simon, Generically stable and smooth measures in NIP theories,   arXiv:1009.0252,
    Trans. Amer. Math. Soc. 365 (2013); no. 5, 2341--2366.
    
    
\bibitem{HW} Hrushovski, E. Wagner, F., Counting and dimensions.  
  Model theory with applications to algebra and analysis. Vol. 2,  161--176,
London Math. Soc. Lecture Note Ser., 350, Cambridge Univ. Press, Cambridge, 2008. 

 
 
\bibitem{ibarlucia}  
Tomas Ibarlucia,   Infinite-dimensional polish groups and property (T), Inventiones mathematicae, 223(2):725--757, 2021.

\bibitem{jahel}  Jahel, Colin,   Des progr\'es sur le probl\'eme d’unique   
ergodicit\'e, Thèse de doctorat, Institut Camille Jordan, \url{  https://tu-dresden.de/mn/math/algebra/das-institut/beschaeftigte/colin-jahel/ressourcen/dateien/These.pdf?lang=en }

\bibitem{JT}  Jahel, Colin and Tsankov, Todor,
Invariant measures on products and on the space of linear orders, J. \'Ec. polytech. Math.9(2022), 155--176.
 
  \bibitem{keisler77}  Keisler,  H. Jerome, Hyperfinite model theory, Logic Colloquium 76, 1977, pp. 5-110.
 
 
 \bibitem{keisler} Keisler,  H. Jerome, Probability Quantifiers. Model Theoretic Languages, Springer-Verlag, pages 509-556 in Model Theoretic Logics, edited by J. Barwise and S. Feferman, 1985.
 

 
 
  
\bibitem{kim-pillay}  Kim, Byunghan; Pillay, Anand Simple theories. Joint AILA-KGS Model Theory Meeting (Florence, 1995).  Ann. Pure Appl. Logic  88  (1997),  no. 2-3, 149--164.

\bibitem{KP2}   Kim, Byunghan; Pillay, Anand,  From stability to simplicity,
Bull. Symbolic Logic 4 (1998), no. 1, 17--36.

\bibitem{kobayashi}  Shoshichi Kobayachi, Transformation Groups in Differential Geometry, Springer 1972.   
 
 \bibitem{krivinemaurey}  J.L. Krivine, B. Maurey, Espaces de Banach stables,Israel J. Math., 39, 4, p. 273-295 (1981).

 \bibitem{kolmogorov}  A.N. Kolmogorov, Foundations of the theory of probability,  tr.
 Nathan Morrison, New York, Chelsea 1956 

\bibitem{LST}    M. Larsen, A. Shalev, and P. H. Tiep, {\em Probabilistic Waring problem for finite simple groups},  Ann. of Math. (2), pp. 561-608, vol. 190(2019), issue 2.   
 

\bibitem{lovasz}  L. Lov\`asz, Large networks and graph limits,  AMS Colloquium Publicatins, vol. 60, 2012.

\bibitem{lovasz-szegedy} L. Lovasz and B. Szegedy: Regularity partitions and the topology of graphons, An Irregular Mind, Szemer\'edi is 70, J. Bolyai Math. Soc and Springer-Verlag (2010), 415-446

\bibitem{lovasz-szegedy2015}
Lov\'asz, L.; Szegedy, B. The automorphism group of a graphon. J. Algebra 421 (2015), 136--166

\bibitem{macpherson-tent} Dugald Macpherson, Katrin Tent, 
Profinite groups with NIP theory and p-adic analytic groups.
 Bulletin of London Math. Soc. 48 (2016) no. 6, 1037--1049.
 
 \bibitem{marimon}   Marimon, Paolo,   Measures and amalgamation properties in $\omega$-categorical structures
doctoral dissertation, section 7.4.3 available in:  \url{  https://spiral.imperial.ac.uk/handle/10044/1/106470 }  

 
 \bibitem{marimon1}  
 Marimon, Paolo. On the non-measurability of $\aleph_0$-categorical Hrushovski constructions. arXiv  2208.06323 (2022), 
 
 \bibitem{marimon2}  Marimon, Paolo. Invariant Keisler measures for omega-categorical structures. arXiv:2211.14628 (2022)

 
 

\bibitem{mvw}  Matthews, C. R.; Vaserstein, L. N., 
Weisfeiler, B., Congruence properties of Zariski-dense subgroups. I.    London Math. Soc. (3) 48 (1984),
no. 3, 514--532
  


 
 

\bibitem{nori} Nori, Madhav V. 
On subgroups of GLn(Fp).
Invent. Math. 88 (1987), no. 2, 257--275. 

\bibitem{pillay-conj} A. Pillay, Type-definability, compact Lie groups, and o-minimality, J. Math. Logic, 4 (2004), 147--162.

\bibitem{razborov}
 A. Razborov: Flag algebras, J. Symbolic Logic 72 (2007), 1239-1282.

  \bibitem{sanders} Tom Sanders, On a non-abelian Balog-Szemeredi-type lemma, arXiv:0912.0306

 \bibitem{lazy}  Saharon Shelah, 
The lazy model-theoretician's guide to stability. 
Comptes Rendus de la Semaine d'\'Etude en Th\'eorie des Mod\`eles (Inst. Math., Univ. Catholique Louvain, Louvain-la-Neuve, 1975). 
Logique et Analyse (N.S.) 18 (1975), no. 71-72, 241-308. 
      
      \bibitem{tao} T. Tao. Hilbert’s fifth problem and related topics. Vol. 153. American Mathematical Soc., 2014.
      
 \bibitem{towsner} Henry Towsner,
A Model Theoretic Proof of Szemeredi's Theorem,arXiv:1002.4456  
 
      \bibitem{vc}  V. N. Vapnik and A. Ya. Chervonenkis,  On the uniform convergence of relative frequencies of events to their probabilities, 
   Theory of probability and its applications,  volume xvi,  number 2, 1971, pp. 264-280.
(Translated by B. Seckler)

 \bibitem{vershik}  A.M. Vershik,  Random and universal metric spaces, 2002.
 
  
\bibitem{yamabe} 
Yamabe, Hidehiko,  On the conjecture of Iwasawa and Gleason. Ann. of Math. (2), 58:48--54, 1953 
\end{thebibliography}
   \end{document}